\pdfoutput=1
\documentclass[11pt]{article}
\usepackage{amsmath,amscd,amssymb,amsthm}
\usepackage[margin=1in]{geometry}
\usepackage{graphicx}
\usepackage{fourier}
\usepackage[sort&compress, comma, square, numbers]{natbib}
\usepackage[small,bf]{caption}
\usepackage{enumerate}
\usepackage[usenames, dvipsnames]{color}
\usepackage[all,import]{xy}
\usepackage{titlesec}
\usepackage{multirow}
\usepackage{authblk}
\usepackage{algorithmic}
\usepackage{algorithm}
\usepackage{grffile}
\usepackage{subfig}
\usepackage{morefloats}
\usepackage{float}
\usepackage[colorlinks,linkcolor=black,bookmarksopen,
bookmarksnumbered,citecolor=black,urlcolor=black]{hyperref}

\graphicspath{{figs/ranking/}}

\def\R{\mathbb{R}}

\def\Z{\mathbb{Z}}

\DeclareMathOperator{\im}{im}

\DeclareMathOperator{\grad}{grad}
\DeclareMathOperator{\curl}{curl}
\DeclareMathOperator{\boldcurl}{\bf curl}
\def\div{\operatorname{div}}

\DeclareMathOperator{\laplacian}{\Delta}

\DeclareMathOperator{\boundary}{\partial}
\DeclareMathOperator{\coboundary}{\delta}

\newcommand{\norm}[1]{\lVert#1\rVert}

\newcommand{\eval}[2]{\langle #1,#2 \rangle}

\newcommand{\innerproduct}[2]{\langle #1, #2 \rangle}

\newcommand{\ainnerproduct}[2]{\langle #1, #2 \rangle}
\newcommand{\aInnerproduct}[2]{\bigl\langle #1, #2 \bigr\rangle}

\long\def\symbolfootnote[#1]#2{\begingroup%
\def\thefootnote{\fnsymbol{footnote}}\footnote[#1]{#2}\endgroup}

\newtheorem{theorem}{Theorem}[section]
\newtheorem*{theorem*}{Theorem}

\newtheorem*{proposition*}{Proposition}
\newtheorem{lemma}[theorem]{Lemma}
\newtheorem*{lemma*}{Lemma}

\newtheorem*{claim*}{Claim}

\newtheorem*{axiom*}{Axiom}

\newtheorem*{conjecture*}{Conjecture}

\newtheorem*{corollary*}{Corollary}

\theoremstyle{definition}

\newtheorem*{definition*}{Definition}
\newtheorem{example}[theorem]{Example}
\newtheorem*{example*}{Example}

\newtheorem*{exercise*}{Exercise}
\newtheorem*{recall*}{Recall}

\theoremstyle{remark}

\newtheorem*{note*}{Note}
\newtheorem{remark}[theorem]{Remark}
\newtheorem*{remark*}{Remark}

\newtheorem*{notation*}{Notation}

\newtheorem*{question*}{Question}
\newtheorem{fact}[theorem]{Fact}
\newtheorem*{fact*}{Fact}

\theoremstyle{theorem}

\theoremstyle{definition}

\theoremstyle{remark}

\newcommand{\lambdamin}{\lambda_{\text{min}}}
\newcommand{\lambdamax}{\lambda_{\text{max}}}

\newcommand{\ER}{Erd\H{o}s-R\'{e}nyi}
\newcommand{\BA}{Barab\'asi-Albert}
\newcommand{\WS}{Watts-Strogatz}

\def\imagetop#1{\vtop{\null\hbox{#1}}}

\newenvironment{xyoverpic}[3]
{%
\begin{xy}
\xyimport#1{\includegraphics[#2]{#3}}
}{\end{xy}}

\allowdisplaybreaks

\begin{document}

\title{\vspace*{-1.2cm} Least Squares Ranking on Graphs\footnote{A
    preliminary version of this paper was originally posted on arXiv
    as ``Least Squares Ranking on Graphs, Hodge Laplacians, Time
    Optimality, and Iterative Methods''. The present version is much
    expanded in scope and includes many new numerical experiments.}}

\author{Anil N.~Hirani\thanks{Author for correspondence:
    \href{mailto:hirani@cs.illinois.edu}{hirani@cs.illinois.edu};
    \href{http://www.cs.illinois.edu/hirani}
    {http://www.cs.illinois.edu/hirani}}}
\author{Kaushik Kalyanaraman}
\affil{Department of Computer
  Science, University of Illinois at Urbana-Champaign}
\author{Seth Watts}
\affil{Department of Mech. Sci. \& Eng.,
University of Illinois at Urbana-Champaign}

\date{}
\maketitle

\begin{abstract}

  Given a set of alternatives to be ranked, and some pairwise
  comparison data, ranking is a least squares computation on a
  graph. The vertices are the alternatives, and the edge values
  comprise the comparison data. The basic idea is very simple and old
  -- come up with values on vertices such that their differences match
  the given edge data. Since an exact match will usually be
  impossible, one settles for matching in a least squares sense. This
  formulation was first described by Leake in 1976 for ranking
  football teams and appears as an example in Professor Gilbert
  Strang's classic linear algebra textbook. If one is willing to look
  into the residual a little further, then the problem really comes
  alive, as shown effectively by the remarkable recent paper of Jiang
  et al. With or without this twist, the humble least squares problem
  on graphs has far-reaching connections with many current areas of
  research. These connections are to \emph{theoretical computer
    science} (spectral graph theory, and multilevel methods for graph
  Laplacian systems); \emph{numerical analysis} (algebraic multigrid,
  and finite element exterior calculus); other \emph{mathematics}
  (Hodge decomposition, and random clique complexes); and
  \emph{applications} (arbitrage, and ranking of sports teams). Not
  all of these connections are explored in this paper, but many
  are. The underlying ideas are easy to explain, requiring only the
  four fundamental subspaces from elementary linear algebra. One of
  our aims is to explain these basic ideas and connections, to get
  researchers in many fields interested in this topic. Another aim is
  to use our numerical experiments for guidance on selecting methods
  and exposing the need for further development.  Many classic Krylov
  iterative methods worked well for small to moderate-sized problems,
  with trade-offs described in the paper. Algebraic multigrid on the
  other hand was generally not competitive on these graph problems,
  even without counting the setup costs.


  \smallskip
  \noindent{\bf MSC Classes: } 65F10, 65F20, 58A14, 68T05, 05C50; {\bf
    ACM Classes: } G.1.3, F.2.1, G.2.2

\end{abstract}

\section{Introduction} \label{sec:intro} 

This paper is about ranking of items, of which some pairs have been
compared. The formulation we use (and which we did not invent) leads
to a least squares computation on graphs, and a deeper analysis
requires a second least squares solution.  The topology of the graph
plays a role, in a way that will be made precise later. The usual
graph Laplacian plays a central role in the first least squares
problem. But the key actor in the second problem is another
Laplacian, hardly studied in theoretical computer science, but
well-studied in numerical analysis. 

The formulation as two least squares problems is akin to finding the
gradient part of a vector field and then finding the curl part. That
in turn, is related to solving an elliptic partial differential
equation. The setting for the ranking problem however, is obviously
different from that for vector fields and differential
equations. Instead of domains that are approximated by meshes, one has
\emph{general} graphs as a starting point. We try to convey the
conceptual and algorithmic implications of these connections and
differences. Another, more practical motivation, is to give guidance
on which numerical methods to use for the ranking problem. As a side
benefit, we are able to point out some directions in numerical linear
algebra and other areas that should be developed further.

\subsection{Ranking and pairwise comparisons}

Even without the many connections brought out in this paper, ranking
is an important problem. Human society seems to have a preoccupation
with ranking. Almost nothing is immune from society's urge to rank and
use rankings.  Rankings are used for marketing, for decision-making,
for allocation of resources, or just for boasting or urging
improvements. Some rankings are based on opinions of a person. But
some are based on opinion polls or some other type of numerical data,
allowing the use of simple statistics and computational tools. Ranking
based on data often consists of computing some weighted average
quantity for each item to be ranked. Then sorting on this number
yields the ranking. However, the ranking problem has more interesting
formulations and algorithms when comparisons and links between the
objects are part of the data.  A well-known example is the PageRank
algorithm used by Google for ranking webpages, which uses the linked
structure of the web (in addition to auctioned keywords) to present
webpages ranked by ``importance''~\cite{PaBrMoWi1999}.

Our focus is on \emph{pairwise comparisons}. As such, this is a very
different starting point than PageRank and very different from ranking
each item independently.  Our interest in pairwise comparisons was
initiated by the recent remarkable paper of Jiang et
al.~\cite{JiLiYaYe2011}. By using exterior calculus concepts and Hodge
decomposition, they motivated the ranking problem in a very different
context and manner than what had been done before. (The less common of
these terms will be described later. One of the pedagogical advantages
of the ranking problem is the ease with which this can be done using
graphs.)

The exterior calculus approach of Jiang et al.~resonated very deeply
with us. We have been working with discretizations of exterior
calculus for a long time, albeit never with a mind towards
ranking~\cite{Hirani2003}. Exterior calculus is the language of modern
physics~\cite{AbMaRa1988, Frankel2004} and can be thought of as the
generalization of vector calculus to smooth curved domains
(manifolds). The discretizations of this calculus have made inroads
into computer graphics~\cite{DeKaTo2008}, and numerical partial
differential equations~\cite{ArFaWi2010} with great effect. Thanks to
the work of Jiang et al.~we hope for a similar development in graph
computations, especially in the case of ranking.  This is the
development that we are trying to urge forward in this paper. 

\subsection{Contributions and goals of this paper}

In this paper you will find the important concept of Hodge
decomposition reduced to the basics of linear algebra, beyond which it
probably cannot further be simplified. You will also find many
numerical linear algebra approaches for least squares ranking on
graphs, suitable for serial computer implementation. There are many
numerical experiments, some with expected results, and some with
suggestive results, but most have never appeared in literature before.
Many of the results of these experiments call out for further
investigations in a variety of fields which we point out along the
way. The subject matter is new and timely, and attention needs to be
drawn to it, even if there are no deep theorems yet in the numerical
analysis aspects of the problem. We hope that the reader will
appreciate the novelty of the experiments and their broad implications
that we point out.  The simplicity of the exposition is aimed at wide
dissemination of the ideas. Our main aim is to generate enthusiasm and
further work in the computational mathematics community, in this broad
arena for numerical problems on graphs.

\medskip\noindent{\bf Use of orientation instead of skew symmetry: }
As a point of departure from previous work, we have a very different
implementation than what is implied by Jiang et al. This difference
stems from our experience in discrete and finite element exterior
calculus. Jiang et al. rely on skew-symmetric tensors as the basic
objects of their formulation, whereas we capture the skew-symmetry
using a vector of values associated with oriented objects. This 
results in savings in space and time, and provides a simpler 
formulation which is easy to implement. This is an example of the 
synergy of studying discrete calculus on meshes and graphs in a common
framework.  This synergy is also on display in a recent book by Grady
and Polimeni~\cite{GrPo2010}. Their chapter on ranking does reference
Jiang et al.~but does not discuss the ideas of that paper. In contrast
to their book, our paper is concerned solely with the problem of least
squares ranking on graphs.

\medskip\noindent{\bf Guidance on solvers: } Our focus on a single
fundamental problem allows us to probe the numerical linear algebra
aspects much more deeply. We compare iterative Krylov
methods~\cite{Saad2003, Vandervorst2009} and algebraic
multigrid~\cite{BrMcRu1985,RuSt1987}, and point out some problems ripe
for further development. Specifically, algebraic multigrid is not
competitive for these graph problems.  The recent multilevel solver of
Koutis el al.~\cite{KoMiPe2010} should also be useful for the ranking
problem. We were not able to experiment with it. We could not find a
reliable complete implementation at the time of writing this paper and
could not get our implementation to work. The authors of that solver
should be able to modify our code easily, to compare the performance
of their solver with Krylov and multigrid methods. The solver of
Koutis et al.~is designed for graph Laplacians and so it should be
quite well suited for the first least squares problem, as we will
reason later. We will see that the second least squares problem
involves a different type of Laplacian and the solver of Koutis et
al.~will likely \emph{not} have any special advantage when used for
such problems.

\medskip\noindent{\bf Proposal for a spectral simplicial theory:
} As mentioned earlier, the least squares problems for ranking on
graphs involve two types of Laplacians. The eigenvalues and
eigenvectors of these operators are worthy of study. Laplacians on
manifolds have been studied extensively in
geometry~\cite{Chavel1984}. However, the study of Laplacian spectra on
graphs is more recent. The first least squares problem involves the
graph Laplacian, which is a very well-studied operator. It is the
spectrum of this graph Laplacian (and closely related operators) that
is studied in the very useful and beautiful spectral \emph{graph}
theory~\cite{Chung1997}. We are proposing that a similar spectral
\emph{simplicial} theory should be developed for simplicial
complexes. An example application is that the matrix spectrum is
useful for understanding the performance of iterative solvers. We hope
however, that the payoff from a spectral simplicial theory will go
well beyond that. A glimpse of this potential is in the topic of the
next paragraph. In that application the eigenvectors of the zero
eigenvalue play a role as will be described in detail later.

\medskip\noindent{\bf Numerical experiments on topology of clique
  complexes: } One can use the numerical linear algebra approach
embodied in our paper to conduct experiments in the new field that
studies the topology of random clique complexes~\cite{Kahle2009}. A
clique complex is a graph augmented by its complete subgraphs which
are considered as additional structures on the graph. This field
generalizes the established field of random graph
theory~\cite{ErRe1960, Bollobas2001}. The concepts of connectivity of
a graph generalize naturally to questions about the homology of clique
complexes. We include several examples in which numerical techniques
for ranking on graphs lead to experiments on clique complexes.  In one
of these examples we reproduce some recent theoretical results for
clique complexes of random graphs~\cite{Kahle2009}. In fact, our
experiments on random graphs are also suggestive of new conjectures
and refinements of the existing theorems. We also include some
experiments on the 1-cohomology of scale-free graphs~\cite{BaAl1999},
which is something that has not been analyzed theoretically elsewhere.

\medskip\noindent{\bf Suggestions for Graph 500 benchmarks: } In the
high performance computing community, a move is afoot to establish and
maintain a Graph 500 list. This is like the Top 500 list of 
supercomputers that is regularly updated, but focused on graph 
problems. Benchmark problems are being developed in the process. We 
want to draw attention to least squares ranking on graphs as a source 
of benchmarks. In the field of high performance computing, problems 
like least squares and linear systems for elliptic partial 
differential equations have been an important source of problems. 
These have led to many developments, such as in domain decomposition, 
preconditioners, and iterative methods. The problem of least squares 
ranking on graphs also involves Laplacians, but these are graph 
Laplacians and other Laplacians on \emph{general} graphs. When very 
large ranking problems on diverse architectures are attempted, it is 
likely that new developments will be needed. At the same time, the 
problems are easy to set up and some old codes from differential 
equations can be used right away. Thus the least squares ranking on 
graphs is a good crossover problem and a bridge from Top 500 to Graph
500.

\section{Least Squares and Ranking}\label{sec:lstsqrs_rnkng}

Let us first examine more carefully why two least squares problems
are involved in this formulation of ranking.  Given is a set of items
to be ranked and some real-valued pairwise comparisons. Not all pairs
need to have been compared. Each given pairwise comparison represents
how much one alternative in a pair is preferred over the other.

The data can be represented as a weighted directed graph, where each
comparison between a pair of items is represented by two edges of
opposite direction and equal weight between the items (this is not the
formulation we use). This leads to the skew-symmetric 2-tensor
representation of comparisons used by Jiang et al. However, this is
equivalent to a simple, weighted, undirected graph whose edges are
oriented. But this is exactly an oriented abstract simplicial
1-complex, which we will define precisely in
Section~\ref{sec:preliminaries}. The edge orientations are arbitrary,
and the pairwise score simply changes sign if the opposite edge
orientation is used. Without loss of generality, we will usually only
consider connected graphs. (Multiple component graphs result in
independent ranking problems, one for each component.)

One version of the ranking problem is to find real-valued scores for
the vertices, which implies their global rank order, such that the
values represent the strength of the rank. The task translates to
finding vertex values whose differences are the edge values. However,
it is not always possible to compute this exactly. Every loop in the
graph has the potential to make existence of such vertex values
impossible if the edge values, taken with signs, do not add up to zero
as the loop is traversed. What saves this procedure is that, the
closest possible global ranking is the vertex value assignment whose
differences reproduce the pairwise edge data in the least squares
sense. This is a very simple and an old idea that was used for ranking
football teams by Leake~\cite{Leake1976}. The residual, i.e., the part
of the edge data that could not be matched, represents inconsistencies
in the pairwise data.

\begin{remark} \label{rem:cnsstnt}
  The term \emph{consistent} is used in a technical sense above. The
  edge data is consistent if the sum of edge weights (taking
  orientations into account) around every loop is zero. Consistency is
  equivalent to having a zero residual in the least squares
  sense. Least squares ranking can still be computed for inconsistent
  edge data. For example, the data illustrated in the left graph below
  is consistent while the other two are not. However, a least squares
  ranking is possible on each.  In particular, $A$ will be the winner
  and $C$ the loser in the first two cases, and there will be a
  three-way tie in the last one.
\end{remark}

\[
\xymatrix{ & A\ar[dr]^{-2} & & & A\ar[dr]^{-1} & & & A\ar[dr]^1 & \\
B\ar[ur]^1 & & C\ar[ll]^1 & B\ar[ur]^1 & & C\ar[ll]^1 & 
B\ar[ur]^1 & & C\ar[ll]^1}
\]

A recent extension of Leake's idea was given by Jiang et
al.~\cite{JiLiYaYe2011} who examine the residual, decomposing it into
local and global inconsistencies, using a second least squares
formulation. This time a 3-tensor, i.e., a 3-dimensional matrix, is
involved as an unknown in their formulation. As we will show,
equivalent least squares problems can be formulated using vectors
rather than matrices or 3-tensors to represent the data and unknowns.

We will however follow the point-of-view introduced by Jiang et al.,
who posed the ranking problem as a discrete \emph{Hodge decomposition}
of the pairwise data treated as a \emph{cochain} on a \emph{simplicial
  2-complex}. After a discussion of ranking methods in
Section~\ref{sec:other_ranking}, in Sections~\ref{sec:preliminaries}
and \ref{sec:hodge} we will define the terms emphasized in the
previous sentence. This will help us clarify the connection to vector
field decomposition, elliptic partial differential equations, and
topology of complexes.

We emphasize that Leake's idea was to use least squares to find values
for each vertex in the graph, given (generally inconsistent) edge
data. This is very different from the use of least squares to fit data
to a model equation in the sense of statistical regression. It is
better to think of the methods of Leake and Jiang et al.~as being
projections or decompositions, which is the viewpoint that we will take
in this paper.

\section{Other Ranking Methods} \label{sec:other_ranking}

The problem of ranking has been addressed in many areas such as social
choice theory, sports analysis, ranking of web pages and machine
learning. In this section, we compare and contrast many of these
methods with least squares ranking on graphs.

\subsection{Social choice theory}
In social choice theory, the goal is to rank alternatives based on
preferences of a number of voters. Often, the methods are concerned
with \emph{ordering} alternatives but not quantification of their
relative strengths.  A common assumption is that voters indicate their
preferences for all alternatives~\cite{JiLiYaYe2011}. This results in
a complete graph for each voter with only edge orientations specified.
However, these may not be reasonable for many real data sets such as
internet and e-commerce applications, and ranking of sports teams.

In Condorcet's method~\cite{Borgers2010}, a winning (or highest
ranked) alternative is one which is preferred over every other in
pairwise comparisons. If no such alternative exists, some tie-breaking
rule is required . This is equivalent to selecting a winner based on
edge orientations of the ranking graph alone and not their weights. On
the other hand, least squares ranking uses both the direction and
weight of edges to rank all alternatives.  The Borda count
\cite{Borgers2010} method ranks alternatives based on a weighted vote
scheme. In contrast to least squares ranking, it does not use pairwise
comparisons. The Kemeny rule \cite{YoLe1978} works by testing all $n!$
possible orderings of $n$ alternatives against input edge orientations
and weights to arrive at a ranking.  Consequently, it is
NP-hard~\cite{BaToTr1989} while least squares ranking is polynomial
time and additionally provides relative strengths of alternatives.
Tideman ranking~\cite{Tideman1987} generates a connected, acyclic
graph by preferentially retaining edges with highest weights, which do
not form a cycle, circumventing inconsistency.  In graphs which are
originally connected and acyclic, least squares ranking produces the
same results as Tideman ranking. However, for general graphs, least
squares ranking accounts for inconsistency.  Tideman ranking is also
ambiguous when all edges have equal weights resulting in a ranking
that depends on which edges are removed.

\subsection{Sports analysis}

Ranking of sports teams is another rich source for ranking
methodologies. In addition to satisfying requirements of fairness,
rankings are often used to also predict future outcomes amongst teams,
which requires some measure of ``accuracy'' of the rankings. (A method
which consistently ranks losing teams highly is not accurate.)  Elo
ranking \cite{Elo1978,Glickman1999} treats each alternative as a
random variable with a given mean and variance. When two alternatives
are compared (e.g., two players in a game of chess) their strengths
and variances are updated based on a set of rules. The rules take into
account both the outcome and the strength of the alternatives before
the comparison. A global ranking is then a sorting of their mean
strengths. Thus, the Elo ranking works exclusively via local
interactions whereas least squares ranking always results in a
globally optimal consistent ranking.  Random walk methods
\cite{CaMuPo2007, KvSo2006, BrSo2010} rank alternatives by finding the
steady state of a Markov chain. The Markov chain models voter opinions
on alternatives when pairwise comparisons and rules for changing votes
are given. There may be some connection between these methods and
least squares ranking given the connections between random walks and
solving Laplace's equation \cite{Reynolds1965}.  Keener
\cite{Keener1993} proposed ranking alternatives such that their ranks
were proportional to their strengths, as defined by a preference
matrix which records the results of pairwise comparisons. The rankings
are found as the Perron-Frobenius eigenvector of the preference
matrix. (Subsequent models also defined the ranks as the fixed point
of a nonlinear strength function and using maximum likelihood methods
amongst others.) Least squares ranking finds ranks which optimally
account for the input data, rather than requiring
proportionality. Additionally, it is easier to recompute ranking on a
given graph for new data by solving the least squares equations with a
new right hand side. In contrast, Keener's method requires creating a
new preference matrix and finding its Perron-Frobenius
eigenvector. Finally, for this eigenvector to exist, the preference
matrix must be irreducible which constrains admissible data on the
graph but there is no such requirement for least squares ranking.
Chartier et al. \cite{ChKrLaPe2011} analyze several popular ranking
methods to determine their stability and sensitivity to perturbations
in their input data. The perturbations are taken about a perfect
season, which is equivalent to ranking on a complete graph with
consistent edge data. Of particular interest is the analysis of the
Massey \cite{Massey1997} ranking method, which is the first least
least squares system (\ref{eq:lsomega1}), which shows that the ordinal
ranks provided by least squares ranking are quite stable in the
presence of perturbations in the input data.

\subsection{Other ranking methods}

PageRank \cite{PaBrMoWi1999} and HITS \cite{Kleinberg1999} are popular
methods for ranking web pages based on hyperlinks between them. These
are different from the pairwise ranking formulation we study, where
there is a value given for how much one alternative is preferred over
the other.  Ranking techniques in machine learning, like in social
choice theory, tend to be concerned with ordering alternatives without
exploring their relative strengths.  Also most such algorithms are for
supervised learning (learning using training data).  For example, the
problem of learning to rank on a graph addressed by
Agarwal~\cite{Agarwal2010} is a supervised learning problem. Moreover,
it is formulated as a constrained quadratic programming problem.  The
constraint function in their case is similar to the least squares
ranking objective function.  The matrix completion method of Gleich
and Lim~\cite{GlLi2011} requires a matrix singular value decomposition
at each iteration, which is more expensive than least squares ranking
using sparse iterative methods.

\section{Preliminaries}\label{sec:preliminaries}

We will see that a convenient language for revealing the connection of
ranking with other fields consists of very basic notions of exterior
calculus and cell complexes, which we will quickly recall in this
section. This is also needed to frame the two least squares problems
as Hodge decomposition, which we do in Section~\ref{sec:rnknghdg}. For
the first least squares problem we need only the graph described in
Section~\ref{sec:lstsqrs_rnkng}. But for the second least squares
problem of ranking we will need to use basic ideas about cell
complexes and functions on them.  We review the irreducible minimum of
the basic terminology and concepts that we need from algebraic
topology and exterior calculus in this section. For more details on
algebraic topology see~\cite{Munkres1984} and for exterior calculus
see~\cite{AbMaRa1988}.

\subsection{From graphs to complexes}

An \emph{abstract simplicial complex} $K$ is a collection of finite
non-empty sets called \emph{simplices} such that if a simplex $\sigma$
is in $K$ then so is every non-empty subset of
$\sigma$~\cite{Munkres1984}. The elements of a simplex are called its
\emph{vertices}. A simplex with $p + 1$ vertices is said to have
\emph{dimension} $p$ and is referred to as a $p$-\emph{simplex}. The
dimension of a complex is the dimension of the highest dimensional
simplex in it. We will refer to a $p$-dimensional abstract simplicial
complex as a $p$-complex.  An \emph{orientation} of a simplex is an
equivalence class of permutations (orderings) of its vertices. All
even permutations fall into one class and the odd ones into the other
class. Thus all simplices of dimension 1 or more have two possible
orientations while a vertex has only one orientation.  An
\emph{oriented} simplex is a simplex along with a choice of an
orientation for it. An oriented abstract simplicial complex is one in
which all simplices have been oriented. For the applications
considered in this paper, the orientations are arbitrary.

Let $G$ be an oriented weighted simple graph, i.e., the edges have
been oriented arbitrarily.  Then $G$ is an abstract simplicial
1-complex.  The vertices of $G$ are the $0$-simplices and the edges
are the $1$-simplices. A $p$-\emph{clique} of $G$ is a complete
subgraph with $p$ vertices. The graph $G$ can be augmented by cliques
to make it an abstract simplicial complex of higher dimension. In
particular, augmenting $G$ by including the $p$-cliques, for all $3
\le p \le d+1$ yields a $d$-dimensional simplicial complex which we
will refer to as the $d$-dimensional \emph{clique complex} of $G$.
For the first least squares problem of ranking we only need $G$ to be
a graph. For the second problem we need the 2-dimensional clique
complex in which the 3-cliques (i.e., triangles) have been oriented
arbitrarily.  Thus we will augment $G$ by including the triangles of
the graph and we will refer to this augmented structure also as~$G$.

\begin{remark}
  The 3-cliques are loops of length three. All the results of this paper
  are valid if we include loops of length $\ge 3$ up to some finite
  length. This yields not a 2-dimensional \emph{simplicial}
  complex but a 2-dimensional \emph{cell}
  complex~\cite{Munkres1984}. In the rest of this paper, the reader
  may substitute the word ``cell'' for ``simplex'' or ``simplicial''
  without changing the results.
\end{remark}

\subsection{Chains and cochains}\label{subsec:chnscchns}

Let $C_p(G;\,\R)$ be the space of real-valued functions on the
oriented $p$-simplices of $G$, such that the function changes sign on
a simplex when its orientation is flipped. (In algebraic topology,
usually one starts with integer valued chains
$C_p(G;\,\Z)$~\cite{Munkres1984}.)  These functions are called
real-valued $p$-\emph{chains}, or $p$-\emph{dimen\-sional}
chains. Since they take values in reals, the space of $p$-chains forms
a vector space.  We will use $C_p(G)$ to abbreviate $C_p(G;\,\R)$. The
\emph{elementary chain basis} for $C_p(G)$ consists of functions that
are 1 on a particular $p$-simplex and 0 on the rest. Thus the
dimension $\dim C_p(G)$ is the number of $p$-simplices in $G$, which
we will refer to by the symbol $N_p$.

The numerical data in the ranking problem are best viewed as
real-valued linear \emph{functionals} on the spaces of chains.  The
reason for using the space of functionals rather than the chains
themselves will become clear when we make the analogy with vector
calculus in Section~\ref{subsec:calculus}. The spaces of functionals
are the vector space duals of $C_p(G)$ and are denoted $C^p(G;\,\R)$,
or $C^p(G)$ and called the spaces of $p$-\emph{cochains}. (Note that
cochain spaces have indices on top.)  The \emph{elementary cochain
  basis} for $C^p(G)$ consists of the cochains that are 1 on an
elementary chain and 0 for the other elementary chains, i.e., it is
the basis dual to the elementary chain basis, where the duality is in
the sense of vector space duality.

\subsection{Boundary and coboundary operators}\label{subsec:operators}

At first sight, the extra structure of chains and cochains of the
previous subsection seems like extra baggage in the ranking
problem. However, the chains and cochains come with the boundary and
coboundary operators that we will now recall. These provide the
scaffolding on which the decomposition or projection view of ranking
is built.

In the exterior calculus view of partial differential equations, the
main objects are often differential forms. In finite element and
discrete exterior calculus these are usually discretized as cochains
on simplicial complexes. The boundary and coboundary operators are
used in defining differential operators, and are the building blocks
of higher order operators like Laplacians. A similar situation holds
for graphs treated as abstract simplicial complexes. One difference
from partial differential equations is the absence of any geometric,
i.e. metric, information. The vertices of the graphs in the ranking
problem need not be placed in any particular geometric location. The
metric information in the differential equations case is captured in
the Hodge star operator, which we will not have occasion to use in the
ranking problem on graphs.

The \emph{boundary} operator $\boundary_p : C_p(G) \to C_{p-1}(G)$ is
usually described by first defining it on $p$-simplices and then
extending it to $C_p(G)$. In the elementary chain basis it takes the
form of a matrix with entries that are either 0 or $\pm 1$. For our
purpose, we take the simpler route and define these directly as
matrices. The first of these is simply the vertex-edge adjacency
matrix of graph theory. This has one column for each edge, with a $-1$
for the starting node and 1 for the ending node of that edge. This is
the matrix form of $\boundary_1$ in the elementary chain
basis. Similarly there is an edge-triangle adjacency matrix
$\boundary_2$. Each column in it corresponds to a triangle and there
is a $\pm 1$ for each of the three edges that appear in the
triangle. The entry is a $+1$ if the triangle and edge orientations
match, and a $-1$ if the orientations do not match. Note that in a
general graph an edge may appear in any number of triangles. This is
different from the simplicial approximation of manifolds, i.e. meshes,
that are used in partial differential equations. This will play an
important role in the performance of linear solvers which were 
originally designed for solving partial differential equations on 
meshes.

The vector spaces and the boundary maps are arranged in what is known
as a \emph{chain complex}
\begin{equation} \label{eq:chncmplx}
\xymatrix{ 0 \ar[r] & C_2(G) \ar[r]^<(.3){\boundary_2} & 
  C_1(G) \ar[r]^<(.3){\boundary_1} &
  C_0(G) \ar[r] & 0} 
\end{equation}
where the first and last maps are zero operators.  The most
important fact about the boundary operators is that $\boundary_p \circ
\boundary_{p+1} = 0$. It is the crucial fact needed in the
decomposition described in Section~\ref{sec:hodge}. Analogous to the
chain complex is the \emph{cochain complex} which is arranged in the
reverse order and uses the \emph{coboundary} operator. Since we are
dealing with real-valued chains and cochains, the coboundary operator
$\coboundary_p$ is simply the dual of the boundary operator
$\boundary_{p+1}$ and is defined by requiring
\[
(\coboundary_p \alpha)(c) = \alpha(\boundary_{p+1} c)\, ,
\]
for all $p$-cochains $\alpha$ and $(p+1)$-chains $c$. This is more
suggestive when the evaluation of a cochain on a chain is written as a
pairing. Then the above relation can be written as
$\eval{\coboundary_p\alpha}{c} =
\eval{\alpha}{\boundary_{p+1}c}$. When the elementary cochain basis is
used, the matrix form of $\coboundary_p$ is simply $\boundary_{p+1}^T$
so we will use the transposed boundary matrix notation rather than the
$\coboundary$ notation. We will write the cochain complex as
\begin{equation} \label{eq:cchncmplx}
\xymatrix{ 0 \ar[r] & C^0(G) \ar[r]^<(.3){\boundary_1^T} & 
  C^1(G) \ar[r]^<(.3){\boundary_2^T} &
  C^2(G) \ar[r] & 0} 
\end{equation}
in which clearly $\boundary_2^T \circ \boundary_1^T = 0$ because of
the analogous property of the boundary matrices.

\subsection{Homology and cohomology}\label{subsec:hmlgychmlgy}

A fundamental problem in topology is to determine if two given spaces
are topologically the same (homeomorphic) or different. We will recall
the definitions of homology and cohomology, which are convenient tools
for distinguishing spaces. Two spaces whose homology or cohomology
differs are not homeomorphic. We need these notions in order to
discuss our experiments on the topology of clique complexes in
Section~\ref{sec:clique}. As in the previous section, we will use
real-valued chains and cochains.

The space of $p$-\emph{cycles} is the space $\ker \boundary_p$ (kernel
of $\boundary_p$), which is a subspace (as a vector space) of
$C_p(G)$.  The image of the boundary map coming from the $p+1$
dimension in a diagram like~\eqref{eq:chncmplx} is $\im
\boundary_{p+1}$ (image of $\boundary_{p+1}$), which is also a
subspace of $C_p(G)$. Their quotient $\ker \boundary_p / \im
\boundary_{p+1}$, in the sense of vector spaces is called the
$p$-dimensional \emph{homology} space and denoted $H_p(G)$. Thus
elements of $H_p(G)$ are equivalence classes of cycles. Cycles $b$ and
$c$ are in the same class if $b-c$ is in $\im \boundary_{p+1}$, and
then $b$ and $c$ are said to be \emph{homologous} to each other.

\begin{remark}\label{rem:betti}
  If the values of the chains (which are called coefficients) need to
  be emphasized, one writes $H_p(G;\, \R)$ or $H_p(G;\ \Z)$ for real
  or integer homology, respectively, and so on. Integer homology
  captures more information than real
  homology~\cite{Munkres1984}. However, the real homology does include
  the one piece of information that is useful for interpreting our
  experiments in Section~\ref{sec:clique}. The number called
  \emph{Betti} number $\beta_p$ for $p$-dimension, which is usually
  defined in integer homology, turns out to be the same as the vector
  space dimension of $H_p(G;\, \R)$. This is a consequence of the
  universal coefficient theorem of algebraic topology~\cite[Chapter
  7]{Munkres1984}, or more simply from~\cite[Theorem
  11.4]{Munkres1984}.
\end{remark}

On the cochain side one has the corresponding $p$-dimensional
\emph{cohomology} space $H^p(G)$, defined as the quotient space $\ker
\boundary_{p+1}^T / \im \boundary_p^T$. As vector spaces, $H^p(G)$ and
$H_p(G)$ are isomorphic. This is because their dimensions are the
same, which follows easily from the rank-nullity theorem of linear
algebra, and the basic facts about the four fundamental subspaces.

\subsection{Laplace-deRham operators} \label{subsec:lplcdrhm}

The cochain and chain complexes can be combined, and excursions in that
diagram lead to various Laplacian operators. The combined diagram that
will suffice for this paper is
\begin{equation}\label{eq:chncchn}
  \xymatrix{C^0(G) \ar[r]^{\boundary_1^T} \ar@{<->}[d] & C^1(G) 
  \ar[r]^{\boundary_2^T} \ar@{<->}[d] & C^2(G) \ar@{<->}[d] \\
  C_0(G) & C_1(G) \ar[l]_{\boundary_1} & C_2(G) 
  \ar[l]_{\boundary_2} \\ }
\end{equation}
%
%
where the vertical arrows are vector space duality isomorphisms. The
excursions used to define the new operators start at one of the
cochains and traverse the box or boxes on the left, right, or both
sides.  In differential geometry and Hodge theory the resulting
operators are known as the \emph{Laplace-deRham} operators
\cite{AbMaRa1988} and denoted $\laplacian_p$ if they act on
differential $p$-forms. In numerical analysis, these operators are
increasingly being referred to as the \emph{Hodge
  Laplacians}~\cite{ArFaWi2006, ArFaWi2010}. The corresponding diagram
in differential geometry consists of differential forms at both
levels. The vertical arrows in that case are the Hodge star operators
which contain the metric information about the manifold.

In the present case, by identifying $C^p(G)$ and $C_p(G)$ via vector
space duality isomorphisms, we can abuse notation and use the identity
operator for the vertical arrows. Three different Laplace-deRham
operators $\laplacian_p : C^p(G) \to C^p(G)$ can be defined for graphs
augmented with triangles and all are of interest. These are
\begin{equation}
\laplacian_0 = \boundary_1\,\boundary_1^T \qquad \qquad
\laplacian_1 = \boundary_1^T\,\boundary_1 +
\boundary_2\,\boundary_2^T \label{eq:lplcdrhm} \qquad \qquad
\laplacian_2 = \boundary_2^T\,\boundary_2\, .
\end{equation}
If cliques with more than 3 vertices were also to be included, then
the definition of $\laplacian_2$ would change to
$\boundary_2^T\boundary_2 + \boundary_3\boundary_3^T$ and there would
be a $\laplacian_3$, $\laplacian_4$ and so on. From the definitions of
these operators it is clear that the matrix form of any Laplace-deRham
operator is square and symmetric.

\begin{remark}\label{rem:hrmncbtt}
  When we study topology of clique complexes, the main objective will
  be to measure the 1-dimensional integer homology Betti number
  $\beta_1$. Using basic linear algebra combined with Hodge
  decomposition of $p$-cochains we give here a simple proof that $\dim
  \ker \laplacian_p = \beta_p$. We will show that 
  \[
  \dim \ker \laplacian_p = \dim H_p(G;\,\R)\, ,
  \]
  as vector spaces. Then by Remark~\ref{rem:betti} the desired result
  follows. Recall that we use $N_p$ for the number of $p$-simplices in
  the clique complex of $G$, and this number is the same as $\dim C_p
  = \dim C^p$. By Hodge decomposition of $p$-cochains (see 
  Section~\ref{sec:hodge}) we have
  \begin{align*}
    \dim \ker \laplacian_p &= \dim C^p  - \dim \im \boundary_p^T -
    \dim \im \boundary_{p+1}\\
    &= N_p - \dim \ker \boundary_p^\perp - \dim \im \boundary_{p+1}\\
    &= N_p - (N_p - \dim \ker \boundary_p) - \dim \im \boundary_{p+1}\\
    & = \dim \ker \boundary_p - \dim \im \boundary_{p+1} = \dim
    H_p(G;\,\R) =  \beta_p\, .
  \end{align*}
\end{remark}

\subsection{Vector calculus analogies}\label{subsec:calculus}

The matrix $\boundary_1^T$ is a graph analog of the gradient and
$\boundary_1$ is the analog of negative divergence. Similarly
$\boundary_2$ is the two-dimensional vector curl and $\boundary_2^T$
is the two-dimensional scalar curl \cite{GiRa1986}. The
two-dimensional vector calculus diagram analogous
to~\eqref{eq:chncchn} is
\begin{equation} \label{eq:deRham2}
\xymatrix{ \text{functions} \ar@<0.5ex>[r]^(.45){\grad} & 
  \ar@<0.5ex>[l]^(.55){-\div}
 \text{vector fields} \ar@<0.5ex>[r]^<(.25){\curl} & 
 \ar@<0.5ex>[l]^(.45){\boldcurl} \text{densities} } 
\end{equation}
since divergence is the negative adjoint of gradient, and scalar and
vector curls are adjoints of each other in two dimensions. This is an
example of a de Rham complex~\cite{BoTu1982, ArFaWi2010}. The more
familiar three-dimensional vector calculus has the de Rham complex
given below.
\begin{equation} \label{eq:deRham3}
\xymatrix{ \text{functions} \ar@<0.5ex>[r]^(.45){\grad} & 
  \ar@<0.5ex>[l]^(.55){-\div}
 \text{vector fields} \ar@<0.5ex>[r]^<(.25){\curl} & 
  \ar@<0.5ex>[l]^(.49){\curl}
 \text{vector fields} \ar@<0.5ex>[r]^<(.25){\div} & 
 \ar@<0.5ex>[l]^(.45){-\grad} \text{densities} } 
\end{equation}
Just as $\boundary_2^T\boundary_1^T = 0$ in the cochain complex,
$\curl\circ \grad = 0$ and $\div \circ \curl = 0$ in these diagrams
above.

The 0-Laplacian in~\eqref{eq:lplcdrhm} is the discrete analog of the
usual scalar function Laplacian in vector calculus and is also the
combinatorial graph Laplacian (without normalization,
see~\cite{Chung1997}). For the de~Rham complex~\eqref{eq:deRham2}
$\laplacian_0 = -\div \circ \grad$. The 1-Laplacian
in~\eqref{eq:lplcdrhm} is the discrete analog of the vector Laplacian
$\laplacian_1 = \boldcurl \circ \curl - \grad \circ \div$ in the de
Rham complex~\eqref{eq:deRham2}. If we did not include the triangles
(or cells) and considered $G$ only as a 1-dimensional complex,
$\boundary_2$ would be the zero matrix. Then the 1-Laplacian would be
$\boundary_1^T\,\boundary_1$ which is sometimes called the edge
Laplacian in graph theory. There is no name for $\laplacian_2$ in
graph theory. But this 2-Laplacian is the graph theoretic analog of
the 2-Laplacian in Hodge theory~\cite{AbMaRa1988} and finite element
exterior calculus~\cite{ArFaWi2010} on a 2-dimensional manifold.

\section{Hodge Decomposition of Vector Spaces}\label{sec:hodge}

Hodge decomposition is an important tool in computer graphics
\cite{ToLoHiDe2003}, engineering~\cite{ChMa1993} and
mathematics~\cite{AbMaRa1988, Morita2001}. It generalizes the
well-known Helmholtz decomposition of vector fields in Euclidean space
\cite{ChMa1993} to differential forms on manifolds. The Helmholtz
decomposition states that every vector field on a compact simply
connected domain can be decomposed into a gradient of a scalar
potential and a curl of a vector potential. The decomposition is
orthogonal and hence unique although the potentials are not. The first
part is curl-free and the second part is divergence-free. If the
domain has nontrivial 1-dimensional homology (e.g., if it is an
annulus, or a torus) then a third component called the harmonic vector
field arises.

For finite-dimensional vector spaces, Hodge decomposition is a really
simple idea. It is just the ``four fundamental spaces'' idea
(popularized in Professor Strang's books) taken one step further. For
a matrix $A$ with $m$ rows and $n$ columns, the four fundamental
subspaces are the column space of $A$, the nullspace of $A$, the row
space of $A$ (which is the column space of $A^T$) and the left
nullspace of $A$ (which is the nullspace of $A^T$) \cite[page
90]{Strang1988}. We slightly prefer the terminology used in algebra
where these would be denoted $\im A$, $\ker A$, $\im A^T$, and $\ker
A^T$ and referred to as the image of $A$, kernel of $A$, image of
$A^T$, and kernel of $A^T$.

Let $U$, $V$ and $W$ be finite-dimensional inner product vector
spaces. Let $A : U \to V$ and $B: V \to W$ be linear maps such that $B
\circ A = 0$. Define $\laplacian := A A^T + B^T B$. The vectors in
$\ker \laplacian$ are called \emph{harmonic}. Pictorially, we have
\begin{equation} \label{eq:UVW}
\xymatrix{ U \ar@<0.5ex>[r]^(.45){A} & \ar@<0.5ex>[l]^(.5){A^T}
 V \ar@<0.5ex>[r]^<(.35){B} & \ar@<0.5ex>[l]^(.5){B^T} W } 
\end{equation}

More formally, the transposes are the adjoint operators which are maps
between the vector space duals, and adjointness requires the presence
of inner products. For example, $A^T : V^\ast \to U^\ast$, where
$V^\ast$ and $U^\ast$ are the vector space duals of the corresponding
spaces. However, our inner products will always be the standard dot
product, and we will identify the vector spaces and their duals. Thus
we can get away with the slightly informal notation used in the
diagram above.

Whenever we have a situation as in~\eqref{eq:UVW} above, the middle
space splits into three subspaces. A splitting of $V$ into two parts
is just a consequence of the fact that $V$ consists of the subspace
$\im A$ and its orthogonal complement $\im A^\perp = \ker A^T$. The
presence of the second map $B$ and the fact that $B\circ A = 0$ is the
crucial ingredient for getting a further splitting of $\ker A^T$.
Just as the first split comes from two of the fundamental subspaces,
the finer splitting of one of the pieces is yet another use of the
fundamental subspaces ideas. These ideas are made more precise in the
following elementary fact.

\begin{fact} \label{lem:Hodge}
  There exists a unique orthogonal decomposition of $V$ (called the Hodge
  decomposition) as:
  \[
  V = \im A \oplus \im B^T \oplus \ker \laplacian \, .
  \]
  Moreover, $\ker \laplacian = \ker B \cap \ker A^T$.
\end{fact}
\begin{proof} We have first the obvious decomposition $V = \im A
  \oplus (\im A)^\perp$, where $(\im A)^\perp$ means the orthogonal
  complement of $\im A$. Thus $V = \im A \oplus \ker A^T$, from which
  follows that $V = \im A \oplus \im B^T \oplus ((\im B^T)^\perp \cap
  \ker A^T)$. This is due to the fact that $A^T\circ B^T = 0$ because
  of which $\im B^T \subset \ker A^T$. This finally yields $V = \im A
  \oplus \im B^T \oplus (\ker B \cap \ker A^T)$.  To prove that $\ker
  \laplacian = \ker B \, \cap \, \ker A^T$, it is trivial to verify
  that $\ker B \, \cap \, \ker A^T \subset \ker \laplacian$.  For the
  other direction, let $h \in \ker \laplacian$. Then $0 =
  \aInnerproduct{\laplacian h}{h} = \aInnerproduct{A^Th}{A^Th} +
  \aInnerproduct{Bh}{Bh} $ from which the result follows. Here the
  three inner products above are on $V$, $U$ and $W$, respectively.
\end{proof}

To be precise, one should write $V \cong \im A \oplus \im B^T \oplus
\ker \laplacian$, since $B^T : W^\ast \to V^\ast$. However, we will
continue to use equality by identifying the dual spaces $V^\ast$ with
the corresponding original vector spaces $V$ etc. as mentioned
earlier.

\section{Least squares ranking and Hodge decomposition}
\label{sec:rnknghdg}

We first consider the case when a given pairwise data, i.e., cochain
$\omega \in C^1(G)$ has components along both $\im \boundary_1^T$ and
$\im\boundary_2$. In this case, the Hodge decomposition, least
squares, normal equations, and the Karush-Kuhn-Tucker equations are
equivalent. This is the content of Theorem~\ref{thm:kernel_none}. The
restrictions can be dropped to prove analogous theorems involving
fewer equations.

In a least squares problem $Ax\simeq b$, to minimize $\norm{b-Ax}_2^2$
as a function of $x$, a necessary condition is that the gradient be
zero which yields the normal equations. Thus residual minimization
implies the normal equations. For the converse, often a sufficient
condition that is described in text books is that the Hessian matrix
(which is $2A^TA$, in this case) be positive definite (see for
example, \cite[page 110]{Heath2002}). This is often useful in the
classical least square case in which $m \ge n$. For then, if $A$ is
full rank $A^T A$ is positive definite. In our case, the matrices $A^T
A$ will be $\boundary_1 \, \boundary_1^T $ or
$\boundary_2^T\,\boundary_2$. In most complexes with interesting
topology we \emph{cannot} rely on these to be nonsingular since
$\boundary_1^T$ and $\boundary_2$ will have nontrivial kernels. The
constant functions on the vertices constitute $\ker \boundary_1^T$ and
all spheres are in $\ker \boundary_2$.  As an alternative, we
will use the following lemma.

\begin{lemma}\label{lem:ne_ls} Given a matrix $A \in \R^{m\times n}$
  and a vector $x_\ast \in \R^n$, if $A^TAx_\ast = A^Tb$ and $x_\ast \notin
  \ker A$ then $x_\ast$ minimizes the residual norm $\norm{b-Ax}_2$
  over all $x\in\R^n$.
\end{lemma}
\begin{proof}
  The dot product $\ainnerproduct{b-Ax_\ast}{Ax_\ast} =
  \ainnerproduct{A^T(b-Ax_\ast)}{x_\ast} = 0$. This means either
  $b-Ax_\ast = 0$, in which case we are done, or the vectors
  $b-Ax_\ast$ and $Ax_\ast$ are orthogonal. The latter means that the
  shortest distance from $b$ to $\im A$ is achieved by $Ax_\ast$.
\end{proof}

\begin{theorem}[Basic Fact]\label{thm:kernel_none}
  Given $\omega$, $h \in C^1(G)$, $\alpha \in C^0(G)$, $\beta \in
  C^2(G)$, with $\alpha \notin \ker \boundary_1^T$ and $\beta
  \notin \ker \boundary_2$, the following are equivalent.
  \begin{enumerate}[(i)]
  \item \label{itm:hodge} Hodge Decomposition (HD)
    \begin{align}
      \omega &= \boundary_1^T\, \alpha + \boundary_2\,
      \beta + h \, ,\label{eq:hodge}\\
      h & \in \ker \laplacian_1\,.\notag
    \end{align}

  \item \label{itm:lstsqrs}Least Squares (LS)

   \noindent $a = \alpha, b = \beta$, and $s=h$ are optimal values of
    the two least squares problems
   \begin{align} 
      \underset{a}{\min} \, \norm{r}_2 \quad
      \text{such that}\quad r & = \omega - \boundary_1^T \;
      a\, ,\label{opt:ls1}\\
      \underset{b}{\min} \, \norm{s}_2 \quad
      \text{such that}\quad s &= r_\ast - \boundary_2 \, b\,,
      \label{opt:ls2}
    \end{align}
    where $r_\ast$ is the minimizing residual for~\eqref{opt:ls1}. In
    least squares short hand notation one would write the two problems
    as $\boundary_1^T a \simeq \omega$ and $\boundary_2\, b \simeq
    r_\ast$.

  \item \label{itm:normal} Normal Equations (NE)

   \noindent $a = \alpha$ and $b = \beta$ are a solution of the two
    linear systems
    \begin{align}
      \boundary_1 \, \boundary_1^T \, a &=
      \boundary_1\, \omega \, ,\label{eq:normal1}\\
      \boundary_2^T\,\boundary_2\, b &= \boundary_2^T\,
      r_\ast\, ,\label{eq:normal2}
    \end{align}
    where $r_\ast$ is the residual $\omega - \boundary_1^T\,\alpha$.
  \item \label{itm:kkt} Karush-Kuhn-Tucker Equations (KKT)

   \noindent $a = \alpha, b = \beta$, and $s = h$ are a solution of
    the two saddle-type systems
   \begin{align}
     \begin{bmatrix}
       I & \boundary_1^T \\
       \boundary_2^T & 0
     \end{bmatrix}\,
     \begin{bmatrix}
       r \\ a
     \end{bmatrix} &=
     \begin{bmatrix}
       \omega\\ 0
     \end{bmatrix}\label{eq:kkt1}\, ,\\
     \begin{bmatrix}
       I & \boundary_2 \\
       \boundary_1 & 0
     \end{bmatrix}\,
     \begin{bmatrix}
       s \\ b
     \end{bmatrix} &=
     \begin{bmatrix}
       r_\ast\\ 0
     \end{bmatrix}\, ,\label{eq:kkt2}
   \end{align}
   where $r_\ast$ is part of the solution for the first system.
 \end{enumerate}
\end{theorem}
\begin{proof}
  Follows from Lemma~\ref{lem:ne_ls}, elementary calculus and
  linear algebra.
\end{proof}

\begin{remark} \label{rem:kernel_some}
  One can prove analogous theorems for the case when $\alpha \in
  \ker\boundary_1^T$ and/or $\beta \in \ker \boundary_2$. This
  would involve skipping the equations corresponding to the term that
  is in the kernel. If both are, then the given data is purely
  harmonic. 
\end{remark}

\begin{remark}
  The existence of the Hodge decomposition comes from
  Fact~\ref{lem:Hodge}. The theorem above states the equivalence of
  Hodge decomposition with least squares, normal and
  Karush-Kuhn-Tucker equations.
\end{remark}

\subsection{Implications of orthogonality}\label{subsec:orthgnlty}

The three terms in the Hodge decomposition~\eqref{eq:hodge} are
mutually orthogonal. This is easy to see. It follows simply from the
fact that $\boundary_1\boundary_2 = 0$ and from the definition of the
harmonic part ($\ker \laplacian_1$). For example, given an $\omega \in
C^1(G)$, if it has a nonzero harmonic part $h$, then $\innerproduct{h}
{\boundary_1^T\, \alpha} = \innerproduct{\boundary_1\, h}{\alpha} = 0$
since $h$ is in $\ker\boundary_1 \cap \ker\boundary_2^T$.

Due to these orthogonality conditions, it is easy to see that the
second least squares problem, which is $\boundary_2\, b \simeq
r_\ast$, can also be written as $\boundary_2 \, b \simeq
\omega$. Similar changes can be made from $r_\ast$ to $\omega$ in the
second systems in all the formulations above. For ease of reference,
below we write the least squares and normal equations using $\omega$
instead of $r_\ast$ all in one place. The least squares systems are
\begin{align}
  \boundary_1^T \, a &\simeq \omega \, ,\label{eq:lsomega1}\\
  \boundary_2\, b &\simeq \omega\, ,\label{eq:lsomega2}
\end{align}
and the corresponding normal equations
\begin{align}
  \boundary_1 \, \boundary_1^T \, a &=
  \boundary_1\, \omega \, ,\label{eq:nrmlomega1}\\
  \boundary_2^T\,\boundary_2\, b &= \boundary_2^T\,
  \omega\, .\label{eq:nrmlomega2}
\end{align}
Note from the definition of the Hodge Laplacians
in~\eqref{eq:lplcdrhm} that the above normal equations can be written
as $\laplacian_0\, a = \boundary_1\, \omega$ and $\laplacian_2\, b =
\boundary_2^T\, \omega$.

\subsection{Connection with optimal homologous chains problem}
\label{subsec:OHCP}

Notice that the least squares formulations in~\eqref{opt:ls1}
and~\eqref{opt:ls2} are exactly analogous to the optimal homologous
chain problem from~\cite{DeHiKr2011}. For example, in~\eqref{opt:ls1}
the given cochain is $\omega$ and one is looking for the smallest
cochain $r$ which is cohomologous to it. In~\cite{DeHiKr2011} it is
the 1-norm of chains that is minimized over all homologous chains. In
contrast, here the 2-norm of \emph{co}chains is minimized over all
\emph{co}homologous cochains. This results in our solving linear
systems in this paper as opposed to linear programming which is used
in~\cite{DeHiKr2011}. In spite of these differences between the two
problems, it is interesting that the problem of least squares ranking
on graphs and a fundamental problem of computational topology
(computing optimal homologous chains) are related.

\subsection{Interpretation in terms of ranking}

Given any pairwise comparison data $\omega \in C^1(G)$ we see from the
Theorem~\ref{thm:kernel_none} that there exists a cochain $\alpha \in
C^0(G)$ (the vertex potential or ranking), a cochain $\beta \in
C^2(G)$, and a harmonic field $h \in \ker \laplacian_1$, such that
$\omega = \boundary_1^T \alpha + \boundary_2 \beta + h$. The $\alpha$
term is the scalar potential that gives the ranking. The $\beta$ term
is defined on cells and captures the local inconsistency in the
data. The harmonic part contains the inconsistency that is present due
to loops longer than the maximum number of sides in the cells. If only
3-cliques (triangles) are considered as the 2-dimensional cells, then
any inconsistency in loops of length four or more will be captured in
the harmonic part. The following example should make some of this more
apparent.

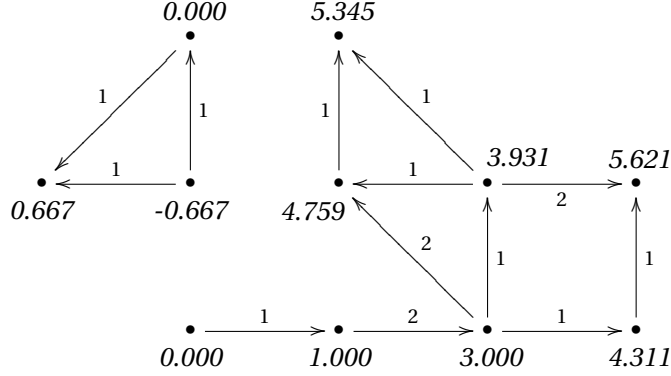
\begin{figure}[ht]
\[
    \xymatrixrowsep{0.625in}
    \xymatrixcolsep{0.625in}
    \xymatrix{ & {\bullet}\ar[ld]_1 \save[]+/u2.0ex/*\txt{
    \emph{0.000}}\restore & {\bullet} \save[]+/u2.0ex/*\txt{\emph{
    5.345}}\restore& & \\
    {\bullet} \save[]+/d2.0ex/*\txt{\emph{0.667}}\restore & {\bullet}
    \ar[l]_1\ar[u]_1 \save[]+/d2.0ex/*\txt{\emph{-0.667}}\restore & 
    {\bullet}\ar[u]^1 \save[]+/dl3ex/*\txt{\emph{4.759}}\restore & 
    {\bullet}\ar[l]_1\ar[ul]_1\ar[r]_2 \save[]+/ur3ex/*\txt{
    \emph{3.931}}\restore & {\bullet} \save[]+/u2.0ex/*\txt{
    \emph{5.621}}\restore \\
     & {\bullet}\ar[r]^1 \save[]+/d2.0ex/*\txt{\emph{0.000}}\restore &
     {\bullet}\ar[r]^2 \save[]+/d2.0ex/*\txt{\emph{1.000}}\restore & 
    {\bullet}\ar[ul]_2\ar[u]_1\ar[r]^1 \save[]+/d2.0ex/*\txt{
    \emph{3.000}}\restore & {\bullet}\ar[u]_1 \save[]+/d2.0ex/*\txt{
    \emph{4.311}}\restore
}
\]
  \caption{The results of the first least squares problem
    \eqref{opt:ls1} on a graph.}
  \label{fig:v11e13}
\end{figure}

\begin{example}
  Figure~\ref{fig:v11e13} shows an example of solving the first least
  squares problem~\eqref{opt:ls1}. It is just as easy to work with
  disconnected graphs, so we show a graph with two components. The
  values on the edges is the given data $\omega$ in $C^1(G)$. The
  vertex potential values $\alpha \in C^0(G)$ are written in
  italics. Note that in the straight line part of the graph which does
  not involve a cycle, it is clear what the vertex potential should be
  (up to an additive constant). There will be no residual in this
  case. The first triangle after the straight line part is consistent
  because the value on the hypotenuse is the sum of the values on the
  other two sides which are oriented appropriately. The other
  triangles and the square loop are all inconsistent. Here only
  triangles are chosen as the 2-dimensional cells, so the $\beta$ part
  will be the inconsistency associated with the triangles if the
  second problem were also to be solved. The harmonic part $h$ would
  be the inconsistency in the square loop. Note that because of two
  connected components the dimension of $\ker \laplacian_0 =
  \boundary_1\,\boundary_1^T$ will be two. Fixing one vertex value in
  each of the two components and deleting the appropriate row and
  column will make the normal equations system~\eqref{eq:normal1}
  nonsingular. We will discuss the issue of nontrivial kernels in the
  context of linear solvers in Section~\ref{sec:solvers}.
\end{example}

\section{Experiments on Topology of Clique
 Complexes}\label{sec:clique}

A recent paper by Kahle~\cite{Kahle2009} explores the homology of
clique complexes arising from Erd\H{o}s-R\'{e}nyi random graphs. The
number of connected components of a graph is the same as the dimension
of its 0-dimensional homology. The connectedness of random graphs has
been explored in literature for many years and the study of higher
dimensional homology can be considered as the new and natural
extension of that line of research. For the ranking problem,
1-homology is of particular interest. When 1-homology of a graph is
trivial, then by Remark~\ref{rem:hrmncbtt} there cannot be any
harmonic component in any 1-cochain. Then if there are no local
inconsistencies (i.e., the curl part is zero) a 1-cochain will be a
pure gradient, hence globally consistent.

Kahle provides bounds on the edge density for Erd\H{o}s-R\'{e}nyi
graphs for which the $p$-dimensional homology is almost always
trivial, and bounds for which it is almost always nontrivial. Since
all results of this type are ``almost always'', in the rest of this
section, we will omit that phrase. We restrict our attention to
1-homology, and the 2-dimensional clique complex of the graph is the
relevant object. In this setting, Kahle states that for an
Erd\H{o}s-R\'{e}nyi graph with $n$ nodes and edge density $\rho$, the
1-homology will be trivial when $\rho < 1/n$ or when $\rho >
1/\sqrt[3]{n}$. The 1-homology will be nontrivial when $1/n < \rho <
1/\sqrt{n}$. When $1/\sqrt{n} < \rho < 1/\sqrt[3]{n}$ there is a
theoretically undetermined transition in homology. (There is a typo in
Jiang et al.'s quotation of the bounds from Kahle's results.)

The numerical framework for exploring least squares ranking on graphs
can be applied to study clique complexes, almost without any
changes. For example, it is satisfying to see Kahle's bounds appear in
the experimental results shown in the first column of
Figure~\ref{fig:clqcmplxs}. But what is more interesting is that
experiments like these can serve as tools for developing new
conjectures and asking more detailed questions than are answered by
the current theory. Once a conjecture looks numerically plausible, one
can set about trying to find a mathematical proof. But before we
mention some such questions, we note that the apparent violation of
Kahle's bounds in Figure~\ref{fig:clqcmplxs} is not really a
violation, because his bound are true in the limit as the number of
vertices $n$ goes to infinity. For example, there are some nonzero
homology points in the region that is supposed to almost always have
trivial homology.

Using this experimental tool one can ask new questions, such as, what
is the behavior of the 1-dimensional Betti number as a function of
$\rho$? By Remark~\ref{rem:hrmncbtt}, the Betti number can be measured
by measuring the dimension of the space of harmonic cochains, i.e., 
$\dim \ker \laplacian_1$. Since $\ker \laplacian_1 = \ker \boundary_1
\cap \ker \boundary_2^T$, one can measure the Betti number by stacking 
the matrices for $\boundary_1$ and $\boundary_2^T$, one on top of the
other and the kernel dimension of this matrix yields the desired Betti
number. The kernel dimension can be found by computing a singular
value decomposition and counting the number of zero singular values. A
faster alternative is to compute the number of zero eigenvalues of
$\laplacian_1$ by using a sparse eigensolver. Usually, the distinction
between what should be considered nonzero and what should be
considered zero is very evident in our experiments.

Results about the Betti number give a more nuanced picture than the
presence or absence of homology. For example, in terms of Betti
number, one can ask if the transition from nonzero homology region to
the zero homology region on the right is sudden. Our experiments on
\ER~graphs shed some light on these questions. In
Figure~\ref{fig:clqcmplxs}, the bottom graph in the first column shows
how the Betti number varies as a function of $\rho$. A clear trend is
visible and the Betti number appears to peak at the start of the
transition zone where the theory is silent (the bounds of Kahle are
marked as dashed vertical lines).  The quantification of homology can
also take another form in these experiments. One can investigate how
much of the norm of a 1-cochain is contained in the harmonic
component, when the homology is nontrivial. Or, how does this harmonic
component vary as a function of $\rho$? The results are in the top two
graphs in the left column of Figure~\ref{fig:clqcmplxs}. The right
column shows results of experiments on topology of clique complex of
\BA~scale-free graphs, for which no theory has yet been
developed. Although Kahle's bounds were not developed for \BA~graphs,
we have drawn those for the plots corresponding to \BA~graphs in
order to provide context and comparison with the \ER~results.

\section{Comparing Linear Solvers}\label{sec:solvers}

In most areas of numerical analysis, the computation that sits at the
heart of the solution method is usually the solution of a linear
system. Least squares ranking on graphs is no different. This section
is about our numerical experiments for testing the accuracy and speed
of linear system solvers for the first and second least squares
problems of ranking on graphs. 

A small part of our experimental work was on graphs with special
structure, such as path, cycle, star, and wheel graphs. However, our
main focus has been on two popular random graph models, and on
scale-free graphs. The earliest and perhaps most studied random graph
model is the one of Erd\H{o}s and R\'{e}nyi~\cite{ErRe1960,
  Bollobas2001}. The two parameters for this model are the number of
nodes and the probability that any pair of nodes is connected by an
edge. We will refer to this probability as the \emph{edge
  density}. This model does not account for the phenomena of
clustering that is seen in many networks in societies and this
shortcoming is overcome by the model of Watts and Strogatz
model~\cite{NeStWa2001}. However, neither of these account for the
power-law shape of degree distributions that are seen in many
real-life graphs such as the World Wide Web, reaction networks of
molecules in a cell, airline flight networks and so on. These are
called scale-free networks or graphs and this feature is captured by
the model of Barab\'asi and Albert~\cite{BaAl1999}. The generative
model is often implemented as a random process, although in practice
these are typically not random graphs.

In the numerical analysis community there is a lot of accumulated
experience on solving linear systems that arise from partial
differential equations. Studies of systems arising from graphs are
less common but appearing with increasing frequency. For least squares
ranking on general graphs there is no guidance available in the
literature. Ours is a first step in an attempt to fill that gap. The
fact that the underlying problem is coming from a graph introduces
some new challenges as we will demonstrate via our experiments in this
section.

We only consider iterative linear solvers here, though sparse direct
solvers might be worth considering. We used a variety of iterative
Krylov methods suitable for symmetric systems and one that is suitable
for rectangular systems. We also used algebraic multigrid using
smoothed aggregation and Lloyd aggregation. The results are discussed
in Sections~\ref{subsec:krylov}--\ref{subsec:schur}. An especially
attractive feature of all these methods is their ability to ignore the
kernel of the operator involved. We do not know how the direct methods
could be made to do that for the second least squares problem which
can have a large dimensional kernel depending on the graph
topology. In direct solvers, for the first least squares problem the
nontrivial kernel can be handled by fixing the value at a single
vertex, just as is done by fixing pressure at a point in fluid
problems.

\begin{remark}
  For an arbitrary graph $G$, by Theorem~\ref{thm:kernel_none} (or its
  special cases mentioned in Remark~\ref{rem:kernel_some}), the first
  least squares problem of ranking~\eqref{eq:lsomega1} is equivalent
  to the the normal equation~\eqref{eq:nrmlomega1}. But the matrix
  involved is then $\boundary_1\,\boundary_1^T$ which is the
  combinatorial graph Laplacian $\laplacian_0$ and hence it is
  symmetric and diagonally dominant. Thus, it can be solved by the
  method of Koutis et al.~\cite{KoMiPe2010}, which can solve such
  systems in time approaching optimality. However, at the time of
  writing, no reliable implementation of the method of Koutis et 
  al.~was available. The system matrix $\laplacian_2$ for a general 
  graph need not be diagonally dominant. Consider for example the 
  complete graph $K_5$ in which every 3-clique is taken to be a 
  triangle. It is easy to verify that $\laplacian_2$ is not diagonally
  dominant. It is a $10\times 10$ matrix with 3s along the diagonal 
  and with each row containing six off-diagonal entries that are 
  $\pm 1$ (with four entries that are 1 and two that are $-1$). Thus, 
  for a general graph, the Koutis et al.~solver cannot be used for 
  solving the second least squares problem of ranking due to lack of 
  diagonal dominance.
\end{remark}

\subsection{Methodology}\label{sctn:mthdlgy}

All numerical experiments were done using the Python programming
language. The Krylov linear solvers used were those provided in the
SciPy module~\cite{JoOlPe2001}. The algebraic multigrid used was the one
provided in PyAMG~\cite{BeOlSc2011}. The simplicial complexes and
boundary matrices were created using the PyDEC
module~\cite{BeHi2011}. The errors and times required by various
solvers is generated as an average over multiple trials. This is done
to minimize the influence of transient factors that can affect
performance of a computer program. All experiments were carried out on
an Apple MacBook with a 2.53~GHz Intel Core~2 Duo processor and with
4~GB of memory.

In each case, a graph $G$ with the desired number of nodes $N_0$ and
other desired characteristics (such as edge density in the case of
Erd\H{o}s-R\'{e}nyi graphs) is first generated by a random process.
We then find all the 3-cliques in the graph and create a simplicial
complex data structure for the resulting 2-complex. Let the number of
edges and triangles in $G$ be $N_1$ and $N_2$.

A random ranking problem instance is created for this complex. This
entails creating a random 1-cochain representing the comparison data
on edges. The point of these experiments is to compare the accuracy
and efficiency of the tested methods, and so the Hodge decomposition
of this 1-cochain has to be known in advance. In other words, a random
problem instance is a 1-cochain $\omega$ such that there are random
but known $\alpha \in C^0(G)$, $\beta \in C^2(G)$, and $h \in \ker
\laplacian_1$ with $\omega = \boundary_1^T \alpha + \boundary_2 \beta
+ h$.

It is clear how to create the random gradient part $\boundary_1^T
\alpha$ and the random curl part $\boundary_2 \beta$ -- simply pick a
random vector with $N_0$ entries for $\alpha$ and a random vector with
$N_2$ entries for $\beta$. To compute a random harmonic part, one can
compute the Hodge decomposition of a random 1-cochain $\rho$ by
solving two least squares problems $\boundary_1^T a \simeq \rho$ and
$\boundary_2 b \simeq \rho$ as outlined in
Section~\ref{subsec:orthgnlty}. If $\alpha$ and $\beta$ are the
solutions, then $\rho - \boundary_1^T \alpha - \boundary_2 \beta$ is
a desired random harmonic cochain. 

For variety we show here another method, and this is the one we
used. It relies on the basic fact that the the residual $b-Ax$ in a
least squares problem $Ax \simeq b$ is orthogonal to the $\im A$ and
hence in $\ker A^T$. We pick a random 1-cochain $\rho$, that is, a
random vector with $N_1$ entries. We then solve the single least
squares problem
\begin{equation}\label{eq:rndmhrmnc}
  \begin{bmatrix}
    \boundary_1^T & \boundary_2
  \end{bmatrix}\, x \simeq \rho \, ,
\end{equation}
where the matrix $[\boundary_1^T\,\boundary_2]$ is formed by
horizontally stacking $\boundary_1^T$ and $\boundary_2$ matrices.
If $x$ is the solution of this least squares problem, then the
residual $\rho - [\boundary_1^T \boundary_2]\, x$ is harmonic. This is
because
\[
\rho - [\boundary_1^T \boundary_2]\, x \;\in\; 
\ker\, [\boundary_1^T \boundary_2]^T \;= \;
\ker \begin{bmatrix}  \boundary_1\\  \boundary_2^T \end{bmatrix}\;=\;
\ker \boundary_1 \cap \ker \boundary_2^T\; =\; 
\ker \laplacian_1\, .
\]

\subsection{Spectral analysis} \label{sub:spectral}

Recall that if $A$ is a symmetric positive definite matrix, the number
of conjugate gradient iterations is related to the norm of the error
by the inequality
\begin{equation} \label{eq:bound}
  \dfrac{\norm{e_k}_A}{\norm{e_0}_A} \leq 
  2\left(\dfrac{\sqrt{\kappa} - 1}{\sqrt{\kappa} + 1}\right)^k\, .
\end{equation}
See for example \cite[page 51]{Greenbaum1997}. Here $\norm{x}_A$ is
the $A$-norm of $x$, i.e., $\norm{x}^2_A := \ainnerproduct{x}{Ax}$ and
$\kappa = \lambdamax/\lambdamin$ is the condition number of $A$. Here
$\lambdamax$ and $\lambdamin$ are the largest and smallest magnitude
eigenvalues of $A$, respectively. The same result holds even if $A$ is
singular, as long as it is semidefinite, $\norm{x}_A$ is considered a
seminorm, and $\lambdamin$ is defined to be the smallest (in
magnitude) nonzero eigenvalue of $A$. Given a desired error $\epsilon$
(in the $A$-norm), one can find the number of iterations of conjugate
gradient method required to achieve that error, by substituting
$\epsilon$ for $\norm{e_k}_A$ and solving for the smallest $k$ which
satisfies inequality~\eqref{eq:bound}. We will refer to this number as
\emph{conjugate gradient iterations} required to achieve error
$\epsilon$.  Thus for an arbitrary graph $G$, the conjugate gradient
iterations required for the first least squares problem of ranking is
given by inequality~\eqref{eq:bound} using $\lambdamin$ and
$\lambdamax$ for the graph Laplacian $\laplacian_0$.

Much is known about spectrum of the graph Laplacian for various types
of graphs \cite{Fiedler1973, Mohar1991, Merris1995, LiZh1998},
including random graphs and scale-free networks
\cite{FaLuVu2003}. These results usually involve some graph
property. For example, for various types of special graphs,
$\lambdamin$ is often bounded in terms of edge connectivity -- the
minimum number of edges to be removed to disconnect a graph. Thus one
can make predictions like single iteration convergence of conjugate
gradient method in the case of complete graphs, and this is borne out
by our numerical experiments. For some special graphs, some well known
lower bounds or formulas for $\lambdamin$ are given below. An easy
upper bound for $\lambdamax$ of $\laplacian_0$ is twice the maximum
degree. This follows from Gerschgorin's theorem~\cite{TrBa1997}.  In
the table below, $\eta(G)$ is the edge connectivity, a path graph is
an acyclic graph in which all but two vertices have degree two, and
the two ``end'' vertices each have degree one. A cycle graph is a
2-regular graph on which there is only one cycle of size $n$, and a
star graph is one in which all but one vertex have degree 1, and each
is connected to a ``center'' vertex which has
degree~$n-1$. Figure~\ref{fig:cgbounds1} in
Appendix~\ref{appndx:cgbnds} shows comparisons of iteration bounds and
actual iteration counts that occur in numerical experiments.
\begin{center}
\begin{tabular}{|c|c|}
\hline
Type of graph & $\lambdamin$ \\
\hline
General & $2\,\eta(G)\,(1 - \cos(\pi/n))$ \\
Complete & $n$ \\
Path & $2\,(1 - \cos(\pi/n))$ \\
Cycle & $2\,(1 - \cos(2\pi/n))$ \\
Star & $1$ \\
\hline
\end{tabular}
\end{center}

For \ER, \WS, and \BA~graphs, there are fewer results for the spectral
radii but there exist results bounding these from which iteration
bounds can be obtained. Results of some experiments using such graphs
are shown in Figure~\ref{fig:cgbounds2}.

\begin{remark}
  When the simplicial (or cell) 2-complex representation of $G$ can be
  embedded as a meshing of a compact surface with or without boundary,
  we \emph{can} place bounds on the spectrum of $\laplacian_2$. We do
  this by generating the \emph{dual graph} $G^D$ of $G$, in which
  every triangle of $G$ is a vertex in $G^D$; the dual vertices in
  $G^D$ are connected by edges which are dual to those in $G$. One can
  then bound the spectrum of $\laplacian_2$ in both the boundary and
  boundaryless surface cases using Cauchy's interlacing theorem for
  eigenvalues~\cite{Parlett1998}.  However, while any graph may be
  embedded in a surface of sufficiently high genus, the embedding of
  the simplicial 2-complex is not possible for \emph{general}
  graphs. As an example, consider the complete graph $K_5$ introduced
  above, for which $N_0 = 5$, $N_1 = 10$, and $N_2 = 10$. The Euler
  characteristic of this complex is $\chi = N_0 - N_1 + N_2 =
  5$. Assume that $K_5$ could be embedded on a surface of genus $g$
  and with $b$ disjoint portions of boundary. Then $\chi = 2 - 2g - b$
  which reduces to $2g + b = -3$ which is not true for any nonnegative
  $g$ and $b$ leading to a contradiction.
\end{remark}

\subsection{Iterative Krylov methods}\label{subsec:krylov}

The Krylov solvers that we used in our experiments are conjugate
gradient (CG) and minimal residual (MINRES) for normal equations and
saddle formulation~~\cite{Saad2003, Vandervorst2009}, and
LSQR~\cite{PaSa1982a, PaSa1982b} for the least squares system solved
using the given rectangular matrix without forming the square system.
There are many other Krylov solvers we did not test. For example, we
did not test GMRES since all of our square linear systems are
symmetric. We did not test CGLS or CGNE since they are mathematically
equivalent to LSQR, which tends to be more commonly used. As mentioned
earlier, the Laplacian systems will in general be symmetric positive
semidefinite. Krylov methods have the nice property that they work in
spite of a nontrivial kernel. This feature is especially effective in
the graph problem because the matrices are integer matrices requiring
no quadrature for their construction such as is required in the case
of partial differential equations~\cite{BoLe2005}.

Tables~\ref{tab:erkrylv}--\ref{tab:bakrylv} are timing and error
results for various Krylov solvers on different formulations of the
ranking problem. In each table, edge and triangle densities are with
respect to the number of edges and all possible triangles in the
complete graph. LSQR was used on the least squares equations directly,
i.e., on the rectangular matrices in~\eqref{eq:lsomega1}
and~\eqref{eq:lsomega2}. The other solvers are used on the normal
equations unless the `-K' designator is listed, which indicates the
solver was used on the Karush-Kuhn-Tucker formulation~\eqref{eq:kkt1}
and~\eqref{eq:kkt2}. Reported errors are measured relative to the
known exact solution except in cases identified by an asterisk
($\ast$) in $\norm{h}$ column where the corresponding error is an
absolute one. These are cases for which the homology of the simplicial
2-complex induced from the graph is trivial leading to a zero harmonic
component.  The relative error column reports error in the norm of the
gradient part ($\norm{\boundary_1^T \alpha}$), norm of the curl part
($\norm{\boundary_2 \beta}$), and the norm of the harmonic part
($\norm{h}$). The timing labeled $\alpha$ shows the iterations and
time required for the first least squares problem, and the one labeled
$\beta$ shows these for the second least squares problem.

\subsection{Algebraic multigrid methods}\label{subsec:amg}

Multigrid methods work by creating a hierarchy of linear systems of
decreasing sizes from the given problem.  At any level in the
hierarchy, a larger system is called a \emph{fine grid} and a smaller
system is called a \emph{coarse grid}. The solutions of the different
systems are related via prolongation (coarse to fine) and restriction
(fine to coarse) maps. In geometric multigrid, coarser levels
correspond to a coarser mesh. In \emph{algebraic} multigrid,
coarsening is carried out using only the matrix and an associated
adjacency graph. Coarsening is performed by aggregating those vertices
of this graph which have a strong connection. The strength of
connection is defined in various ways for different
schemes~\cite{Stuben2001}.  In smoothed aggregation, a vertex can
belong fractionally to several aggregates~\cite{VaMaBr1996}. In Lloyd
aggregation, the number of connections between vertices and the
centers of their aggregates is minimized~\cite{Bell2008}.

We used algebraic multigrid with smoothed aggregation and Lloyd
aggregation for the ranking problem on Erd\H{o}s-R\'enyi and
Watts-Strogatz random graphs, and Barab\'asi-Albert scale-free
graphs. The parameters for the various graphs were identical to ones
used in the experiments with Krylov solvers. The results of these
numerical experiments are given in
Tables~\ref{tab:eramg}--\ref{tab:baamg}. The columns are the same as
the ones in the Krylov
tables~\ref{tab:erkrylv}--\ref{tab:bakrylv}. Only the normal equation
formulations~\eqref{eq:nrmlomega1} and~\eqref{eq:nrmlomega2} are used
in these experiments. In the Algorithm/Formulation column, the
notation AMG (SA) and AMG (LA) are used to indicate the smoothed and
Lloyd aggregation schemes. In this column, Schur indicates yet another
variation which is described as Algorithm~\ref{alg:schur} in the next
subsection. As before, in these tables, the cases where the homology
of the simplicial 2-complex is trivial are marked by an asterisk
($\ast$). The cases where PyAMG failed in the setup phase are marked
by a dagger symbol ($\dagger$), and the cases where algebraic
multigrid reached maximum number of specified iterations are marked by
the double dagger symbol ($\ddagger$).  The entries marked with a dash
symbol (--) are those for which Algorithm~\ref{alg:schur} cannot be
applied because $\laplacian_2$ has no simple sparse/dense
partitioning. Tables~\ref{tab:eramgsttstcs}--\ref{tab:baamgsttstcs}
list the measured setup and solve times for each of ranking problems
on the three graph models where solution could be found using
algebraic multigrid. 

\subsection{Schur complement and algebraic multigrid}
\label{subsec:schur}

A third approach is a hybrid one that combines algebraic multigrid
with Krylov or direct solvers. The idea is to partition the problem
into a small dense part and a large sparse part.  The small part can
be solved efficiently, e.g., by Krylov or direct methods, and the
large sparse part can be solved by algebraic multigrid. A general
matrix $A$ is partitioned as
\begin{equation*}
A = \begin{bmatrix} A_{11} & A_{12} \\ A_{21} & A_{22} \end{bmatrix}\, ,
\end{equation*}
where $A_{11}$ is the large sparse part and $A_{22}$ is the small 
dense part. The linear system $A x = b$ can be written as
\[
  \begin{bmatrix} A_{11} & A_{12} \\ A_{21} & A_{22} \end{bmatrix}
  \begin{bmatrix}x_1 \\ x_2\end{bmatrix} = 
  \begin{bmatrix}b_1 \\ b_2\end{bmatrix}\notag \, ,
\]
which can be reduced by row operations to
\begin{equation} \label{eq:schur}
  \begin{bmatrix} A_{11} & A_{12} \\ 0 & A_{22} - 
    A_{21} A_{11}^{-1} A_{12} \end{bmatrix}
  \begin{bmatrix}x_1 \\ x_2\end{bmatrix} = 
  \begin{bmatrix}b_1 \\  b_2 - A_{21} A_{11}^{-1} b_1
  \end{bmatrix} \, . 
\end{equation}
One can then solve the second block for $x_2$ first and use the result
to solve for $x_1$. This is the Schur complement formulation. The hope
is that $A_{22}$ will be small even if it is dense, and that $A_{11}$
will be sparse. Then, if $A_{11}^{-1}$ can be found easily, forming
$A_{22} - A_{21} A_{11}^{-1} A_{12}$ will be easy and $x_2$ can be
found perhaps using a dense but small linear system.

\begin{algorithm}[t]
  \begin{algorithmic}[1]
  \REQUIRE Approximate inverse $\widetilde{A}_{11}^{-1}$ and initial
  guess $x^{(0)}$.

  \STATE Compute $r^{(0)} = \begin{bmatrix} r_1^{(0)} \\
    r_2^{(0)} \end{bmatrix} = \begin{bmatrix} b_1 ^{\phantom{0}} \\
    b_2 ^{\phantom{0}} \end{bmatrix} - \begin{bmatrix} A_{11}^{\phantom{0}} &
    A_{12}^{\phantom{0}} \\ A_{21}^{\phantom{0}} &
    A_{22}^{\phantom{0}} \end{bmatrix} \begin{bmatrix} x_1^{(0)} \\
    x_2^{(0)} \end{bmatrix}$

  \FOR {$i = 0 \, , 1 \, , \dots \, ,$ until convergence}
  
  \STATE Krylov/Direct-Solve: $(A_{22} - A_{21} \widetilde{A}_{11}^{-1} A_{12})
  e_2^{(i)} = r_2^{(i)} - A_{21} \widetilde{A}_{11}^{-1} r_1^{(i)}$
  
  \STATE AMG-Solve: $A_{11} e_{1}^{(i)} = r_1^{(i)} - A_{12}
  e_{2}^{(i)}$
  
  \STATE $x^{(i + 1)} = x^{(i)} + e^{(i)}$

  \STATE $r^{(i + 1)} = b - A x^{(i + 1)}$
  
  \ENDFOR
\end{algorithmic}
\caption{Iterative solution of Schur complement system using
  Krylov/direct and AMG solvers}
\label{alg:schur}
\end{algorithm}

In order to solve the dense block of~\eqref{eq:schur}, one needs the
inverse of the sparse block $A_{11}$ which can be computed using
algebraic multigrid. Another approach is to approximate $A_{11}^{-1}$
and to successively refine it by solving the two block equations
of~\eqref{eq:schur}. A pseudocode description of this method is
provided in Algorithm~\ref{alg:schur}.

In general the matrix $A$ will not already be in this sparse/dense
partition form. Reordering can be performed using a simplistic
technique such as sorting by degree of vertices, by reordering schemes
like reverse Cuthill-McKee or approximate minimum degree ordering, or
sorting by number of nonzeros in the rows. For a general graph
Laplacian (corresponding to the normal equations of the first least
squares system), sorting by degree of vertices is a more natural
choice rather than reordering schemes that work well on meshes or
graphs with uniform degree nodes. Likewise, for the second least
squares problem, sorting by number of nonzeros in the rows is a
natural choice for reordering. However, none of these schemes yield
such a sparse/dense partitioning for the graphs we considered. See
Figures~\ref{fig:erlplcln0rrdrng}--\ref{fig:balplcln2zm}
which show the patterns of nonzeros in $\laplacian_0$ and
$\laplacian_2$.

Due to the nature of the nonzero patterns in the matrices and their
reorderings, the computation of $A_{11}^{-1}$ by algebraic multigrid
is not feasible. This was borne out in our experiments and hence the
details are not reported here. The results of using
Algorithm~\ref{alg:schur} for the first least squares problem of
ranking on Erd\H{o}s-R\'enyi random graphs are shown in
Table~\ref{tab:eramg}. While for this problem
Algorithm~\ref{alg:schur} appears slightly more competitive than
algebraic multigrid with smoothed aggregation or Lloyd aggregation, it
is not at all competitive as compared to Krylov solvers. Thus, we did
not apply it to the ranking problem on the other two graph models.

\subsection{Discussion and comparisons}\label{subsec:comparisons}

From the results in Tables~\ref{tab:erkrylv}--\ref{tab:bakrylv}, it is
clear that amongst Krylov methods conjugate gradient is a good choice
for the first least squares problem. For the second least squares
problem, conjugate gradient performs better on sparser graphs but LSQR
does better on denser graphs.

Algebraic multigrid performs optimally for certain types of elliptic
partial differential equations on meshes. So a naive hope would be
that it would do well with Laplacians on graphs. Our experiments
demonstrate that although algebraic multigrid can sometimes solve the
linear systems, it often performs quite poorly in terms of time as
compared to Krylov methods. This is true even if the setup time
required by algebraic multigrid is excluded. (The setup is the process
of forming the coarser levels before solving can begin.) If the setup
time is also included, the performance becomes much worse compared to
Krylov methods. For example, in the cases corresponding to the last
two rows of Table~\ref{tab:eramg}, algebraic multigrid could not even
be used, whereas Krylov methods performed reasonably well (see
Table~\ref{tab:erkrylv}). In the same table, for the $\beta$ problem
in rows 4 and 5, algebraic multigrid took between 2 and over 400
seconds while Krylov methods took only a fraction of a second. Even
not counting the setup phase, these cases took between over 1 to over
32 seconds for algebraic multigrid, as can be seen in
Table~\ref{tab:eramgsttstcs}. Similar behavior is seen for
\BA~graphs. See for instance rows 4, 6, and 9 in
Table~\ref{tab:baamg}.

We are making available the code, and providing details of the
experiments, so that multigrid experts can help improve its
performance. It is quite possible that the experts will find some
small tweaks (which we missed, not being multigrid experts) that make
algebraic multigrid competitive.  In contrast, we note that the Krylov
methods work well ``out-of-the-box'', for example without any
preconditioning.  One should keep in mind however, that our
experiments are for small- to moderate-sized problems, those than can
be done on a single serial computer. Currently that means a linear
system with nearly a quarter million unknowns on a modest
laptop. Moving to larger graph problems, requiring parallel
computation, will open up a whole new field which we have not explored
in this paper.

If the Schur complement method is used with the inverse matrix
$A_{11}^{-1}$ determined by algebraic multigrid then it is slower than
using Krylov methods alone.  This means that the time penalty in this
method is attributable to either the matrix partitioning or to the
algebraic multigrid solution for the inverse matrix, and not to the
Krylov solution of the dense partition. The results of using
Algorithm~\ref{alg:schur}, which uses an approximate inverse, are
shown in Table~\ref{tab:eramg}. It can be seen that for smaller graphs
it is slower that just using algebraic multigrid and for larger graphs
it is marginally faster. However, in either case it is much slower
than using Krylov methods alone.

\begin{figure}[t]
  \centering
  \begin{tabular}{ccc}
    \includegraphics[scale=0.5, trim=1.8in 0.9in 1.8in 0.9in, clip]
    {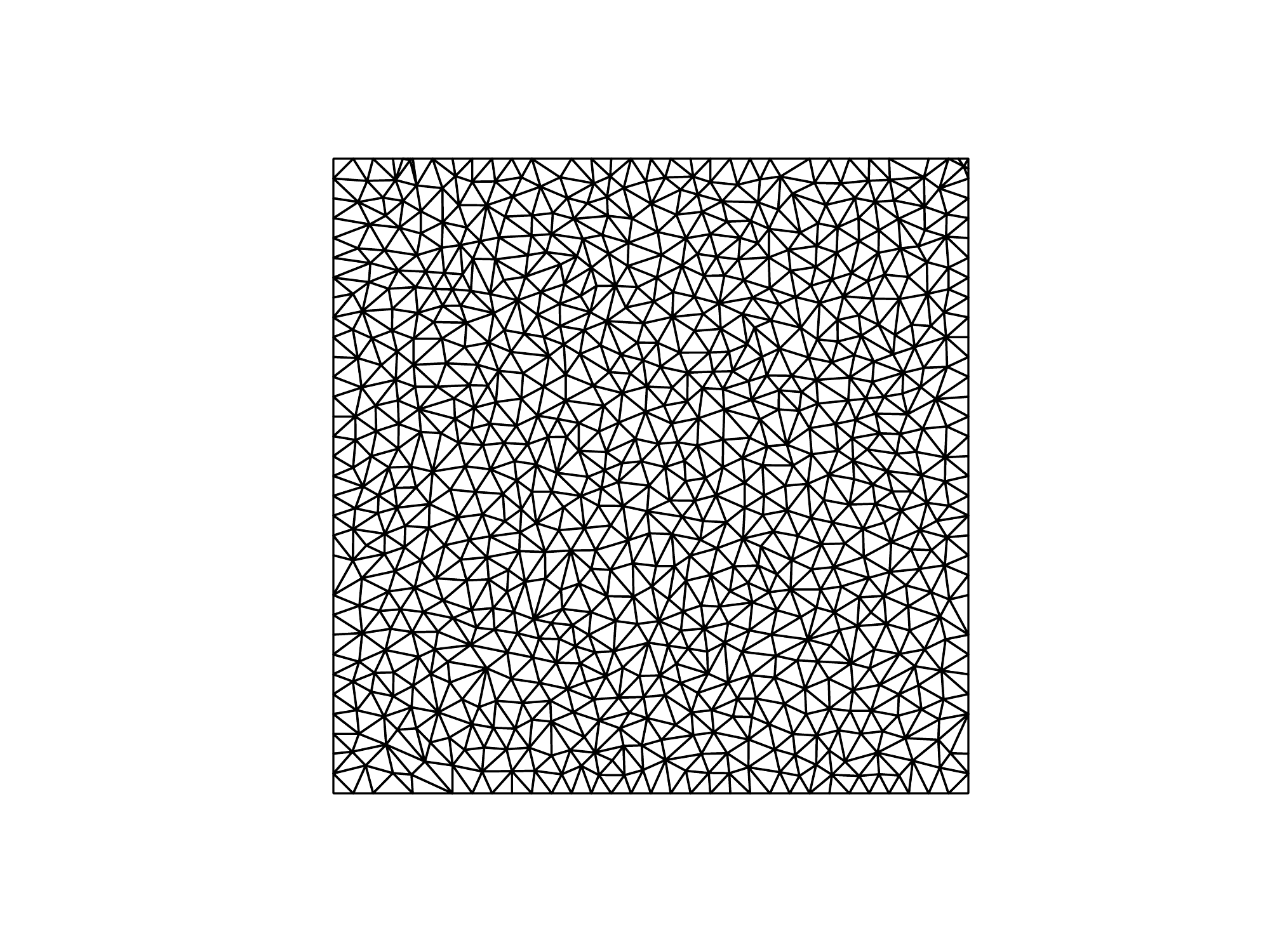} &
    \includegraphics[scale=0.35, trim=1.4in 0.5in 1.4in 0.5in, clip]
    {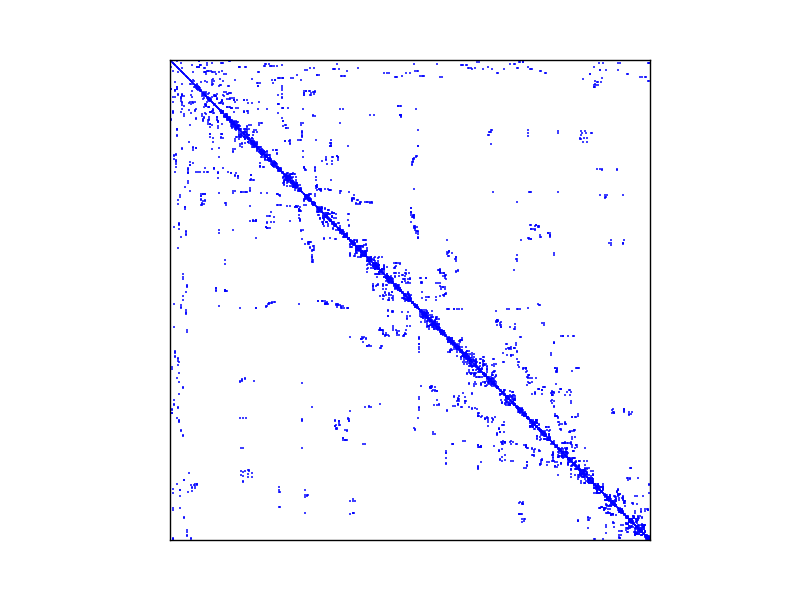} &
    \includegraphics[scale=0.35, trim=1.4in 0.5in 1.4in 0.5in, clip]
    {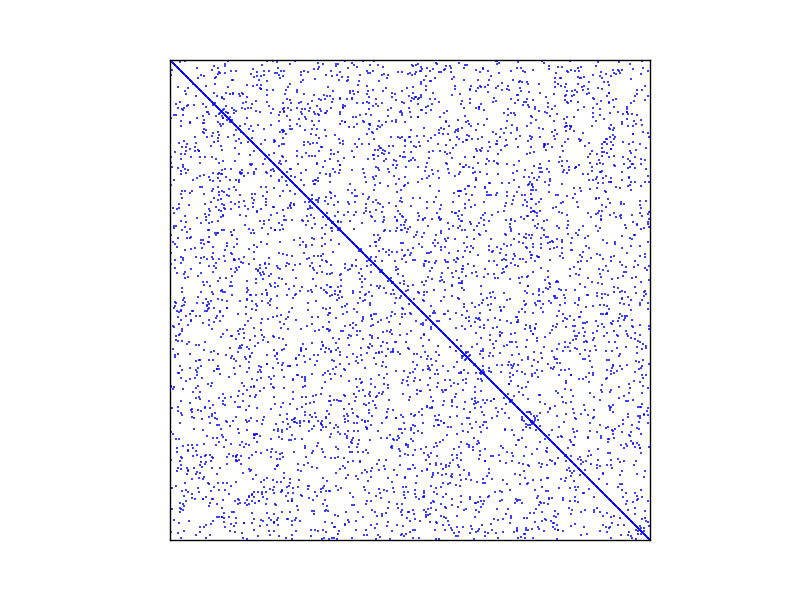}
  \end{tabular}
  \caption{For the mesh shown on the left, the middle plot shows the
    nonzeros in $\laplacian_0$. The rightmost plot shows the nonzeros
    in $\laplacian_0$ for an \ER~graph constructed using the same
    number of vertices and edges as are in the mesh. The
    $\laplacian_0$ on the mesh has locality whereas the one on the
    graph does not.}
  \label{fig:mshvrssgrph}
\end{figure}

Tables~\ref{tab:ernnz}--\ref{tab:bannz} show the number of nonzero
components in the linear systems being solved by the different
methods. The matrix sizes can be determined based on the number of
vertices, edges, and triangles. These tables show the details of all
levels for multigrid methods. For the first least squares problem, the
storage requirements for $\boundary_1$ and $\laplacian_0$ are very
similar, with only an additional $N_0$ nonzeros required in
$\laplacian_0$ to store the vertex degrees.  However, for the second
least squares problem, the storage requirements for $\laplacian_2$
grow at a faster rate than for
$\boundary_2$. Figures~\ref{fig:erlplcln0rrdrng}--\ref{fig:balplcln2zm}
show a visual representation of the nonzeros in most of the Laplacian
matrices. The original matrices as well as various reorderings are
displayed in these figures. These figures indicate that aggregation
schemes for algebraic multigrids may have to be rethought for these
graph problems. There is no apparent (or very little) structure or
locality visible in these figures. This is very different from what
would be observed in the case of Laplacians on meshes. See for example
Figure~\ref{fig:mshvrssgrph}.

\section{Conclusions and Future Work}\label{sec:conclusions}

Other researchers have formulated least squares ranking on graphs and
studied its theory, applications, and connections with other
fields. In this paper, we have presented the first detailed numerical
studies of the two least squares problems involved. We also took this
opportunity to introduce the problem in an elementary way and
highlighted its many connections to different fields. For example, we
have shown that in setting up the least squares ranking problem on
graphs it is natural to make excursions into areas such as elementary
exterior calculus and elementary algebraic topology. In fact, due to
the absence of any geometry ranking provides an easy introduction to
some basic aspects of exterior calculus and algebraic topology. The
lack of geometry arises from the absence of location information about
the vertices and the absence of a metric on the graphs.  Since the
problem and its formulation is so easy to state and work with, the
first author has also used it in undergraduate projects with an aim to
introducing undergraduate students to research. In contrast, the very
closely related problems in elliptic partial differential equations
are much harder to introduce at a comparable elementary level.

The same features which make it easy to state and formulate the
problem also make it easy to implement the linear systems (i.e.,
construct the matrices involved) for these least squares problems. We
have seen that the system matrices only involve the boundary matrices
and their transposes, and these are easily constructed.  Hodge
decomposition and harmonic cochains on graphs are studied in this
paper, but these can also be studied on a simplicial approximation of
a manifold. There these problems arise in the context of numerical
methods for elliptic partial differential equations. The presence of a
metric in that case requires integration and the introduction of a
matrix approximation of the Hodge star. This takes one into the realm
of finite element exterior calculus or discrete exterior calculus. The
creation of matrices (that are analogous to the ones studied in this
paper) requires more work and functional analysis is required to fully
appreciate the numerical issues like stability and convergence to
smooth solution. In graph problems, there is no issue of convergence
to a smooth problem because there is no such problem that we are
approximating.  The simple code developed for the ranking problem is
also easy to use for experimental studies on topology of random clique
complexes. This can be used to formulate conjectures, or as we have
shown, to get more detailed information than the current theory in
that field can provide.

Our detailed numerical studies lead to some interesting
conclusions. Algebraic multigrid method is a topic of intense research
currently since it has been very successful for certain elliptic
partial differential equations. However, it performs poorly in the
graph problems that arise in least squares ranking on graphs. At first
this is surprising. After all, the two problems are Poisson's
equations, one involving the scalar Laplacian and the other involving
the Laplace-deRham operator on 2-cochains. However, as alluded to
earlier and as can be seen from
Figures~\ref{fig:erlplcln0rrdrng}--\ref{fig:balplcln2zm}, the linear
systems in these graph problems suffer from lack of locality and
structure which is usually apparent in partial differential equation
problems on meshes. 

For the first least squares problem, our results indicate that using
the conjugate gradient solver on the normal equations is always the
best option, generating the lowest errors and the fastest times for
each of the three types of graphs we considered.  For the second least
squares problem, conjugate gradient is the better option when the
triangle density (and thus density of $\laplacian_2$) is low, but as the
triangle density increases, LSQR becomes the faster option.

An interesting area for future research is the decomposition of these
problems into smaller pieces which can be solved with minimal
interaction. One then also needs to combine the solutions thus
obtained into a global solution. It is not at all obvious how this can
be done. For partial differential equations, smaller subproblems are
sometimes obtained by domain decomposition methods. Analogously,
appropriate graph partitioning methods and their role in decomposing
these graph problems may be worth studying. These might also be
relevant in creating new aggregation schemes for algebraic
multigrid. The connection between the solvers for diagonally dominant
systems~\cite{KoMiPe2010} and algebraic multigrid would also be
interesting to study further. Once reliable software implementing the
complete algorithm of~\cite{KoMiPe2010} becomes available, comparisons
between that and the various methods studied here would also be
relevant. Graphs on surfaces have a special structure which can be
exploited for efficiency. It would be interesting to study examples of
ranking problems which are natural to study on surfaces. There are
some hints that assignment of ranking work to committees might
naturally fit into this framework.

\phantomsection
\section*{Acknowledgement} The work of ANH and KK was supported in
part by NSF Grant No.~DMS-0645604. We thank Nathan Dunfield, Xiaoye
Jiang, Rich Lehoucq, Luke Olson, Jacob Schroder, and Han Wang for
useful discussions.  

\addcontentsline{toc}{section}{Acknowledgement}

\phantomsection
\bibliographystyle{acmdoi}
\bibliography{hirani,ranking}
\addcontentsline{toc}{section}{References}

\clearpage

\appendix

\section{Topology computations on clique complexes}
\label{appndx:clique}

\begin{figure}[H]
  \centering
  \begin{tabular}{cc}
    \includegraphics[scale=0.4]{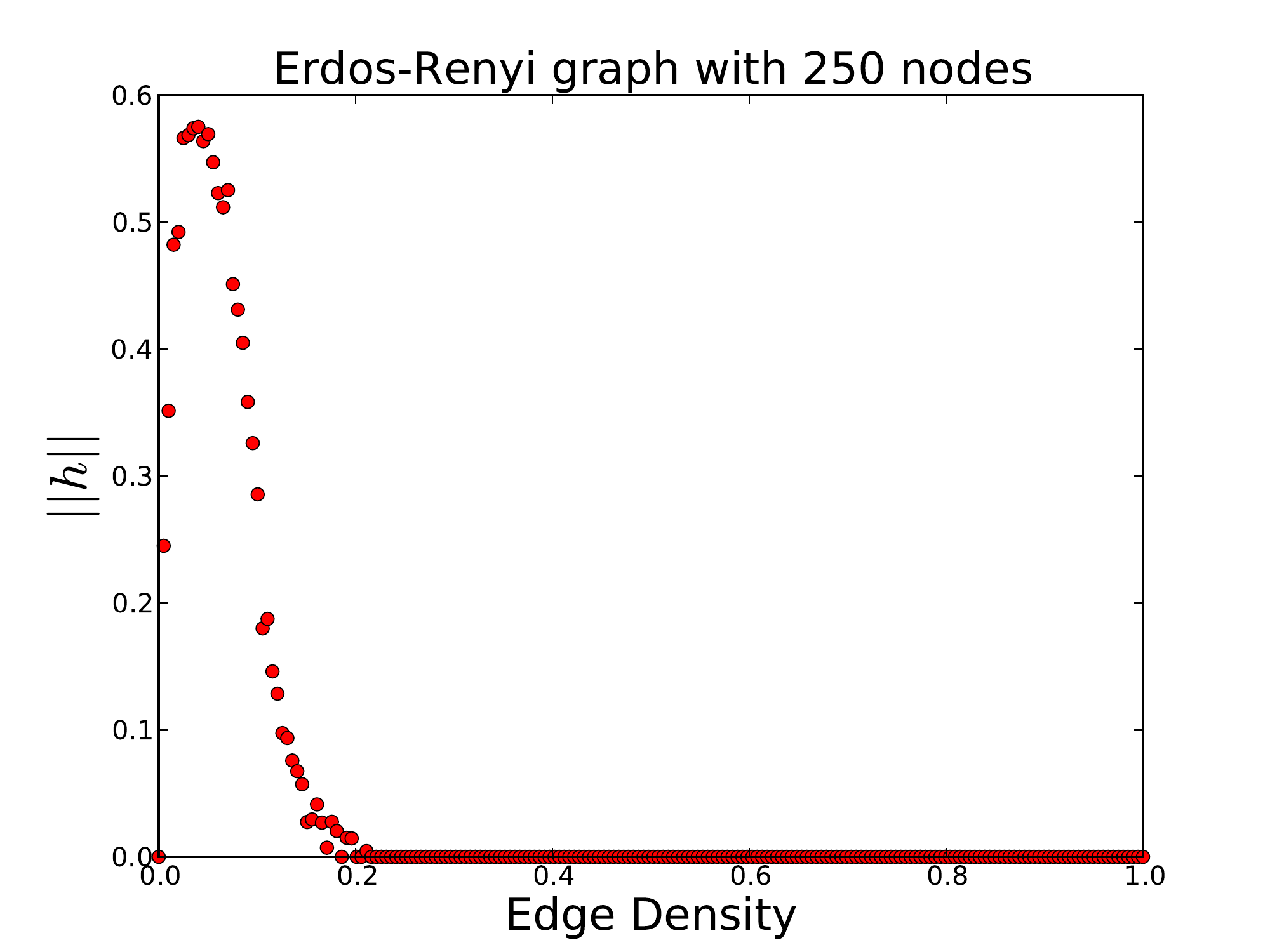} &
    \includegraphics[scale=0.4]{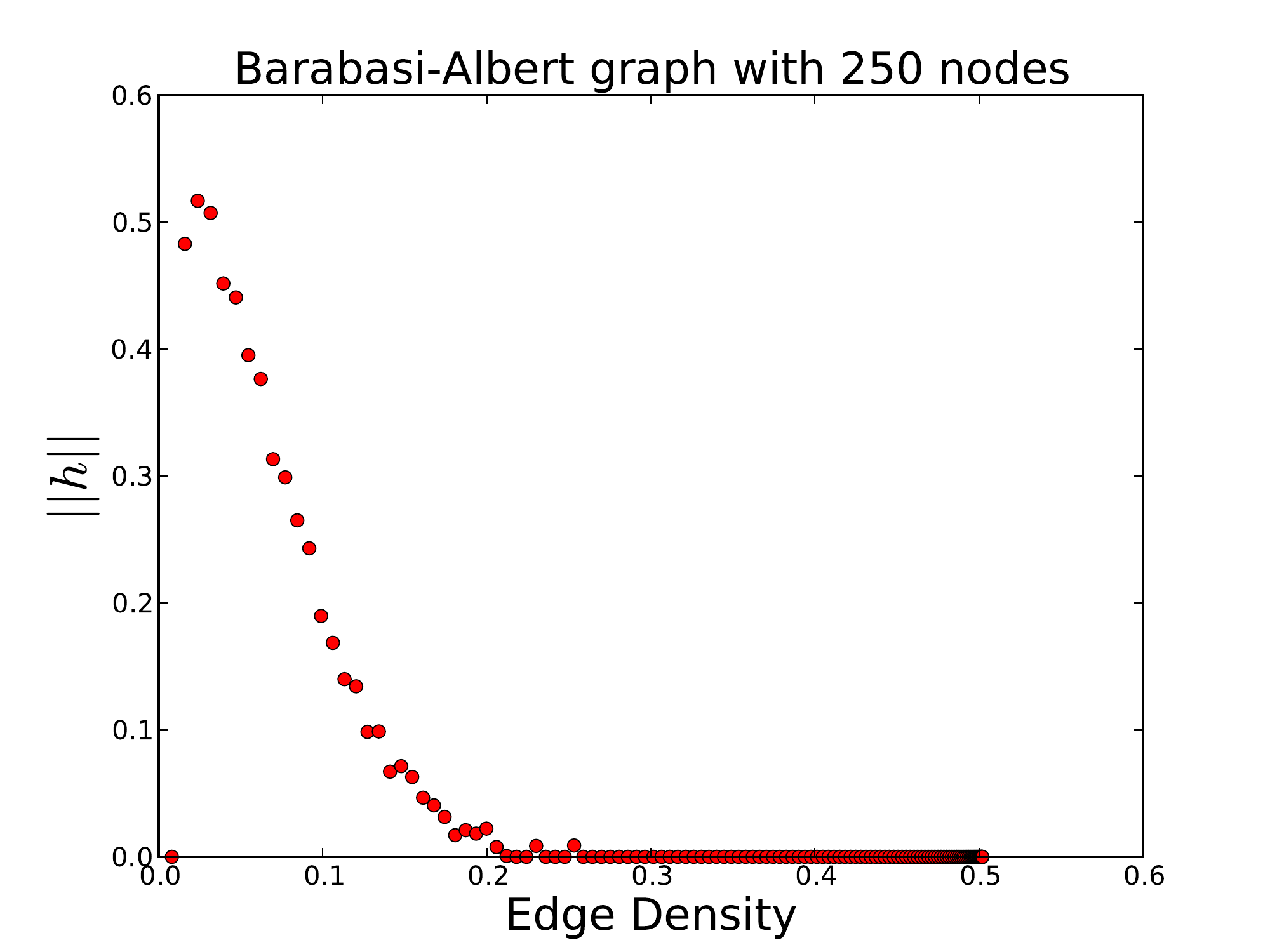} \\
    \includegraphics[scale=0.4]{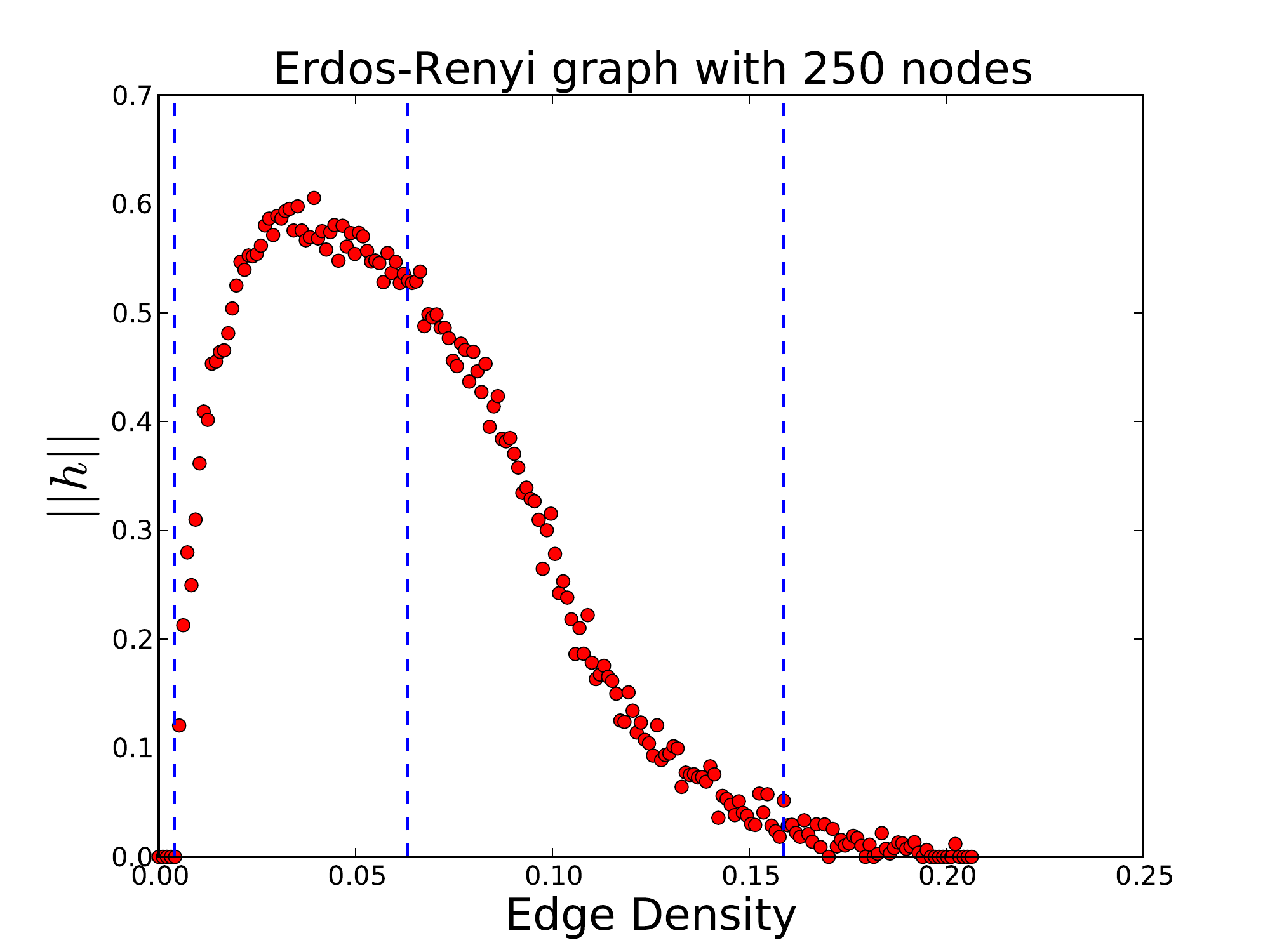} &
    \includegraphics[scale=0.4]{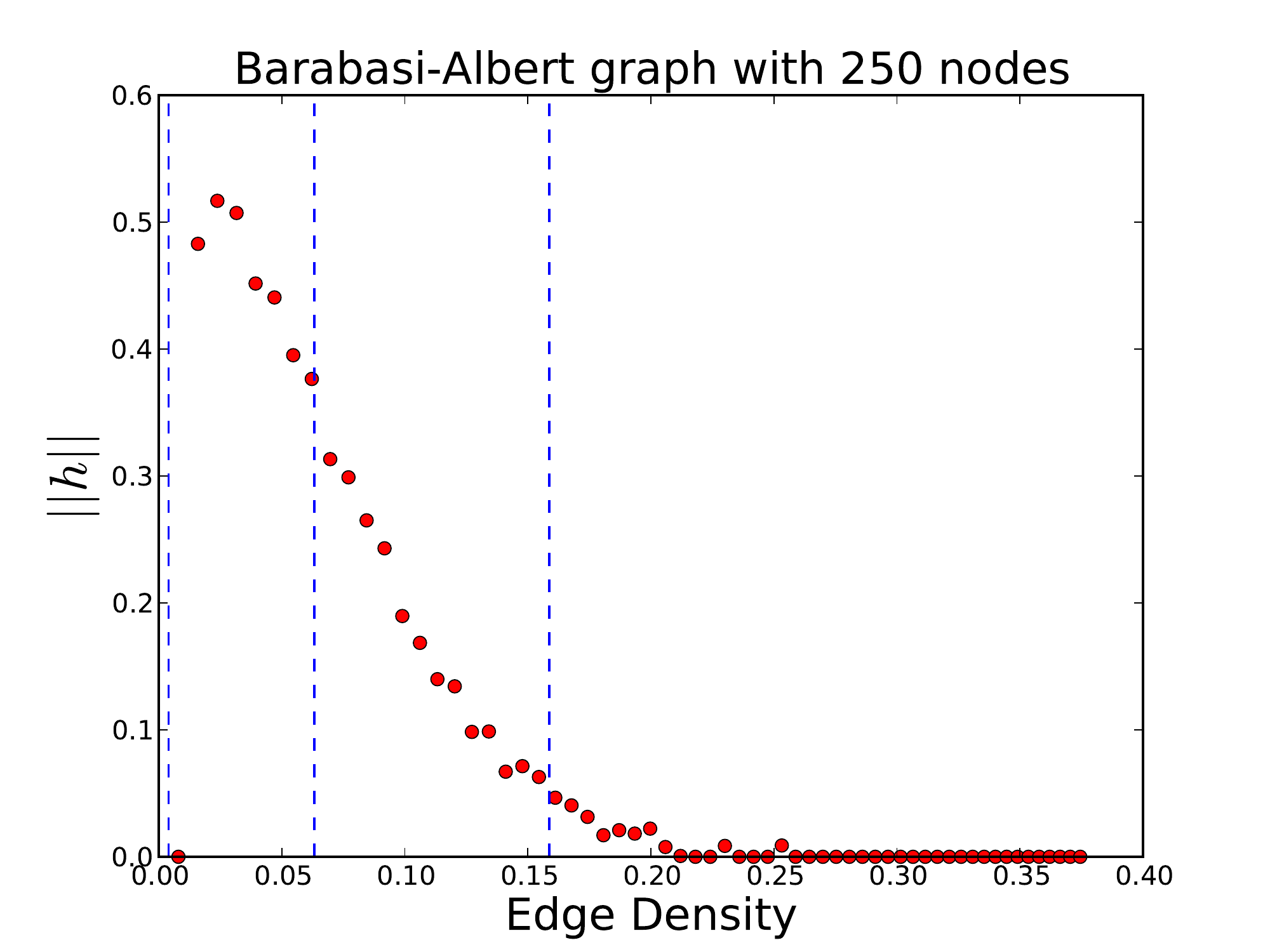} \\
    \includegraphics[scale=0.4]{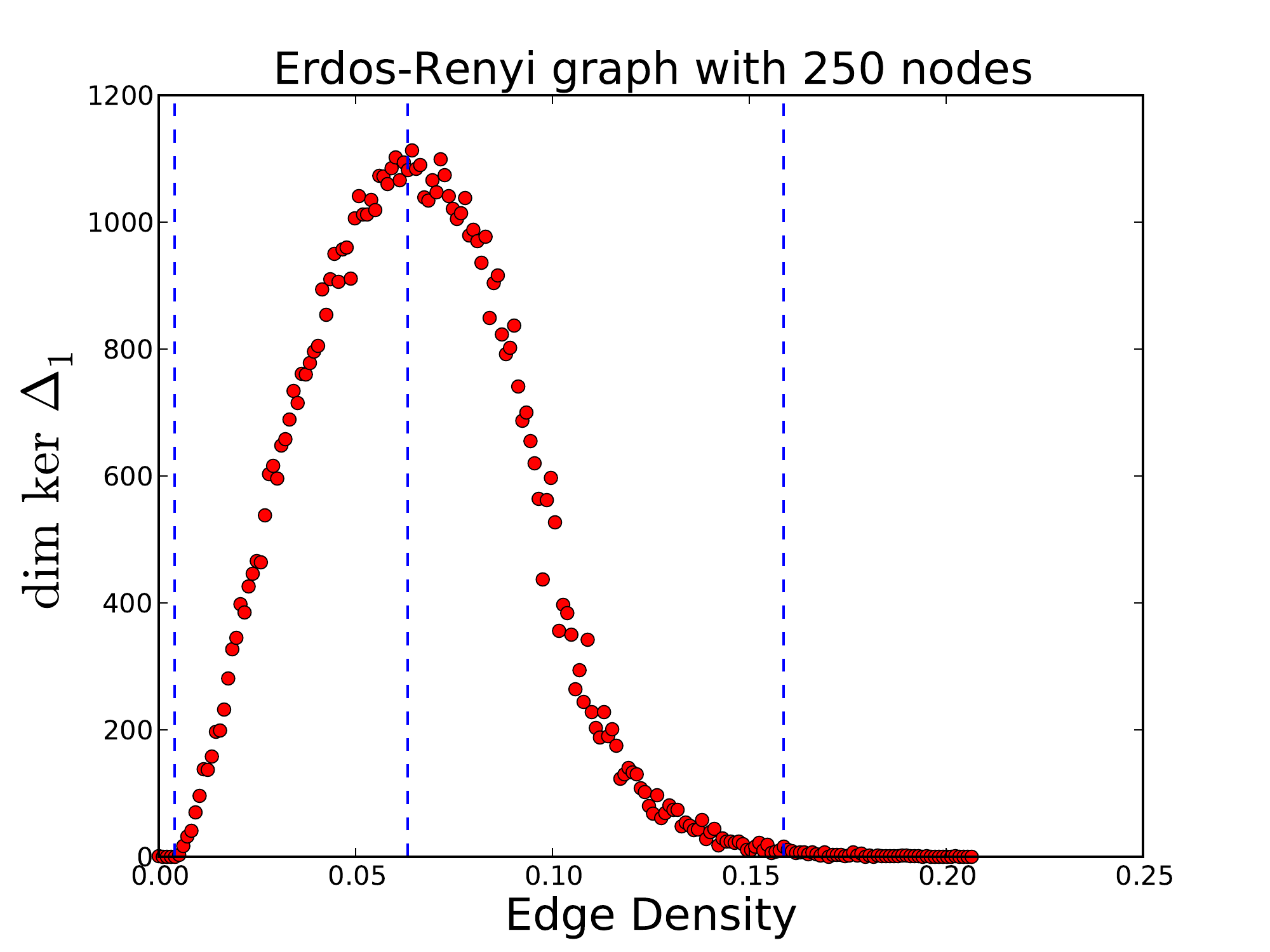} &
    \includegraphics[scale=0.4]{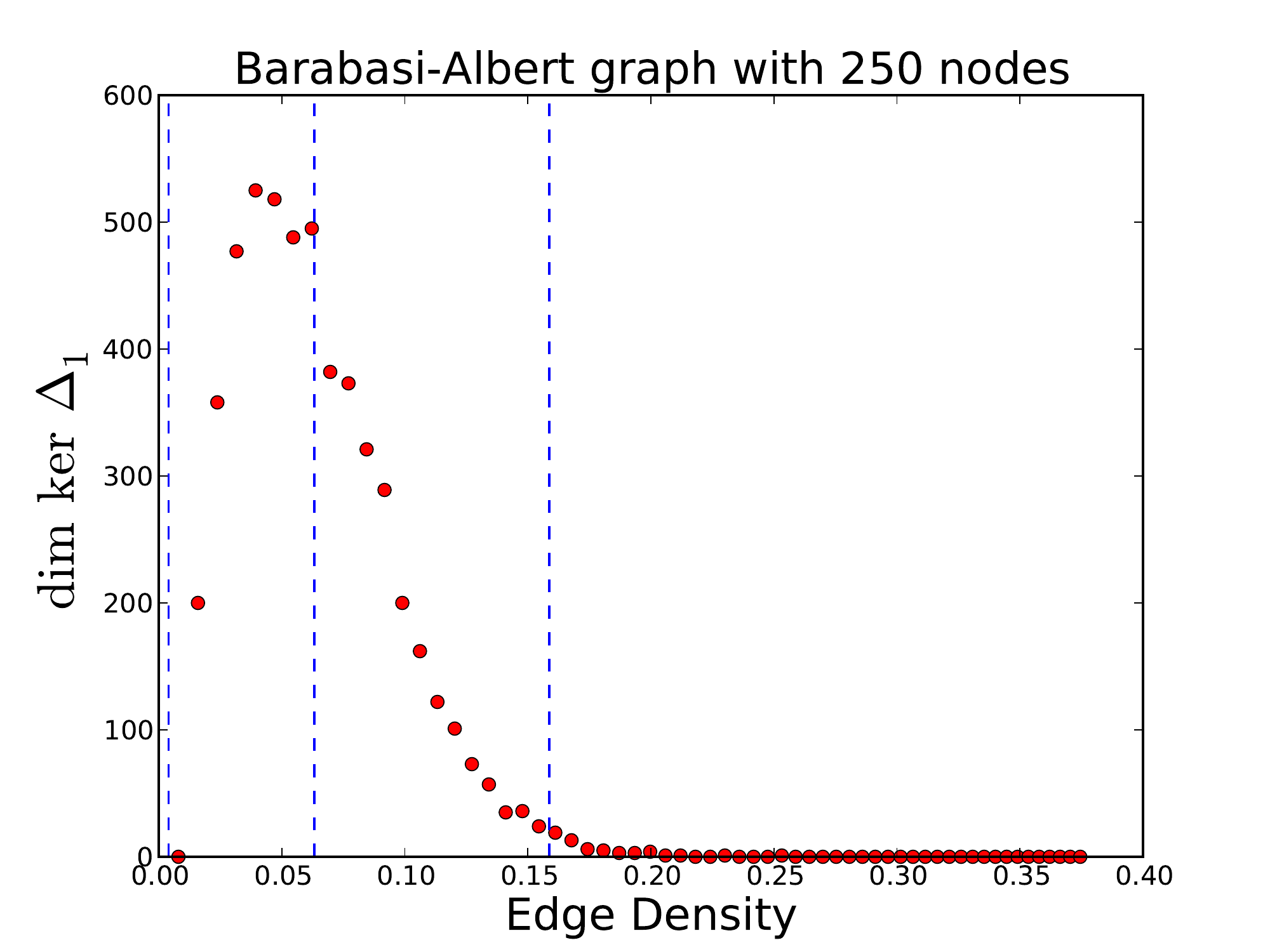}
  \end{tabular}
  \caption{\textit{Left column : }\ER~random graphs; \textit{Right
      column : }\BA~scale-free graphs. The top and middle rows show
    the norm of the harmonic component as a function of edge density.
    The bottom row shows the trends for the 1-dimensional Betti
    number.  The bounds of Kahle~\cite{Kahle2009} are the dashed
    vertical lines. They are also drawn for the right column graphs
    (although the theory is only for \ER~graphs) to provide a
    comparison. See Section~\ref{sec:clique} for more details.}
 \label{fig:clqcmplxs}
\end{figure}

\clearpage

\section{Bounds on conjugate gradient
  iterations} \label{appndx:cgbnds}

\begin{figure}[h]
  \centering
 \includegraphics[scale=0.4]{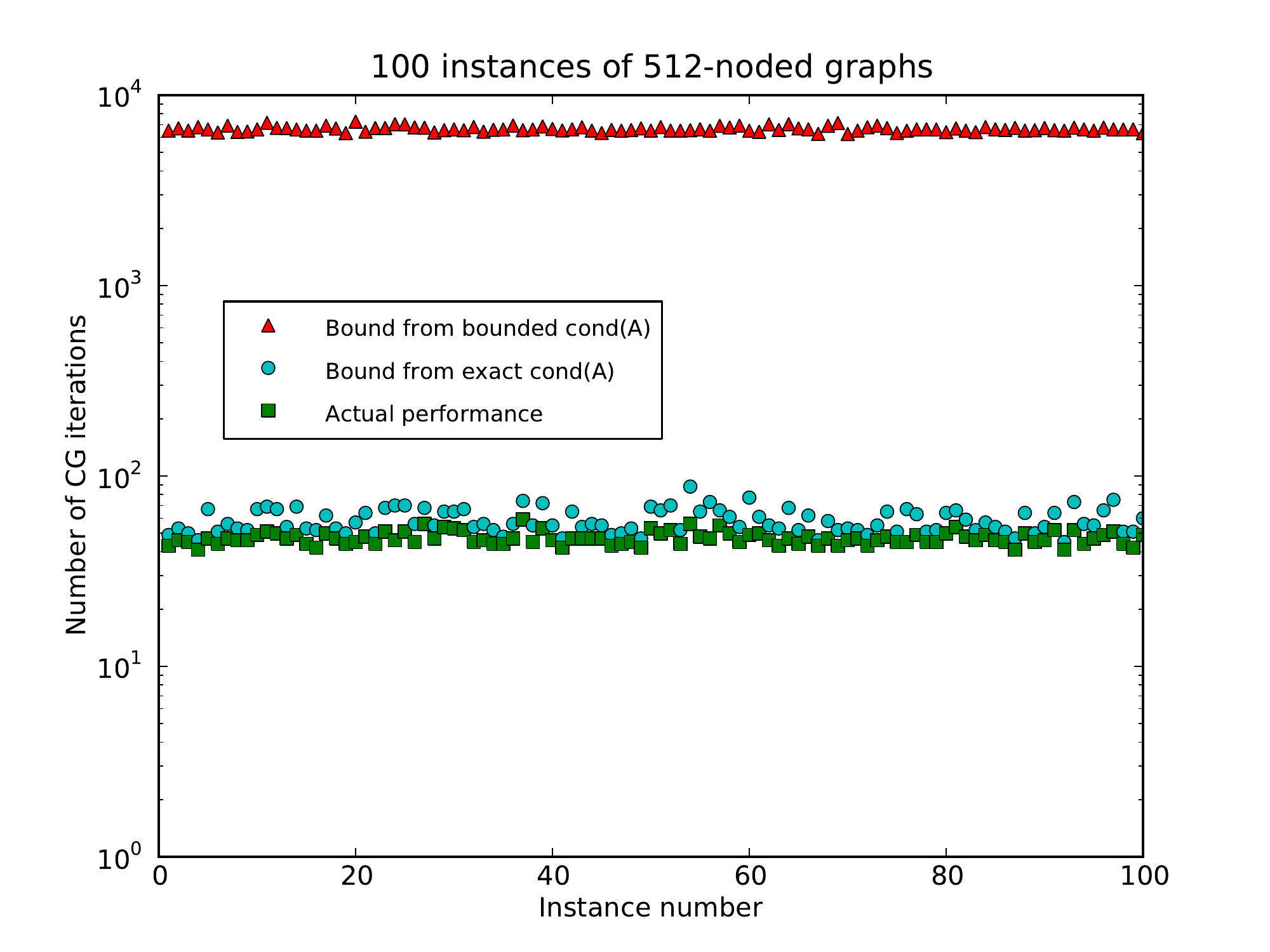}
  \includegraphics[scale=0.4]{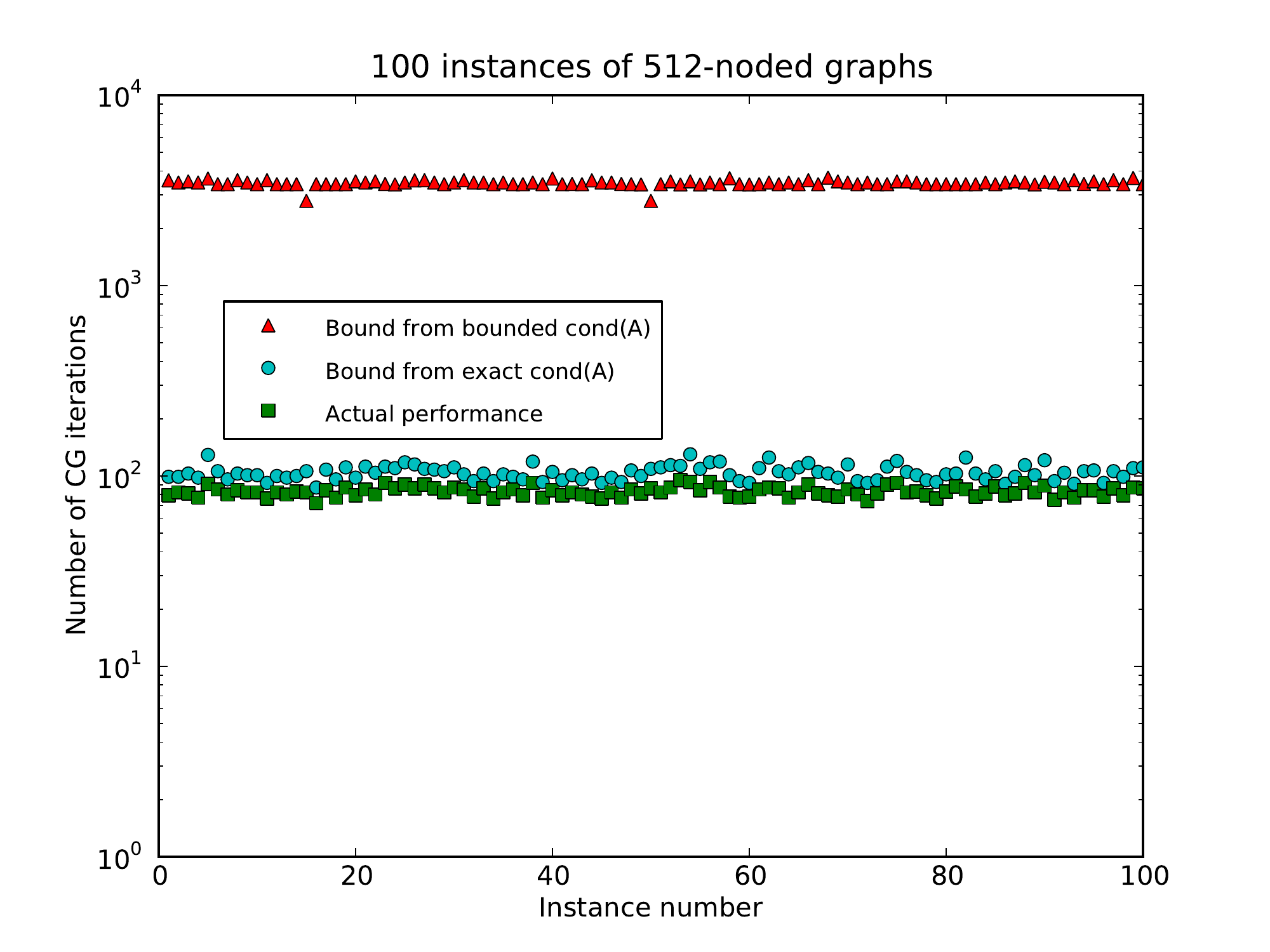}
  \includegraphics[scale=0.4]{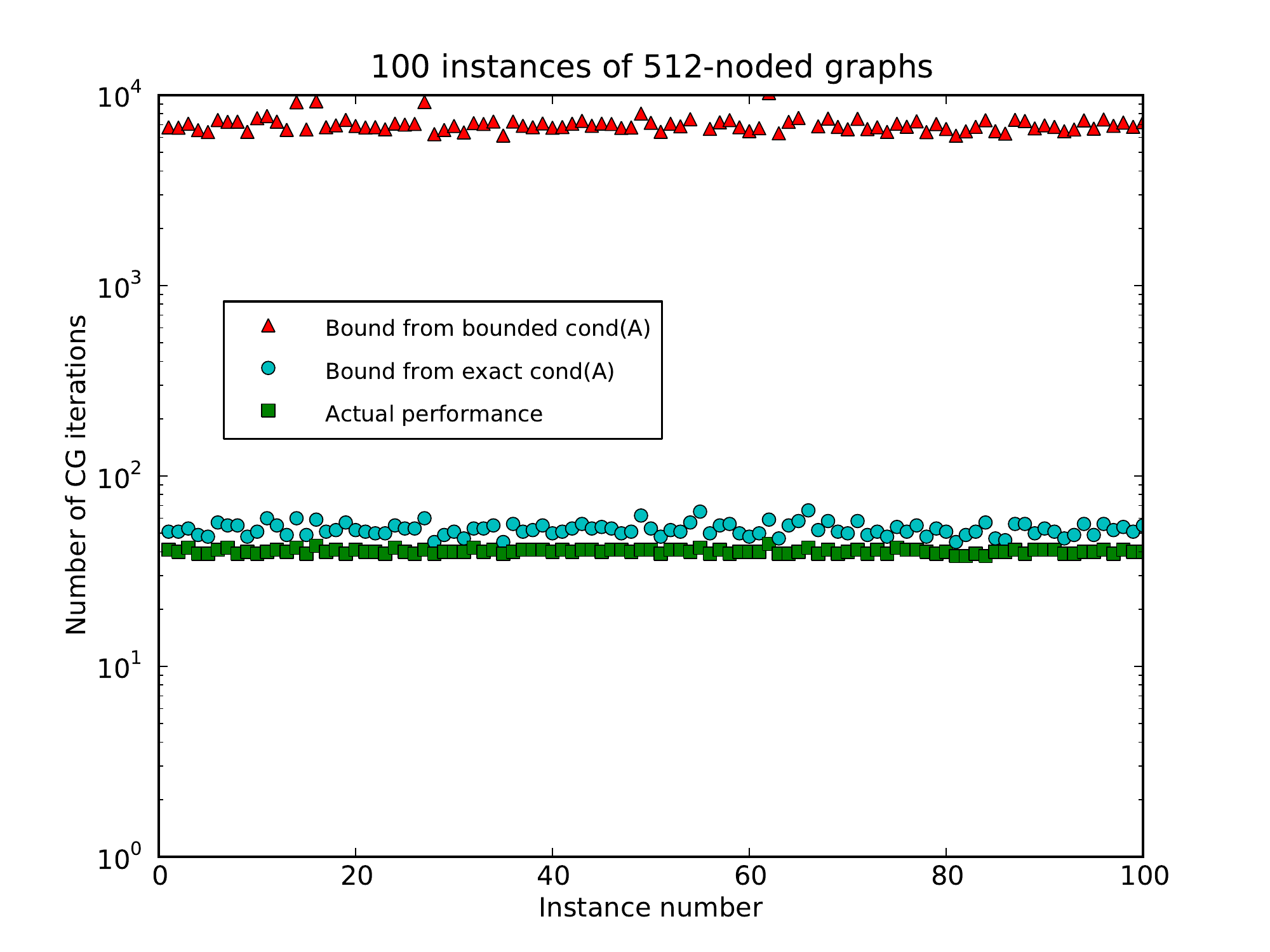}
  \caption{Conjugate gradient iteration bounds compared with actual
    number of iterations required to reduce the error
    $\norm{e_k}_{\laplacian_0}$ to $10^{-8}$ for a variety of
    graphs. \emph{Top row :} results for
    Erd\H{o}s-R\'enyi random graphs and Watts-Strogatz graphs.
    \emph{Bottom row :} results for Barab\'asi-Albert graphs.}
  \label{fig:cgbounds2}
\end{figure}

\begin{figure}[h]
  \centering
  \includegraphics[scale=0.4]{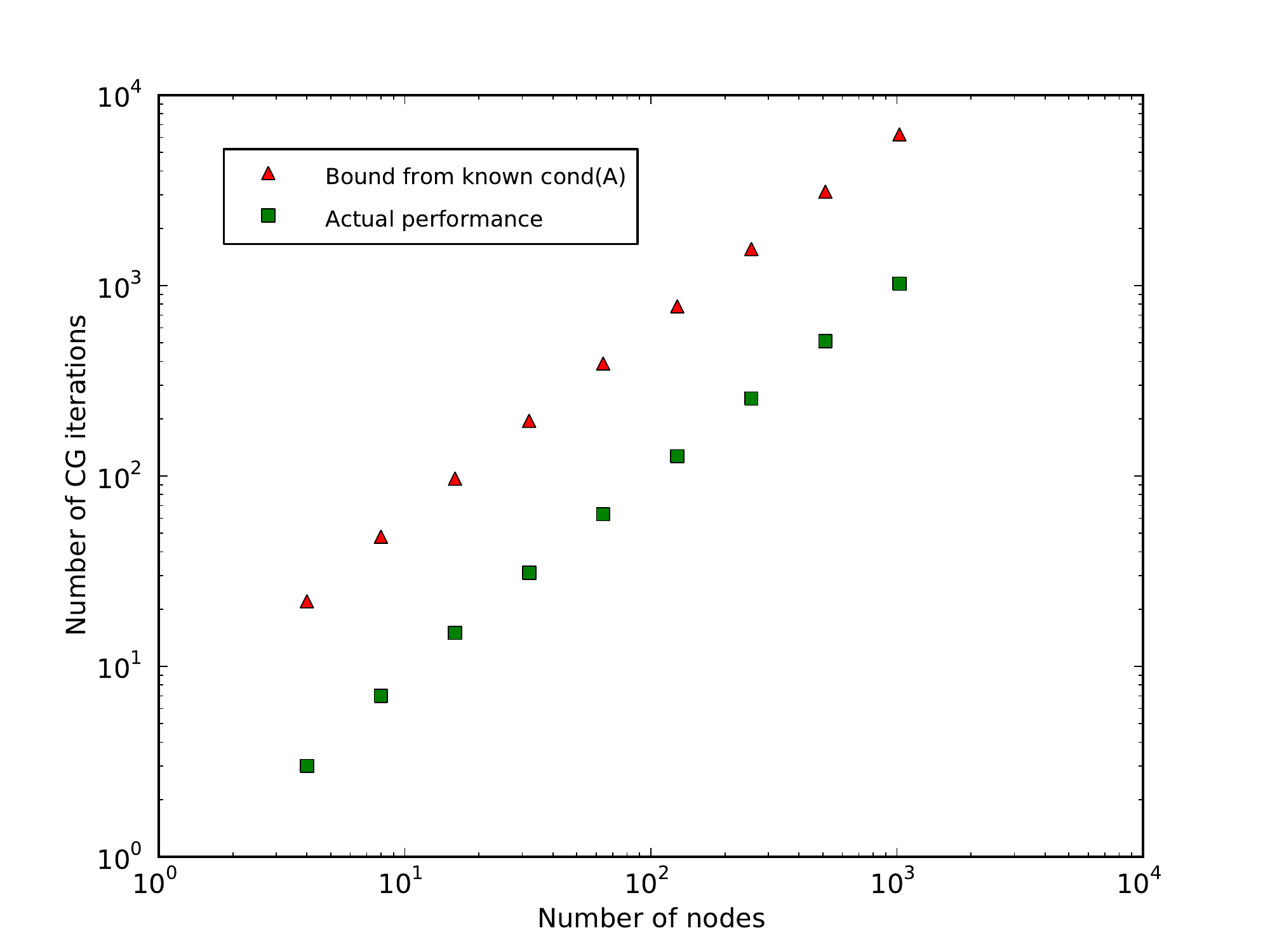}
  \includegraphics[scale=0.4]{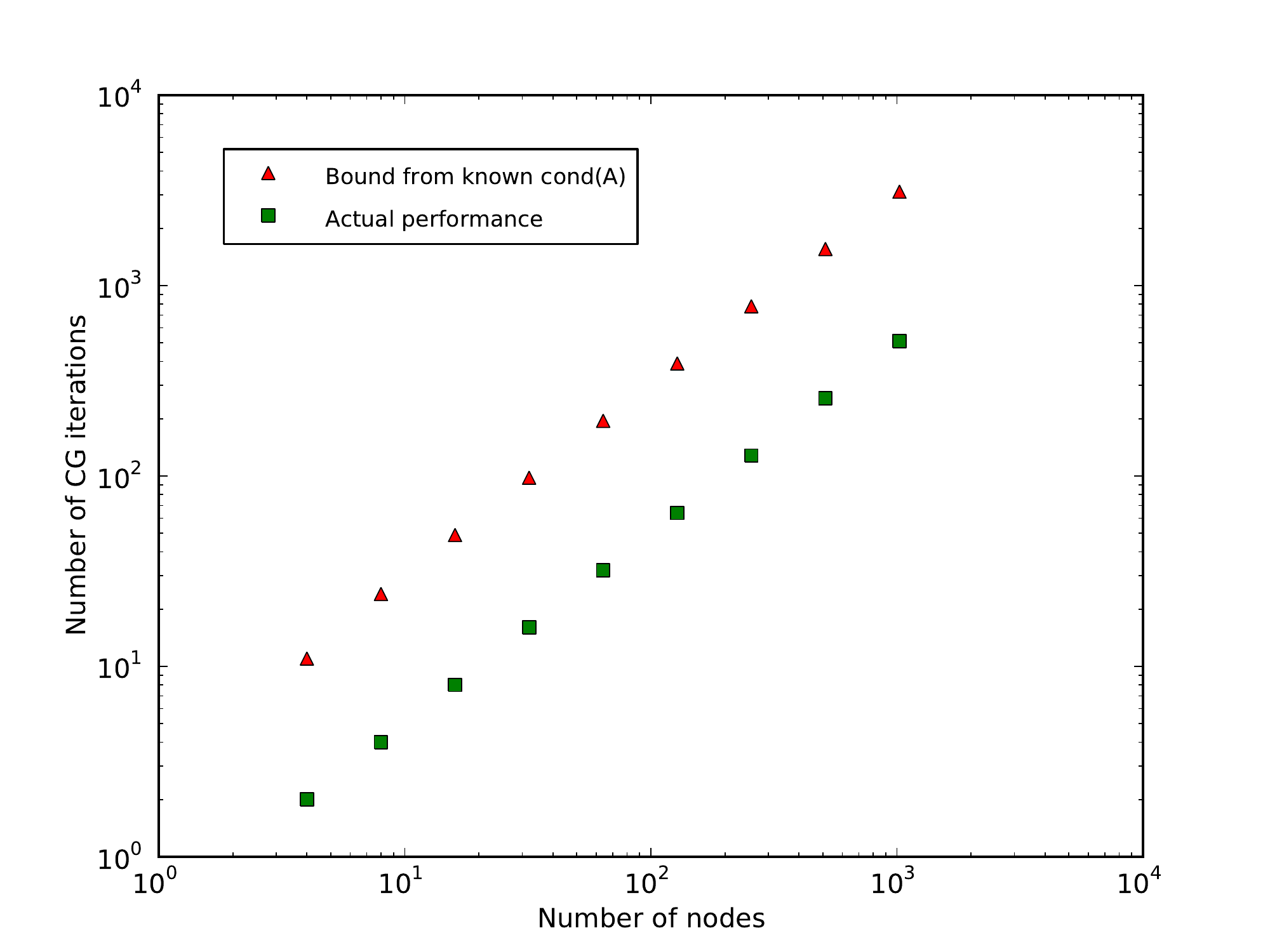}
  \includegraphics[scale=0.4]{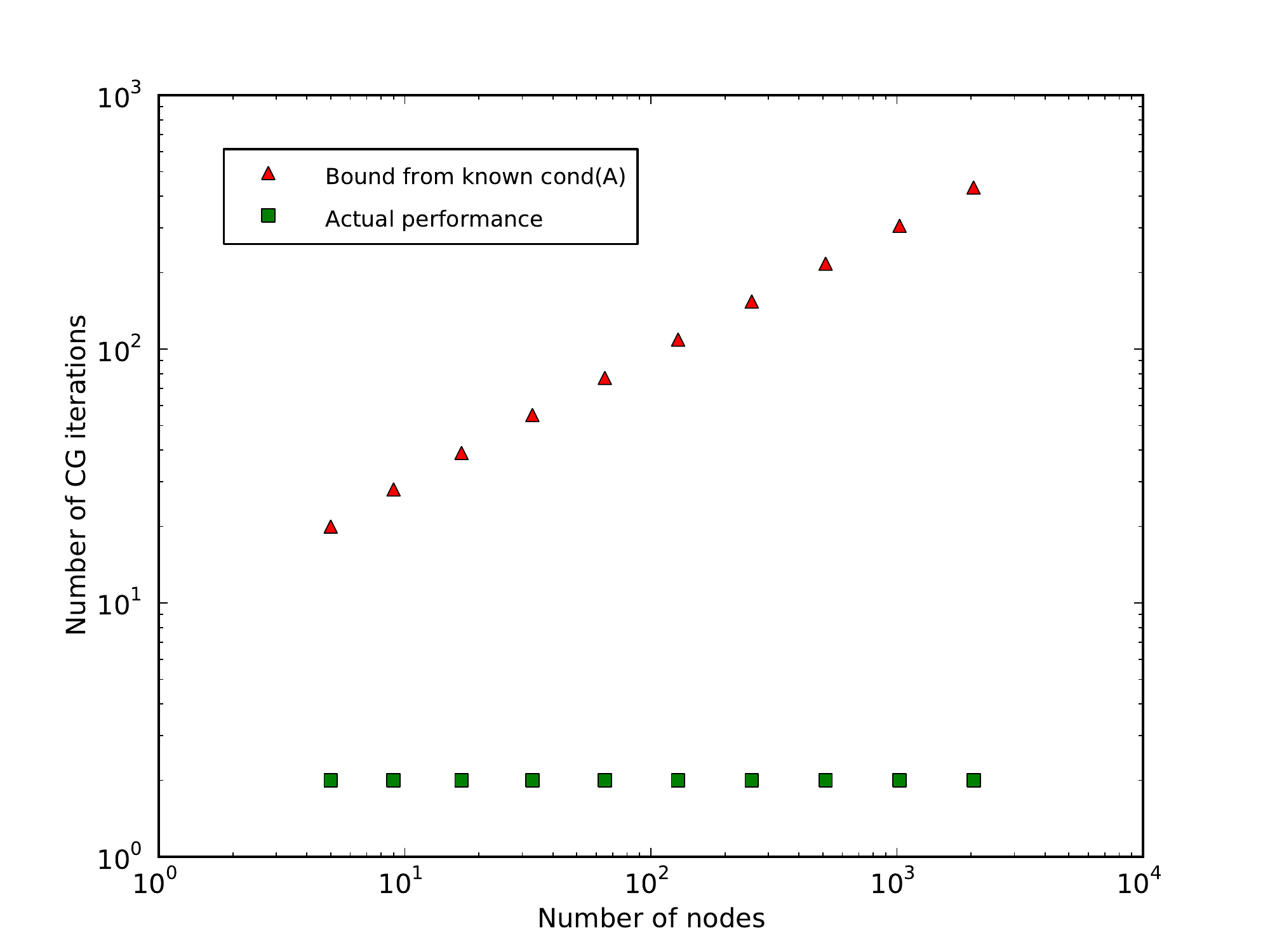}
  \includegraphics[scale=0.4]{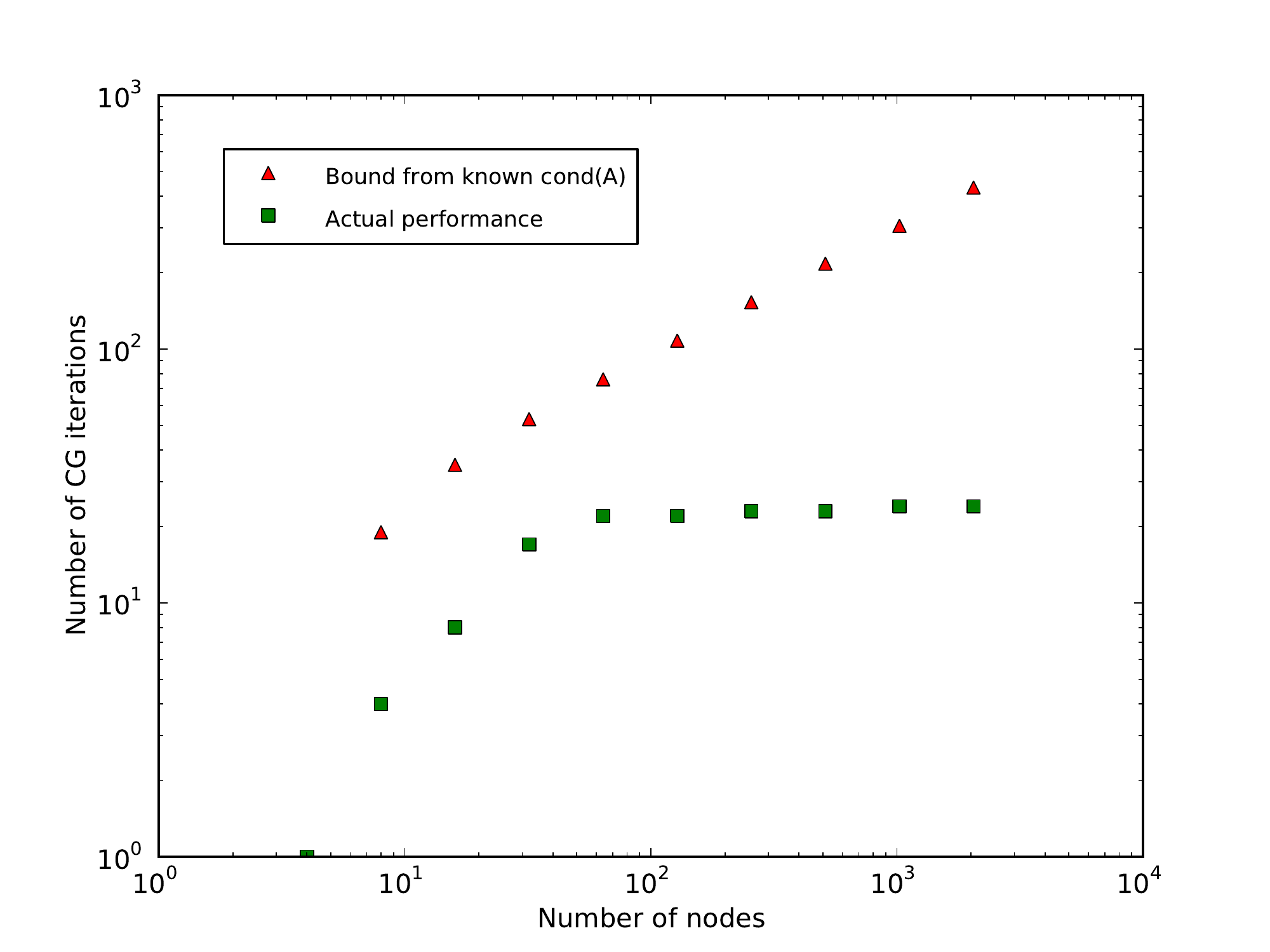}
  \includegraphics[scale=0.4]{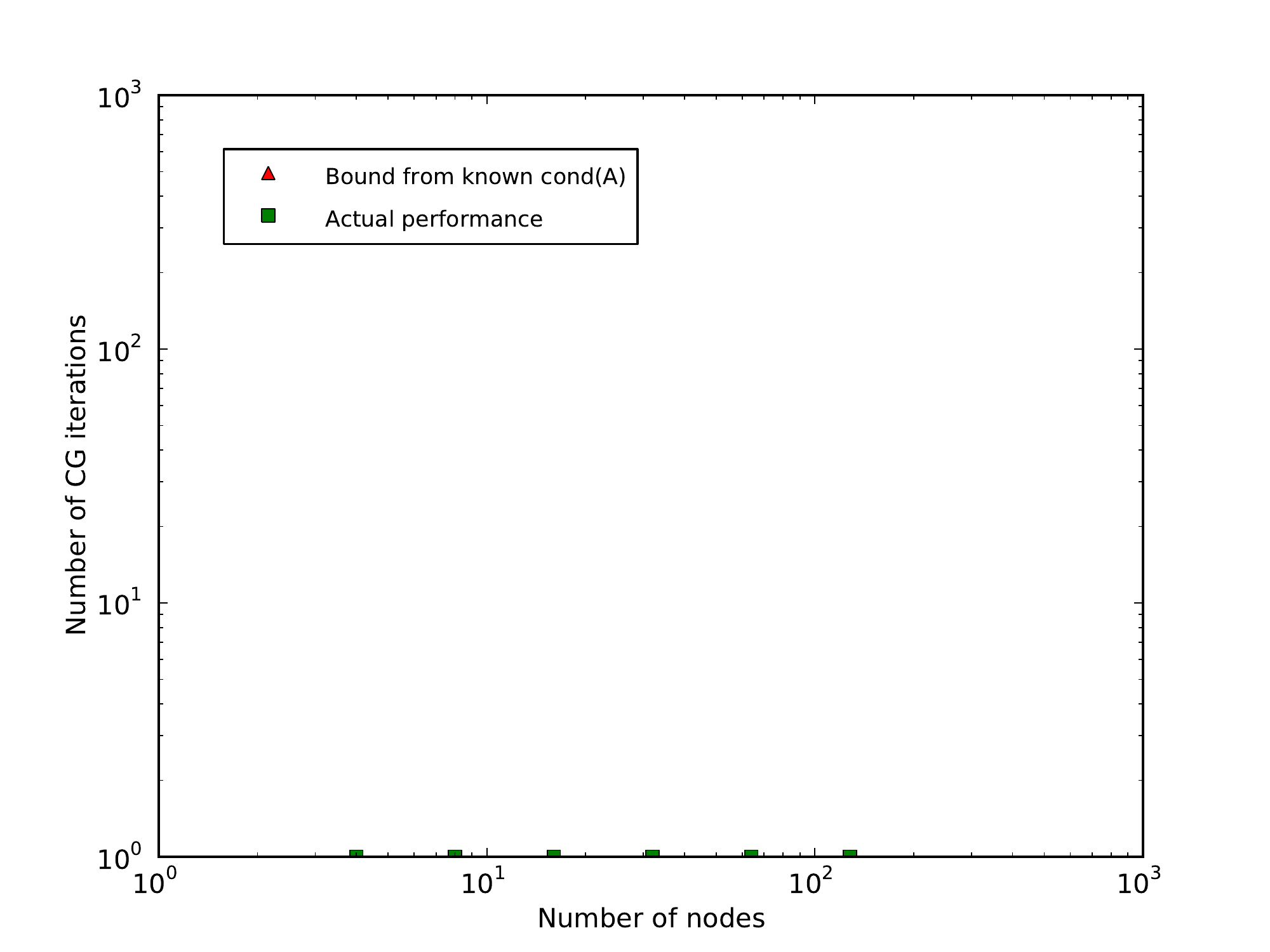}
  \caption{Conjugate gradient iteration bounds compared with actual
    number of iterations required to reduce the error
    $\norm{e_k}_{\laplacian_0}$ to $10^{-8}$ for a variety of
    graphs. \emph{Top row :} results for path and cycle graphs.
    \emph{Middle row :} results for star and wheel
    graphs. \emph{Bottom row :} results for complete graphs. In the
    complete graph case, the markers for bounds are not visible in the
    plot because they coincide with the other markers.  }
 \label{fig:cgbounds1}
\end{figure}

Figures~\ref{fig:cgbounds2} and~\ref{fig:cgbounds1} show the conjugate
gradient iteration estimates and actual number of iterations in
experiments. The condition number of the Laplacian used is modulo the
kernel, i.e., it is the ratio of the magnitude of the largest
eigenvalue to the magnitude of the smallest nonzero eigenvalue.  For
certain special graphs, the condition number is explicitly known and
the iteration count bound obtained from these is labeled ``Bound from
known cond($A$)'' in the figures. For \ER, \WS, and \BA~graphs there
are known bounds on the condition number and the iteration count
obtained from such bounds is denoted ``Bound from bounded
cond($A$)''. For such graphs we also measured the condition number
(modulo kernel) to obtain iteration bounds that are labeled ``Bound
from exact cond($A$)''.

\clearpage

\section{Timing and error using Krylov solvers} \label{appndx:krylov}
The next three tables show the results of numerical experiments on
simplicial complexes \ER, \WS, and \BA~graphs.  The cases where the
homology of the simplicial 2-complex is trivial are marked by an
asterisk ($\ast$). In these cases it is the absolute error which is
shown. See Section~\ref{subsec:krylov} for details.
\bigskip
\begin{table}[h] \scriptsize
 \centering
 \begin{tabular}{|ccccc|c|ccc|cc|cc|}
\hline
\multirow{2}{*}{$N_0$} & \multirow{2}{*}{$N_1$} & \multirow{2}{*}{$N_2$} & Edge & Triangle & Algorithm / & \multicolumn{3}{|c|}{Relative Error} & \multicolumn{2}{|c|}{$\alpha$} & \multicolumn{2}{|c|}{$\beta$} \\
  &  &  & Density & Density & Formulation & $\|\boundary_1^T \alpha\|$ & $\|\boundary_2\beta\|$ & $\|h\|$ & iter. & sec. & iter. & sec. \\ 
 \hline 
\multirow{5}{*}{100} & \multirow{5}{*}{380} & \multirow{5}{*}{52} & \multirow{5}{*}{7.68e-02} & \multirow{5}{*}{3.22e-04} & \textcolor{Red}{CG} & \textcolor{Red}{1.1e-08} & \textcolor{Red}{6.3e-09} & \textcolor{Red}{2.0e-07} & \textcolor{Red}{27} & \textcolor{Red}{0.0041} & \textcolor{Red}{17} & \textcolor{Red}{0.0027} \\ 
  &  &  &  &  & \textcolor{Cerulean}{MINRES} & \textcolor{Cerulean}{3.0e-08} & \textcolor{Cerulean}{1.3e-08} & \textcolor{Cerulean}{4.8e-07} & \textcolor{Cerulean}{26} & \textcolor{Cerulean}{0.0084} & \textcolor{Cerulean}{16} & \textcolor{Cerulean}{0.0053} \\ 
  &  &  &  &  & \textcolor{OliveGreen}{CG-K} & \textcolor{OliveGreen}{1.0e-08} & \textcolor{OliveGreen}{6.3e-09} & \textcolor{OliveGreen}{1.3e-07} & \textcolor{OliveGreen}{57} & \textcolor{OliveGreen}{0.0090} & \textcolor{OliveGreen}{35} & \textcolor{OliveGreen}{0.0055} \\ 
  &  &  &  &  & \textcolor{Violet}{MINRES-K} & \textcolor{Violet}{3.0e-08} & \textcolor{Violet}{1.3e-08} & \textcolor{Violet}{4.8e-07} & \textcolor{Violet}{53} & \textcolor{Violet}{0.0174} & \textcolor{Violet}{33} & \textcolor{Violet}{0.0108} \\ 
  &  &  &  &  & \textcolor{Magenta}{LSQR} & \textcolor{Magenta}{2.8e-08} & \textcolor{Magenta}{1.3e-08} & \textcolor{Magenta}{4.7e-07} & \textcolor{Magenta}{26} & \textcolor{Magenta}{0.0166} & \textcolor{Magenta}{16} & \textcolor{Magenta}{0.0103} \\ 
\hline 
\multirow{5}{*}{100} & \multirow{5}{*}{494} & \multirow{5}{*}{144} & \multirow{5}{*}{9.98e-02} & \multirow{5}{*}{8.91e-04} & \textcolor{Red}{CG} & \textcolor{Red}{1.1e-08} & \textcolor{Red}{1.0e-08} & \textcolor{Red}{2.9e-07} & \textcolor{Red}{21} & \textcolor{Red}{0.0037} & \textcolor{Red}{30} & \textcolor{Red}{0.0053} \\ 
  &  &  &  &  & \textcolor{Cerulean}{MINRES} & \textcolor{Cerulean}{5.2e-08} & \textcolor{Cerulean}{1.1e-07} & \textcolor{Cerulean}{2.8e-06} & \textcolor{Cerulean}{19} & \textcolor{Cerulean}{0.0066} & \textcolor{Cerulean}{27} & \textcolor{Cerulean}{0.0095} \\ 
  &  &  &  &  & \textcolor{OliveGreen}{CG-K} & \textcolor{OliveGreen}{4.8e-09} & \textcolor{OliveGreen}{6.5e-09} & \textcolor{OliveGreen}{1.5e-07} & \textcolor{OliveGreen}{45} & \textcolor{OliveGreen}{0.0082} & \textcolor{OliveGreen}{438} & \textcolor{OliveGreen}{0.0788} \\ 
  &  &  &  &  & \textcolor{Violet}{MINRES-K} & \textcolor{Violet}{5.2e-08} & \textcolor{Violet}{5.3e-08} & \textcolor{Violet}{1.5e-06} & \textcolor{Violet}{39} & \textcolor{Violet}{0.0137} & \textcolor{Violet}{57} & \textcolor{Violet}{0.0201} \\ 
  &  &  &  &  & \textcolor{Magenta}{LSQR} & \textcolor{Magenta}{4.9e-08} & \textcolor{Magenta}{2.0e-08} & \textcolor{Magenta}{8.9e-07} & \textcolor{Magenta}{19} & \textcolor{Magenta}{0.0121} & \textcolor{Magenta}{29} & \textcolor{Magenta}{0.0183} \\ 
\hline 
\multirow{5}{*}{100} & \multirow{5}{*}{1212} & \multirow{5}{*}{2359} & \multirow{5}{*}{2.45e-01} & \multirow{5}{*}{1.46e-02} & \textcolor{Red}{CG} & \textcolor{Red}{3.9e-09} & \textcolor{Red}{1.6e-08} & \textcolor{Red}{3.5e-05} & \textcolor{Red}{14} & \textcolor{Red}{0.0023} & \textcolor{Red}{60} & \textcolor{Red}{0.0191} \\ 
  &  &  &  &  & \textcolor{Cerulean}{MINRES} & \textcolor{Cerulean}{1.9e-08} & \textcolor{Cerulean}{1.1e-06} & \textcolor{Cerulean}{2.5e-03} & \textcolor{Cerulean}{13} & \textcolor{Cerulean}{0.0044} & \textcolor{Cerulean}{46} & \textcolor{Cerulean}{0.0229} \\ 
  &  &  &  &  & \textcolor{OliveGreen}{CG-K} & \textcolor{OliveGreen}{3.9e-09} & \textcolor{OliveGreen}{4.6e-09} & \textcolor{OliveGreen}{1.0e-05} & \textcolor{OliveGreen}{29} & \textcolor{OliveGreen}{0.0053} & \textcolor{OliveGreen}{129} & \textcolor{OliveGreen}{0.0342} \\ 
  &  &  &  &  & \textcolor{Violet}{MINRES-K} & \textcolor{Violet}{3.0e-07} & \textcolor{Violet}{2.7e-07} & \textcolor{Violet}{6.1e-04} & \textcolor{Violet}{23} & \textcolor{Violet}{0.0084} & \textcolor{Violet}{105} & \textcolor{Violet}{0.0466} \\ 
  &  &  &  &  & \textcolor{Magenta}{LSQR} & \textcolor{Magenta}{7.6e-08} & \textcolor{Magenta}{3.0e-08} & \textcolor{Magenta}{7.6e-05} & \textcolor{Magenta}{12} & \textcolor{Magenta}{0.0082} & \textcolor{Magenta}{58} & \textcolor{Magenta}{0.0429} \\ 
\hline 
\multirow{5}{*}{100} & \multirow{5}{*}{2530} & \multirow{5}{*}{21494} & \multirow{5}{*}{5.11e-01} & \multirow{5}{*}{1.33e-01} & \textcolor{Red}{CG} & \textcolor{Red}{2.5e-09} & \textcolor{Red}{9.2e-09} & \textcolor{Red}{3.0e-06$^{\ast}$} & \textcolor{Red}{10} & \textcolor{Red}{0.0017} & \textcolor{Red}{17} & \textcolor{Red}{0.1033} \\ 
  &  &  &  &  & \textcolor{Cerulean}{MINRES} & \textcolor{Cerulean}{2.0e-08} & \textcolor{Cerulean}{3.1e-07} & \textcolor{Cerulean}{1.0e-04$^{\ast}$} & \textcolor{Cerulean}{9} & \textcolor{Cerulean}{0.0032} & \textcolor{Cerulean}{14} & \textcolor{Cerulean}{0.0897} \\ 
  &  &  &  &  & \textcolor{OliveGreen}{CG-K} & \textcolor{OliveGreen}{1.9e-08} & \textcolor{OliveGreen}{9.7e-09} & \textcolor{OliveGreen}{2.1e-14$^{\ast}$} & \textcolor{OliveGreen}{19} & \textcolor{OliveGreen}{0.0044} & \textcolor{OliveGreen}{34} & \textcolor{OliveGreen}{0.0438} \\ 
  &  &  &  &  & \textcolor{Violet}{MINRES-K} & \textcolor{Violet}{1.0e-06} & \textcolor{Violet}{2.9e-07} & \textcolor{Violet}{2.1e-06$^{\ast}$} & \textcolor{Violet}{15} & \textcolor{Violet}{0.0062} & \textcolor{Violet}{28} & \textcolor{Violet}{0.0476} \\ 
  &  &  &  &  & \textcolor{Magenta}{LSQR} & \textcolor{Magenta}{1.3e-07} & \textcolor{Magenta}{2.9e-09} & \textcolor{Magenta}{2.9e-06$^{\ast}$} & \textcolor{Magenta}{8} & \textcolor{Magenta}{0.0059} & \textcolor{Magenta}{18} & \textcolor{Magenta}{0.0362} \\ 
\hline 
\multirow{5}{*}{100} & \multirow{5}{*}{3706} & \multirow{5}{*}{67865} & \multirow{5}{*}{7.49e-01} & \multirow{5}{*}{4.20e-01} & \textcolor{Red}{CG} & \textcolor{Red}{2.1e-09} & \textcolor{Red}{6.6e-09} & \textcolor{Red}{5.6e-06$^{\ast}$} & \textcolor{Red}{8} & \textcolor{Red}{0.0015} & \textcolor{Red}{11} & \textcolor{Red}{0.4759} \\ 
  &  &  &  &  & \textcolor{Cerulean}{MINRES} & \textcolor{Cerulean}{2.5e-08} & \textcolor{Cerulean}{7.5e-07} & \textcolor{Cerulean}{6.4e-04$^{\ast}$} & \textcolor{Cerulean}{7} & \textcolor{Cerulean}{0.0026} & \textcolor{Cerulean}{8} & \textcolor{Cerulean}{0.3529} \\ 
  &  &  &  &  & \textcolor{OliveGreen}{CG-K} & \textcolor{OliveGreen}{1.2e-07} & \textcolor{OliveGreen}{6.8e-09} & \textcolor{OliveGreen}{4.7e-14$^{\ast}$} & \textcolor{OliveGreen}{15} & \textcolor{OliveGreen}{0.0041} & \textcolor{OliveGreen}{22} & \textcolor{OliveGreen}{0.1054} \\ 
  &  &  &  &  & \textcolor{Violet}{MINRES-K} & \textcolor{Violet}{3.8e-06} & \textcolor{Violet}{7.4e-07} & \textcolor{Violet}{2.5e-06$^{\ast}$} & \textcolor{Violet}{11} & \textcolor{Violet}{0.0050} & \textcolor{Violet}{16} & \textcolor{Violet}{0.0914} \\ 
  &  &  &  &  & \textcolor{Magenta}{LSQR} & \textcolor{Magenta}{2.9e-07} & \textcolor{Magenta}{2.3e-10} & \textcolor{Magenta}{6.5e-06$^{\ast}$} & \textcolor{Magenta}{6} & \textcolor{Magenta}{0.0047} & \textcolor{Magenta}{13} & \textcolor{Magenta}{0.0946} \\ 
\hline 
\multirow{5}{*}{500} & \multirow{5}{*}{1290} & \multirow{5}{*}{21} & \multirow{5}{*}{1.03e-02} & \multirow{5}{*}{1.01e-06} & \textcolor{Red}{CG} & \textcolor{Red}{1.4e-08} & \textcolor{Red}{1.9e-09} & \textcolor{Red}{3.6e-07} & \textcolor{Red}{43} & \textcolor{Red}{0.0072} & \textcolor{Red}{3} & \textcolor{Red}{0.0007} \\ 
  &  &  &  &  & \textcolor{Cerulean}{MINRES} & \textcolor{Cerulean}{1.6e-07} & \textcolor{Cerulean}{1.9e-09} & \textcolor{Cerulean}{4.0e-06} & \textcolor{Cerulean}{38} & \textcolor{Cerulean}{0.0129} & \textcolor{Cerulean}{3} & \textcolor{Cerulean}{0.0013} \\ 
  &  &  &  &  & \textcolor{OliveGreen}{CG-K} & \textcolor{OliveGreen}{6.3e-09} & \textcolor{OliveGreen}{1.9e-09} & \textcolor{OliveGreen}{1.4e-07} & \textcolor{OliveGreen}{93} & \textcolor{OliveGreen}{0.0183} & \textcolor{OliveGreen}{7} & \textcolor{OliveGreen}{0.0013} \\ 
  &  &  &  &  & \textcolor{Violet}{MINRES-K} & \textcolor{Violet}{6.6e-08} & \textcolor{Violet}{1.9e-09} & \textcolor{Violet}{1.7e-06} & \textcolor{Violet}{81} & \textcolor{Violet}{0.0302} & \textcolor{Violet}{7} & \textcolor{Violet}{0.0026} \\ 
  &  &  &  &  & \textcolor{Magenta}{LSQR} & \textcolor{Magenta}{2.2e-08} & \textcolor{Magenta}{1.9e-09} & \textcolor{Magenta}{5.5e-07} & \textcolor{Magenta}{42} & \textcolor{Magenta}{0.0278} & \textcolor{Magenta}{3} & \textcolor{Magenta}{0.0023} \\ 
\hline 
\multirow{5}{*}{500} & \multirow{5}{*}{12394} & \multirow{5}{*}{20315} & \multirow{5}{*}{9.94e-02} & \multirow{5}{*}{9.81e-04} & \textcolor{Red}{CG} & \textcolor{Red}{3.3e-09} & \textcolor{Red}{2.7e-08} & \textcolor{Red}{8.9e-05} & \textcolor{Red}{13} & \textcolor{Red}{0.0030} & \textcolor{Red}{105} & \textcolor{Red}{0.2155} \\ 
  &  &  &  &  & \textcolor{Cerulean}{MINRES} & \textcolor{Cerulean}{7.2e-08} & \textcolor{Cerulean}{3.9e-06} & \textcolor{Cerulean}{1.3e-02} & \textcolor{Cerulean}{11} & \textcolor{Cerulean}{0.0045} & \textcolor{Cerulean}{72} & \textcolor{Cerulean}{0.1755} \\ 
  &  &  &  &  & \textcolor{OliveGreen}{CG-K} & \textcolor{OliveGreen}{3.3e-09} & \textcolor{OliveGreen}{8.0e-09} & \textcolor{OliveGreen}{2.7e-05} & \textcolor{OliveGreen}{27} & \textcolor{OliveGreen}{0.0154} & \textcolor{OliveGreen}{227} & \textcolor{OliveGreen}{0.3955} \\ 
  &  &  &  &  & \textcolor{Violet}{MINRES-K} & \textcolor{Violet}{3.4e-07} & \textcolor{Violet}{1.4e-06} & \textcolor{Violet}{4.7e-03} & \textcolor{Violet}{21} & \textcolor{Violet}{0.0161} & \textcolor{Violet}{161} & \textcolor{Violet}{0.3604} \\ 
  &  &  &  &  & \textcolor{Magenta}{LSQR} & \textcolor{Magenta}{7.0e-08} & \textcolor{Magenta}{5.6e-08} & \textcolor{Magenta}{2.0e-04} & \textcolor{Magenta}{11} & \textcolor{Magenta}{0.0106} & \textcolor{Magenta}{100} & \textcolor{Magenta}{0.2348} \\ 
\hline 
\multirow{5}{*}{500} & \multirow{5}{*}{24788} & \multirow{5}{*}{162986} & \multirow{5}{*}{1.99e-01} & \multirow{5}{*}{7.87e-03} & \textcolor{Red}{CG} & \textcolor{Red}{1.9e-09} & \textcolor{Red}{1.3e-08} & \textcolor{Red}{1.0e-05$^{\ast}$} & \textcolor{Red}{11} & \textcolor{Red}{0.0032} & \textcolor{Red}{25} & \textcolor{Red}{1.447} \\ 
  &  &  &  &  & \textcolor{Cerulean}{MINRES} & \textcolor{Cerulean}{6.8e-08} & \textcolor{Cerulean}{1.8e-06} & \textcolor{Cerulean}{1.4e-03$^{\ast}$} & \textcolor{Cerulean}{9} & \textcolor{Cerulean}{0.0043} & \textcolor{Cerulean}{18} & \textcolor{Cerulean}{1.149} \\ 
  &  &  &  &  & \textcolor{OliveGreen}{CG-K} & \textcolor{OliveGreen}{1.4e-08} & \textcolor{OliveGreen}{7.5e-09} & \textcolor{OliveGreen}{5.7e-14$^{\ast}$} & \textcolor{OliveGreen}{21} & \textcolor{OliveGreen}{0.0303} & \textcolor{OliveGreen}{52} & \textcolor{OliveGreen}{0.7681} \\ 
  &  &  &  &  & \textcolor{Violet}{MINRES-K} & \textcolor{Violet}{2.4e-06} & \textcolor{Violet}{1.6e-06} & \textcolor{Violet}{1.4e-04$^{\ast}$} & \textcolor{Violet}{15} & \textcolor{Violet}{0.0257} & \textcolor{Violet}{37} & \textcolor{Violet}{0.8263} \\ 
  &  &  &  &  & \textcolor{Magenta}{LSQR} & \textcolor{Magenta}{6.7e-08} & \textcolor{Magenta}{3.3e-09} & \textcolor{Magenta}{5.1e-06$^{\ast}$} & \textcolor{Magenta}{9} & \textcolor{Magenta}{0.0123} & \textcolor{Magenta}{27} & \textcolor{Magenta}{0.6047} \\ 
\hline 
\multirow{5}{*}{1000} & \multirow{5}{*}{49690} & \multirow{5}{*}{163767} & \multirow{5}{*}{9.95e-02} & \multirow{5}{*}{9.86e-04} & \textcolor{Red}{CG} & \textcolor{Red}{1.9e-09} & \textcolor{Red}{1.4e-08} & \textcolor{Red}{2.5e-03} & \textcolor{Red}{11} & \textcolor{Red}{0.0049} & \textcolor{Red}{52} & \textcolor{Red}{1.916} \\ 
  &  &  &  &  & \textcolor{Cerulean}{MINRES} & \textcolor{Cerulean}{7.7e-08} & \textcolor{Cerulean}{4.2e-06} & \textcolor{Cerulean}{7.4e-01} & \textcolor{Cerulean}{9} & \textcolor{Cerulean}{0.0057} & \textcolor{Cerulean}{34} & \textcolor{Cerulean}{1.43} \\ 
  &  &  &  &  & \textcolor{OliveGreen}{CG-K} & \textcolor{OliveGreen}{2.2e-09} & \textcolor{OliveGreen}{4.5e-09} & \textcolor{OliveGreen}{8.0e-04} & \textcolor{OliveGreen}{23} & \textcolor{OliveGreen}{0.0731} & \textcolor{OliveGreen}{113} & \textcolor{OliveGreen}{1.976} \\ 
  &  &  &  &  & \textcolor{Violet}{MINRES-K} & \textcolor{Violet}{4.8e-07} & \textcolor{Violet}{1.1e-06} & \textcolor{Violet}{1.9e-01} & \textcolor{Violet}{17} & \textcolor{Violet}{0.0629} & \textcolor{Violet}{77} & \textcolor{Violet}{2.307} \\ 
  &  &  &  &  & \textcolor{Magenta}{LSQR} & \textcolor{Magenta}{7.6e-08} & \textcolor{Magenta}{1.4e-08} & \textcolor{Magenta}{3.1e-03} & \textcolor{Magenta}{9} & \textcolor{Magenta}{0.0252} & \textcolor{Magenta}{52} & \textcolor{Magenta}{1.408} \\ 
\hline 
\end{tabular} 

 \caption{\ER~graphs.} 
 \label{tab:erkrylv}
\end{table}

\begin{table}[h] \scriptsize
 \centering
 \begin{tabular}{|ccccc|c|ccc|cc|cc|}
\hline
\multirow{2}{*}{$N_0$} & \multirow{2}{*}{$N_1$} & \multirow{2}{*}{$N_2$} & Edge & Triangle & Algorithm / & \multicolumn{3}{|c|}{Relative Error} & \multicolumn{2}{|c|}{$\alpha$} & \multicolumn{2}{|c|}{$\beta$} \\
  &  &  & Density & Density & Formulation & $\|\boundary_1^T \alpha\|$ & $\|\boundary_2\beta\|$ & $\|h\|$ & iter. & sec. & iter. & sec. \\ 
 \hline 
\multirow{5}{*}{100} & \multirow{5}{*}{500} & \multirow{5}{*}{729} & \multirow{5}{*}{1.01e-01} & \multirow{5}{*}{4.51e-03} & \textcolor{Red}{CG} & \textcolor{Red}{1.1e-08} & \textcolor{Red}{1.3e-08} & \textcolor{Red}{6.7e-07} & \textcolor{Red}{28} & \textcolor{Red}{0.0043} & \textcolor{Red}{39} & \textcolor{Red}{0.0075} \\ 
  &  &  &  &  & \textcolor{Cerulean}{MINRES} & \textcolor{Cerulean}{3.8e-08} & \textcolor{Cerulean}{3.3e-07} & \textcolor{Cerulean}{1.6e-05} & \textcolor{Cerulean}{27} & \textcolor{Cerulean}{0.0088} & \textcolor{Cerulean}{33} & \textcolor{Cerulean}{0.0122} \\ 
  &  &  &  &  & \textcolor{OliveGreen}{CG-K} & \textcolor{OliveGreen}{1.1e-08} & \textcolor{OliveGreen}{7.4e-09} & \textcolor{OliveGreen}{2.6e-07} & \textcolor{OliveGreen}{57} & \textcolor{OliveGreen}{0.0091} & \textcolor{OliveGreen}{99} & \textcolor{OliveGreen}{0.0182} \\ 
  &  &  &  &  & \textcolor{Violet}{MINRES-K} & \textcolor{Violet}{4.7e-07} & \textcolor{Violet}{1.5e-07} & \textcolor{Violet}{9.5e-06} & \textcolor{Violet}{49} & \textcolor{Violet}{0.0164} & \textcolor{Violet}{71} & \textcolor{Violet}{0.0254} \\ 
  &  &  &  &  & \textcolor{Magenta}{LSQR} & \textcolor{Magenta}{2.0e-07} & \textcolor{Magenta}{2.3e-08} & \textcolor{Magenta}{2.9e-06} & \textcolor{Magenta}{25} & \textcolor{Magenta}{0.0160} & \textcolor{Magenta}{38} & \textcolor{Magenta}{0.0251} \\ 
\hline 
\multirow{5}{*}{100} & \multirow{5}{*}{1000} & \multirow{5}{*}{3655} & \multirow{5}{*}{2.02e-01} & \multirow{5}{*}{2.26e-02} & \textcolor{Red}{CG} & \textcolor{Red}{2.6e-09} & \textcolor{Red}{1.6e-08} & \textcolor{Red}{2.3e-06} & \textcolor{Red}{18} & \textcolor{Red}{0.0029} & \textcolor{Red}{44} & \textcolor{Red}{0.0252} \\ 
  &  &  &  &  & \textcolor{Cerulean}{MINRES} & \textcolor{Cerulean}{1.5e-08} & \textcolor{Cerulean}{6.0e-07} & \textcolor{Cerulean}{9.0e-05} & \textcolor{Cerulean}{17} & \textcolor{Cerulean}{0.0057} & \textcolor{Cerulean}{35} & \textcolor{Cerulean}{0.0264} \\ 
  &  &  &  &  & \textcolor{OliveGreen}{CG-K} & \textcolor{OliveGreen}{1.5e-08} & \textcolor{OliveGreen}{4.3e-09} & \textcolor{OliveGreen}{7.2e-07} & \textcolor{OliveGreen}{35} & \textcolor{OliveGreen}{0.0062} & \textcolor{OliveGreen}{97} & \textcolor{OliveGreen}{0.0295} \\ 
  &  &  &  &  & \textcolor{Violet}{MINRES-K} & \textcolor{Violet}{1.7e-06} & \textcolor{Violet}{1.5e-07} & \textcolor{Violet}{4.2e-05} & \textcolor{Violet}{29} & \textcolor{Violet}{0.0103} & \textcolor{Violet}{79} & \textcolor{Violet}{0.0381} \\ 
  &  &  &  &  & \textcolor{Magenta}{LSQR} & \textcolor{Magenta}{8.2e-08} & \textcolor{Magenta}{1.6e-08} & \textcolor{Magenta}{2.9e-06} & \textcolor{Magenta}{16} & \textcolor{Magenta}{0.0106} & \textcolor{Magenta}{44} & \textcolor{Magenta}{0.0346} \\ 
\hline 
\multirow{5}{*}{100} & \multirow{5}{*}{1500} & \multirow{5}{*}{8354} & \multirow{5}{*}{3.03e-01} & \multirow{5}{*}{5.17e-02} & \textcolor{Red}{CG} & \textcolor{Red}{7.4e-09} & \textcolor{Red}{1.7e-08} & \textcolor{Red}{6.7e-05} & \textcolor{Red}{13} & \textcolor{Red}{0.0021} & \textcolor{Red}{109} & \textcolor{Red}{0.2024} \\ 
  &  &  &  &  & \textcolor{Cerulean}{MINRES} & \textcolor{Cerulean}{7.4e-09} & \textcolor{Cerulean}{2.4e-06} & \textcolor{Cerulean}{9.6e-03} & \textcolor{Cerulean}{13} & \textcolor{Cerulean}{0.0044} & \textcolor{Cerulean}{84} & \textcolor{Cerulean}{0.1743} \\ 
  &  &  &  &  & \textcolor{OliveGreen}{CG-K} & \textcolor{OliveGreen}{8.8e-09} & \textcolor{OliveGreen}{7.2e-09} & \textcolor{OliveGreen}{2.9e-05} & \textcolor{OliveGreen}{27} & \textcolor{OliveGreen}{0.0053} & \textcolor{OliveGreen}{227} & \textcolor{OliveGreen}{0.1091} \\ 
  &  &  &  &  & \textcolor{Violet}{MINRES-K} & \textcolor{Violet}{2.9e-06} & \textcolor{Violet}{5.2e-07} & \textcolor{Violet}{2.3e-03} & \textcolor{Violet}{21} & \textcolor{Violet}{0.0079} & \textcolor{Violet}{185} & \textcolor{Violet}{0.1249} \\ 
  &  &  &  &  & \textcolor{Magenta}{LSQR} & \textcolor{Magenta}{6.3e-08} & \textcolor{Magenta}{1.7e-08} & \textcolor{Magenta}{7.2e-05} & \textcolor{Magenta}{12} & \textcolor{Magenta}{0.0083} & \textcolor{Magenta}{109} & \textcolor{Magenta}{0.1067} \\ 
\hline 
\multirow{5}{*}{100} & \multirow{5}{*}{2000} & \multirow{5}{*}{15530} & \multirow{5}{*}{4.04e-01} & \multirow{5}{*}{9.60e-02} & \textcolor{Red}{CG} & \textcolor{Red}{7.1e-09} & \textcolor{Red}{2.0e-08} & \textcolor{Red}{4.8e-06$^{\ast}$} & \textcolor{Red}{11} & \textcolor{Red}{0.0019} & \textcolor{Red}{68} & \textcolor{Red}{0.2842} \\ 
  &  &  &  &  & \textcolor{Cerulean}{MINRES} & \textcolor{Cerulean}{7.1e-09} & \textcolor{Cerulean}{2.9e-06} & \textcolor{Cerulean}{6.9e-04$^{\ast}$} & \textcolor{Cerulean}{11} & \textcolor{Cerulean}{0.0038} & \textcolor{Cerulean}{53} & \textcolor{Cerulean}{0.2332} \\ 
  &  &  &  &  & \textcolor{OliveGreen}{CG-K} & \textcolor{OliveGreen}{7.1e-09} & \textcolor{OliveGreen}{8.4e-09} & \textcolor{OliveGreen}{1.7e-14$^{\ast}$} & \textcolor{OliveGreen}{23} & \textcolor{OliveGreen}{0.0049} & \textcolor{OliveGreen}{142} & \textcolor{OliveGreen}{0.1167} \\ 
  &  &  &  &  & \textcolor{Violet}{MINRES-K} & \textcolor{Violet}{5.0e-06} & \textcolor{Violet}{4.7e-07} & \textcolor{Violet}{6.8e-05$^{\ast}$} & \textcolor{Violet}{17} & \textcolor{Violet}{0.0067} & \textcolor{Violet}{118} & \textcolor{Violet}{0.1401} \\ 
  &  &  &  &  & \textcolor{Magenta}{LSQR} & \textcolor{Magenta}{5.0e-08} & \textcolor{Magenta}{6.1e-09} & \textcolor{Magenta}{1.7e-06$^{\ast}$} & \textcolor{Magenta}{10} & \textcolor{Magenta}{0.0071} & \textcolor{Magenta}{71} & \textcolor{Magenta}{0.0977} \\ 
\hline 
\multirow{5}{*}{500} & \multirow{5}{*}{2500} & \multirow{5}{*}{3720} & \multirow{5}{*}{2.00e-02} & \multirow{5}{*}{1.80e-04} & \textcolor{Red}{CG} & \textcolor{Red}{1.8e-08} & \textcolor{Red}{1.1e-08} & \textcolor{Red}{1.2e-06} & \textcolor{Red}{38} & \textcolor{Red}{0.0066} & \textcolor{Red}{44} & \textcolor{Red}{0.0169} \\ 
  &  &  &  &  & \textcolor{Cerulean}{MINRES} & \textcolor{Cerulean}{1.8e-07} & \textcolor{Cerulean}{6.3e-07} & \textcolor{Cerulean}{6.2e-05} & \textcolor{Cerulean}{33} & \textcolor{Cerulean}{0.0114} & \textcolor{Cerulean}{36} & \textcolor{Cerulean}{0.0203} \\ 
  &  &  &  &  & \textcolor{OliveGreen}{CG-K} & \textcolor{OliveGreen}{1.8e-08} & \textcolor{OliveGreen}{4.8e-09} & \textcolor{OliveGreen}{6.9e-07} & \textcolor{OliveGreen}{77} & \textcolor{OliveGreen}{0.0184} & \textcolor{OliveGreen}{93} & \textcolor{OliveGreen}{0.0326} \\ 
  &  &  &  &  & \textcolor{Violet}{MINRES-K} & \textcolor{Violet}{1.7e-06} & \textcolor{Violet}{2.2e-07} & \textcolor{Violet}{5.2e-05} & \textcolor{Violet}{57} & \textcolor{Violet}{0.0237} & \textcolor{Violet}{77} & \textcolor{Violet}{0.0410} \\ 
  &  &  &  &  & \textcolor{Magenta}{LSQR} & \textcolor{Magenta}{2.3e-07} & \textcolor{Magenta}{2.8e-08} & \textcolor{Magenta}{7.0e-06} & \textcolor{Magenta}{32} & \textcolor{Magenta}{0.0222} & \textcolor{Magenta}{42} & \textcolor{Magenta}{0.0346} \\ 
\hline 
\multirow{5}{*}{500} & \multirow{5}{*}{5000} & \multirow{5}{*}{16948} & \multirow{5}{*}{4.01e-02} & \multirow{5}{*}{8.18e-04} & \textcolor{Red}{CG} & \textcolor{Red}{1.2e-08} & \textcolor{Red}{1.6e-08} & \textcolor{Red}{4.5e-06} & \textcolor{Red}{32} & \textcolor{Red}{0.0059} & \textcolor{Red}{34} & \textcolor{Red}{0.0898} \\ 
  &  &  &  &  & \textcolor{Cerulean}{MINRES} & \textcolor{Cerulean}{1.0e-07} & \textcolor{Cerulean}{6.2e-07} & \textcolor{Cerulean}{1.7e-04} & \textcolor{Cerulean}{28} & \textcolor{Cerulean}{0.0100} & \textcolor{Cerulean}{26} & \textcolor{Cerulean}{0.0771} \\ 
  &  &  &  &  & \textcolor{OliveGreen}{CG-K} & \textcolor{OliveGreen}{2.2e-08} & \textcolor{OliveGreen}{4.0e-09} & \textcolor{OliveGreen}{1.4e-06} & \textcolor{OliveGreen}{63} & \textcolor{OliveGreen}{0.0201} & \textcolor{OliveGreen}{75} & \textcolor{OliveGreen}{0.0795} \\ 
  &  &  &  &  & \textcolor{Violet}{MINRES-K} & \textcolor{Violet}{5.2e-06} & \textcolor{Violet}{2.0e-07} & \textcolor{Violet}{2.2e-04} & \textcolor{Violet}{43} & \textcolor{Violet}{0.0215} & \textcolor{Violet}{57} & \textcolor{Violet}{0.0841} \\ 
  &  &  &  &  & \textcolor{Magenta}{LSQR} & \textcolor{Magenta}{4.7e-07} & \textcolor{Magenta}{1.1e-08} & \textcolor{Magenta}{1.9e-05} & \textcolor{Magenta}{25} & \textcolor{Magenta}{0.0187} & \textcolor{Magenta}{35} & \textcolor{Magenta}{0.0594} \\ 
\hline 
\multirow{5}{*}{500} & \multirow{5}{*}{12500} & \multirow{5}{*}{110507} & \multirow{5}{*}{1.00e-01} & \multirow{5}{*}{5.34e-03} & \textcolor{Red}{CG} & \textcolor{Red}{8.4e-09} & \textcolor{Red}{4.0e-08} & \textcolor{Red}{4.5e-05} & \textcolor{Red}{23} & \textcolor{Red}{0.0051} & \textcolor{Red}{161} & \textcolor{Red}{6.469} \\ 
  &  &  &  &  & \textcolor{Cerulean}{MINRES} & \textcolor{Cerulean}{2.0e-08} & \textcolor{Cerulean}{1.8e-05} & \textcolor{Cerulean}{2.1e-02} & \textcolor{Cerulean}{22} & \textcolor{Cerulean}{0.0088} & \textcolor{Cerulean}{86} & \textcolor{Cerulean}{3.632} \\ 
  &  &  &  &  & \textcolor{OliveGreen}{CG-K} & \textcolor{OliveGreen}{1.9e-08} & \textcolor{OliveGreen}{7.6e-09} & \textcolor{OliveGreen}{8.5e-06} & \textcolor{OliveGreen}{45} & \textcolor{OliveGreen}{0.0254} & \textcolor{OliveGreen}{367} & \textcolor{OliveGreen}{3.337} \\ 
  &  &  &  &  & \textcolor{Violet}{MINRES-K} & \textcolor{Violet}{1.3e-05} & \textcolor{Violet}{3.6e-06} & \textcolor{Violet}{4.1e-03} & \textcolor{Violet}{29} & \textcolor{Violet}{0.0223} & \textcolor{Violet}{217} & \textcolor{Violet}{2.617} \\ 
  &  &  &  &  & \textcolor{Magenta}{LSQR} & \textcolor{Magenta}{9.5e-07} & \textcolor{Magenta}{2.6e-08} & \textcolor{Magenta}{7.4e-05} & \textcolor{Magenta}{18} & \textcolor{Magenta}{0.0176} & \textcolor{Magenta}{166} & \textcolor{Magenta}{1.996} \\ 
\hline 
\multirow{5}{*}{1000} & \multirow{5}{*}{5000} & \multirow{5}{*}{7386} & \multirow{5}{*}{1.00e-02} & \multirow{5}{*}{4.44e-05} & \textcolor{Red}{CG} & \textcolor{Red}{1.4e-08} & \textcolor{Red}{1.2e-08} & \textcolor{Red}{1.8e-06} & \textcolor{Red}{43} & \textcolor{Red}{0.0085} & \textcolor{Red}{44} & \textcolor{Red}{0.0278} \\ 
  &  &  &  &  & \textcolor{Cerulean}{MINRES} & \textcolor{Cerulean}{3.0e-07} & \textcolor{Cerulean}{7.5e-07} & \textcolor{Cerulean}{1.0e-04} & \textcolor{Cerulean}{36} & \textcolor{Cerulean}{0.0135} & \textcolor{Cerulean}{34} & \textcolor{Cerulean}{0.0280} \\ 
  &  &  &  &  & \textcolor{OliveGreen}{CG-K} & \textcolor{OliveGreen}{2.2e-08} & \textcolor{OliveGreen}{5.6e-09} & \textcolor{OliveGreen}{1.1e-06} & \textcolor{OliveGreen}{85} & \textcolor{OliveGreen}{0.0284} & \textcolor{OliveGreen}{93} & \textcolor{OliveGreen}{0.0522} \\ 
  &  &  &  &  & \textcolor{Violet}{MINRES-K} & \textcolor{Violet}{2.9e-06} & \textcolor{Violet}{3.4e-07} & \textcolor{Violet}{1.2e-04} & \textcolor{Violet}{61} & \textcolor{Violet}{0.0315} & \textcolor{Violet}{73} & \textcolor{Violet}{0.0554} \\ 
  &  &  &  &  & \textcolor{Magenta}{LSQR} & \textcolor{Magenta}{3.6e-07} & \textcolor{Magenta}{2.8e-08} & \textcolor{Magenta}{1.4e-05} & \textcolor{Magenta}{35} & \textcolor{Magenta}{0.0268} & \textcolor{Magenta}{42} & \textcolor{Magenta}{0.0437} \\ 
\hline 
\multirow{5}{*}{1000} & \multirow{5}{*}{10000} & \multirow{5}{*}{33022} & \multirow{5}{*}{2.00e-02} & \multirow{5}{*}{1.99e-04} & \textcolor{Red}{CG} & \textcolor{Red}{1.7e-08} & \textcolor{Red}{1.4e-08} & \textcolor{Red}{5.2e-06} & \textcolor{Red}{32} & \textcolor{Red}{0.0072} & \textcolor{Red}{39} & \textcolor{Red}{0.1953} \\ 
  &  &  &  &  & \textcolor{Cerulean}{MINRES} & \textcolor{Cerulean}{2.1e-07} & \textcolor{Cerulean}{1.3e-06} & \textcolor{Cerulean}{4.8e-04} & \textcolor{Cerulean}{27} & \textcolor{Cerulean}{0.0108} & \textcolor{Cerulean}{28} & \textcolor{Cerulean}{0.1508} \\ 
  &  &  &  &  & \textcolor{OliveGreen}{CG-K} & \textcolor{OliveGreen}{2.3e-08} & \textcolor{OliveGreen}{5.0e-09} & \textcolor{OliveGreen}{2.1e-06} & \textcolor{OliveGreen}{65} & \textcolor{OliveGreen}{0.0324} & \textcolor{OliveGreen}{85} & \textcolor{OliveGreen}{0.2098} \\ 
  &  &  &  &  & \textcolor{Violet}{MINRES-K} & \textcolor{Violet}{6.7e-06} & \textcolor{Violet}{5.3e-07} & \textcolor{Violet}{4.2e-04} & \textcolor{Violet}{43} & \textcolor{Violet}{0.0297} & \textcolor{Violet}{63} & \textcolor{Violet}{0.1852} \\ 
  &  &  &  &  & \textcolor{Magenta}{LSQR} & \textcolor{Magenta}{3.0e-07} & \textcolor{Magenta}{1.4e-08} & \textcolor{Magenta}{1.7e-05} & \textcolor{Magenta}{26} & \textcolor{Magenta}{0.0230} & \textcolor{Magenta}{39} & \textcolor{Magenta}{0.1365} \\ 
\hline 
\multirow{5}{*}{1000} & \multirow{5}{*}{25000} & \multirow{5}{*}{220002} & \multirow{5}{*}{5.01e-02} & \multirow{5}{*}{1.32e-03} & \textcolor{Red}{CG} & \textcolor{Red}{1.0e-08} & \textcolor{Red}{2.4e-08} & \textcolor{Red}{3.4e-05} & \textcolor{Red}{27} & \textcolor{Red}{0.0080} & \textcolor{Red}{85} & \textcolor{Red}{7.04} \\ 
  &  &  &  &  & \textcolor{Cerulean}{MINRES} & \textcolor{Cerulean}{2.1e-07} & \textcolor{Cerulean}{1.2e-05} & \textcolor{Cerulean}{1.8e-02} & \textcolor{Cerulean}{23} & \textcolor{Cerulean}{0.0110} & \textcolor{Cerulean}{48} & \textcolor{Cerulean}{4.52} \\ 
  &  &  &  &  & \textcolor{OliveGreen}{CG-K} & \textcolor{OliveGreen}{2.4e-08} & \textcolor{OliveGreen}{4.6e-09} & \textcolor{OliveGreen}{6.9e-06} & \textcolor{OliveGreen}{53} & \textcolor{OliveGreen}{0.0688} & \textcolor{OliveGreen}{189} & \textcolor{OliveGreen}{3.7} \\ 
  &  &  &  &  & \textcolor{Violet}{MINRES-K} & \textcolor{Violet}{2.6e-05} & \textcolor{Violet}{2.3e-06} & \textcolor{Violet}{4.0e-03} & \textcolor{Violet}{31} & \textcolor{Violet}{0.0554} & \textcolor{Violet}{119} & \textcolor{Violet}{3.588} \\ 
  &  &  &  &  & \textcolor{Magenta}{LSQR} & \textcolor{Magenta}{5.9e-07} & \textcolor{Magenta}{1.2e-08} & \textcolor{Magenta}{5.5e-05} & \textcolor{Magenta}{21} & \textcolor{Magenta}{0.0383} & \textcolor{Magenta}{89} & \textcolor{Magenta}{2.712} \\ 
\hline 
\end{tabular} 

 \caption{\WS~graphs.}
 \label{tab:wskrylv} 
\end{table}

\begin{table}[h] \scriptsize
  \centering
  \begin{tabular}{|ccccc|c|ccc|cc|cc|}
\hline
\multirow{2}{*}{$N_0$} & \multirow{2}{*}{$N_1$} & \multirow{2}{*}{$N_2$} & Edge & Triangle & Algorithm / & \multicolumn{3}{|c|}{Relative Error} & \multicolumn{2}{|c|}{$\alpha$} & \multicolumn{2}{|c|}{$\beta$} \\
  &  &  & Density & Density & Formulation & $\|\boundary_1^T \alpha\|$ & $\|\boundary_2\beta\|$ & $\|h\|$ & iter. & sec. & iter. & sec. \\ 
 \hline 
\multirow{5}{*}{100} & \multirow{5}{*}{475} & \multirow{5}{*}{301} & \multirow{5}{*}{9.60e-02} & \multirow{5}{*}{1.86e-03} & \textcolor{Red}{CG} & \textcolor{Red}{1.1e-08} & \textcolor{Red}{2.4e-08} & \textcolor{Red}{1.3e-06} & \textcolor{Red}{25} & \textcolor{Red}{0.0038} & \textcolor{Red}{79} & \textcolor{Red}{0.0127} \\ 
  &  &  &  &  & \textcolor{Cerulean}{MINRES} & \textcolor{Cerulean}{2.7e-08} & \textcolor{Cerulean}{3.2e-07} & \textcolor{Cerulean}{1.8e-05} & \textcolor{Cerulean}{24} & \textcolor{Cerulean}{0.0078} & \textcolor{Cerulean}{70} & \textcolor{Cerulean}{0.0233} \\ 
  &  &  &  &  & \textcolor{OliveGreen}{CG-K} & \textcolor{OliveGreen}{8.3e-09} & \textcolor{OliveGreen}{7.3e-09} & \textcolor{OliveGreen}{4.1e-07} & \textcolor{OliveGreen}{55} & \textcolor{OliveGreen}{0.0088} & \textcolor{OliveGreen}{171} & \textcolor{OliveGreen}{0.0282} \\ 
  &  &  &  &  & \textcolor{Violet}{MINRES-K} & \textcolor{Violet}{1.6e-07} & \textcolor{Violet}{2.0e-07} & \textcolor{Violet}{1.2e-05} & \textcolor{Violet}{45} & \textcolor{Violet}{0.0149} & \textcolor{Violet}{145} & \textcolor{Violet}{0.0497} \\ 
  &  &  &  &  & \textcolor{Magenta}{LSQR} & \textcolor{Magenta}{1.5e-07} & \textcolor{Magenta}{7.0e-08} & \textcolor{Magenta}{4.9e-06} & \textcolor{Magenta}{22} & \textcolor{Magenta}{0.0142} & \textcolor{Magenta}{75} & \textcolor{Magenta}{0.0476} \\ 
\hline 
\multirow{5}{*}{100} & \multirow{5}{*}{900} & \multirow{5}{*}{1701} & \multirow{5}{*}{1.82e-01} & \multirow{5}{*}{1.05e-02} & \textcolor{Red}{CG} & \textcolor{Red}{1.1e-08} & \textcolor{Red}{2.6e-08} & \textcolor{Red}{1.0e-05} & \textcolor{Red}{21} & \textcolor{Red}{0.0033} & \textcolor{Red}{96} & \textcolor{Red}{0.0281} \\ 
  &  &  &  &  & \textcolor{Cerulean}{MINRES} & \textcolor{Cerulean}{3.3e-08} & \textcolor{Cerulean}{1.1e-06} & \textcolor{Cerulean}{4.2e-04} & \textcolor{Cerulean}{20} & \textcolor{Cerulean}{0.0067} & \textcolor{Cerulean}{74} & \textcolor{Cerulean}{0.0351} \\ 
  &  &  &  &  & \textcolor{OliveGreen}{CG-K} & \textcolor{OliveGreen}{1.5e-08} & \textcolor{OliveGreen}{7.5e-09} & \textcolor{OliveGreen}{2.6e-06} & \textcolor{OliveGreen}{45} & \textcolor{OliveGreen}{0.0079} & \textcolor{OliveGreen}{211} & \textcolor{OliveGreen}{0.0494} \\ 
  &  &  &  &  & \textcolor{Violet}{MINRES-K} & \textcolor{Violet}{6.0e-07} & \textcolor{Violet}{3.8e-07} & \textcolor{Violet}{1.6e-04} & \textcolor{Violet}{35} & \textcolor{Violet}{0.0123} & \textcolor{Violet}{165} & \textcolor{Violet}{0.0682} \\ 
  &  &  &  &  & \textcolor{Magenta}{LSQR} & \textcolor{Magenta}{2.4e-07} & \textcolor{Magenta}{4.6e-08} & \textcolor{Magenta}{2.5e-05} & \textcolor{Magenta}{18} & \textcolor{Magenta}{0.0119} & \textcolor{Magenta}{93} & \textcolor{Magenta}{0.0653} \\ 
\hline 
\multirow{5}{*}{100} & \multirow{5}{*}{1600} & \multirow{5}{*}{8105} & \multirow{5}{*}{3.23e-01} & \multirow{5}{*}{5.01e-02} & \textcolor{Red}{CG} & \textcolor{Red}{6.7e-09} & \textcolor{Red}{1.3e-08} & \textcolor{Red}{2.2e-06$^{\ast}$} & \textcolor{Red}{20} & \textcolor{Red}{0.0032} & \textcolor{Red}{50} & \textcolor{Red}{0.0967} \\ 
  &  &  &  &  & \textcolor{Cerulean}{MINRES} & \textcolor{Cerulean}{2.0e-08} & \textcolor{Cerulean}{7.8e-07} & \textcolor{Cerulean}{1.3e-04$^{\ast}$} & \textcolor{Cerulean}{19} & \textcolor{Cerulean}{0.0063} & \textcolor{Cerulean}{39} & \textcolor{Cerulean}{0.0832} \\ 
  &  &  &  &  & \textcolor{OliveGreen}{CG-K} & \textcolor{OliveGreen}{2.3e-08} & \textcolor{OliveGreen}{8.1e-09} & \textcolor{OliveGreen}{9.3e-11$^{\ast}$} & \textcolor{OliveGreen}{39} & \textcolor{OliveGreen}{0.0077} & \textcolor{OliveGreen}{106} & \textcolor{OliveGreen}{0.0506} \\ 
  &  &  &  &  & \textcolor{Violet}{MINRES-K} & \textcolor{Violet}{1.4e-06} & \textcolor{Violet}{3.4e-07} & \textcolor{Violet}{1.2e-05$^{\ast}$} & \textcolor{Violet}{31} & \textcolor{Violet}{0.0116} & \textcolor{Violet}{82} & \textcolor{Violet}{0.0545} \\ 
  &  &  &  &  & \textcolor{Magenta}{LSQR} & \textcolor{Magenta}{4.7e-07} & \textcolor{Magenta}{9.3e-09} & \textcolor{Magenta}{7.0e-06$^{\ast}$} & \textcolor{Magenta}{16} & \textcolor{Magenta}{0.0108} & \textcolor{Magenta}{51} & \textcolor{Magenta}{0.0499} \\ 
\hline 
\multirow{5}{*}{100} & \multirow{5}{*}{2400} & \multirow{5}{*}{24497} & \multirow{5}{*}{4.85e-01} & \multirow{5}{*}{1.51e-01} & \textcolor{Red}{CG} & \textcolor{Red}{5.2e-09} & \textcolor{Red}{1.2e-08} & \textcolor{Red}{4.3e-06$^{\ast}$} & \textcolor{Red}{18} & \textcolor{Red}{0.0030} & \textcolor{Red}{27} & \textcolor{Red}{0.2435} \\ 
  &  &  &  &  & \textcolor{Cerulean}{MINRES} & \textcolor{Cerulean}{5.0e-08} & \textcolor{Cerulean}{6.8e-07} & \textcolor{Cerulean}{2.4e-04$^{\ast}$} & \textcolor{Cerulean}{16} & \textcolor{Cerulean}{0.0054} & \textcolor{Cerulean}{21} & \textcolor{Cerulean}{0.1960} \\ 
  &  &  &  &  & \textcolor{OliveGreen}{CG-K} & \textcolor{OliveGreen}{1.8e-08} & \textcolor{OliveGreen}{7.4e-09} & \textcolor{OliveGreen}{2.2e-14$^{\ast}$} & \textcolor{OliveGreen}{35} & \textcolor{OliveGreen}{0.0078} & \textcolor{OliveGreen}{56} & \textcolor{OliveGreen}{0.0817} \\ 
  &  &  &  &  & \textcolor{Violet}{MINRES-K} & \textcolor{Violet}{4.6e-06} & \textcolor{Violet}{3.1e-07} & \textcolor{Violet}{6.1e-06$^{\ast}$} & \textcolor{Violet}{25} & \textcolor{Violet}{0.0100} & \textcolor{Violet}{44} & \textcolor{Violet}{0.0831} \\ 
  &  &  &  &  & \textcolor{Magenta}{LSQR} & \textcolor{Magenta}{4.2e-07} & \textcolor{Magenta}{3.2e-09} & \textcolor{Magenta}{8.0e-06$^{\ast}$} & \textcolor{Magenta}{14} & \textcolor{Magenta}{0.0097} & \textcolor{Magenta}{29} & \textcolor{Magenta}{0.0632} \\ 
\hline 
\multirow{5}{*}{500} & \multirow{5}{*}{4900} & \multirow{5}{*}{4740} & \multirow{5}{*}{3.93e-02} & \multirow{5}{*}{2.29e-04} & \textcolor{Red}{CG} & \textcolor{Red}{1.0e-08} & \textcolor{Red}{3.4e-08} & \textcolor{Red}{9.1e-06} & \textcolor{Red}{32} & \textcolor{Red}{0.0059} & \textcolor{Red}{165} & \textcolor{Red}{0.0956} \\ 
  &  &  &  &  & \textcolor{Cerulean}{MINRES} & \textcolor{Cerulean}{4.7e-08} & \textcolor{Cerulean}{1.8e-06} & \textcolor{Cerulean}{4.8e-04} & \textcolor{Cerulean}{30} & \textcolor{Cerulean}{0.0107} & \textcolor{Cerulean}{118} & \textcolor{Cerulean}{0.0904} \\ 
  &  &  &  &  & \textcolor{OliveGreen}{CG-K} & \textcolor{OliveGreen}{2.2e-08} & \textcolor{OliveGreen}{8.2e-09} & \textcolor{OliveGreen}{2.2e-06} & \textcolor{OliveGreen}{73} & \textcolor{OliveGreen}{0.0231} & \textcolor{OliveGreen}{367} & \textcolor{OliveGreen}{0.1707} \\ 
  &  &  &  &  & \textcolor{Violet}{MINRES-K} & \textcolor{Violet}{6.8e-07} & \textcolor{Violet}{8.2e-07} & \textcolor{Violet}{2.2e-04} & \textcolor{Violet}{53} & \textcolor{Violet}{0.0262} & \textcolor{Violet}{257} & \textcolor{Violet}{0.1669} \\ 
  &  &  &  &  & \textcolor{Magenta}{LSQR} & \textcolor{Magenta}{3.1e-07} & \textcolor{Magenta}{1.2e-07} & \textcolor{Magenta}{3.8e-05} & \textcolor{Magenta}{27} & \textcolor{Magenta}{0.0202} & \textcolor{Magenta}{150} & \textcolor{Magenta}{0.1378} \\ 
\hline 
\multirow{5}{*}{500} & \multirow{5}{*}{9600} & \multirow{5}{*}{25016} & \multirow{5}{*}{7.70e-02} & \multirow{5}{*}{1.21e-03} & \textcolor{Red}{CG} & \textcolor{Red}{8.4e-09} & \textcolor{Red}{3.0e-08} & \textcolor{Red}{5.8e-05} & \textcolor{Red}{27} & \textcolor{Red}{0.0056} & \textcolor{Red}{215} & \textcolor{Red}{1.043} \\ 
  &  &  &  &  & \textcolor{Cerulean}{MINRES} & \textcolor{Cerulean}{8.2e-08} & \textcolor{Cerulean}{4.7e-06} & \textcolor{Cerulean}{9.1e-03} & \textcolor{Cerulean}{24} & \textcolor{Cerulean}{0.0092} & \textcolor{Cerulean}{133} & \textcolor{Cerulean}{0.6919} \\ 
  &  &  &  &  & \textcolor{OliveGreen}{CG-K} & \textcolor{OliveGreen}{1.1e-08} & \textcolor{OliveGreen}{8.5e-09} & \textcolor{OliveGreen}{1.3e-05} & \textcolor{OliveGreen}{57} & \textcolor{OliveGreen}{0.0267} & \textcolor{OliveGreen}{477} & \textcolor{OliveGreen}{0.9132} \\ 
  &  &  &  &  & \textcolor{Violet}{MINRES-K} & \textcolor{Violet}{1.5e-06} & \textcolor{Violet}{1.7e-06} & \textcolor{Violet}{3.3e-03} & \textcolor{Violet}{41} & \textcolor{Violet}{0.0270} & \textcolor{Violet}{303} & \textcolor{Violet}{0.7096} \\ 
  &  &  &  &  & \textcolor{Magenta}{LSQR} & \textcolor{Magenta}{3.5e-07} & \textcolor{Magenta}{7.9e-08} & \textcolor{Magenta}{1.8e-04} & \textcolor{Magenta}{22} & \textcolor{Magenta}{0.0191} & \textcolor{Magenta}{202} & \textcolor{Magenta}{0.5085} \\ 
\hline 
\multirow{5}{*}{500} & \multirow{5}{*}{18400} & \multirow{5}{*}{133933} & \multirow{5}{*}{1.47e-01} & \multirow{5}{*}{6.47e-03} & \textcolor{Red}{CG} & \textcolor{Red}{4.5e-09} & \textcolor{Red}{2.2e-08} & \textcolor{Red}{2.0e-03} & \textcolor{Red}{23} & \textcolor{Red}{0.0058} & \textcolor{Red}{111} & \textcolor{Red}{7.204} \\ 
  &  &  &  &  & \textcolor{Cerulean}{MINRES} & \textcolor{Cerulean}{6.8e-08} & \textcolor{Cerulean}{5.8e-06} & \textcolor{Cerulean}{5.3e-01} & \textcolor{Cerulean}{20} & \textcolor{Cerulean}{0.0086} & \textcolor{Cerulean}{69} & \textcolor{Cerulean}{4.667} \\ 
  &  &  &  &  & \textcolor{OliveGreen}{CG-K} & \textcolor{OliveGreen}{1.8e-08} & \textcolor{OliveGreen}{3.8e-09} & \textcolor{OliveGreen}{3.6e-04} & \textcolor{OliveGreen}{45} & \textcolor{OliveGreen}{0.0414} & \textcolor{OliveGreen}{251} & \textcolor{OliveGreen}{2.809} \\ 
  &  &  &  &  & \textcolor{Violet}{MINRES-K} & \textcolor{Violet}{1.1e-05} & \textcolor{Violet}{1.2e-06} & \textcolor{Violet}{1.3e-01} & \textcolor{Violet}{29} & \textcolor{Violet}{0.0362} & \textcolor{Violet}{163} & \textcolor{Violet}{2.646} \\ 
  &  &  &  &  & \textcolor{Magenta}{LSQR} & \textcolor{Magenta}{3.4e-07} & \textcolor{Magenta}{1.7e-08} & \textcolor{Magenta}{2.5e-03} & \textcolor{Magenta}{18} & \textcolor{Magenta}{0.0205} & \textcolor{Magenta}{113} & \textcolor{Magenta}{1.867} \\ 
\hline 
\multirow{5}{*}{1000} & \multirow{5}{*}{9900} & \multirow{5}{*}{6264} & \multirow{5}{*}{1.98e-02} & \multirow{5}{*}{3.77e-05} & \textcolor{Red}{CG} & \textcolor{Red}{7.4e-09} & \textcolor{Red}{4.4e-08} & \textcolor{Red}{1.1e-05} & \textcolor{Red}{41} & \textcolor{Red}{0.0091} & \textcolor{Red}{211} & \textcolor{Red}{0.1410} \\ 
  &  &  &  &  & \textcolor{Cerulean}{MINRES} & \textcolor{Cerulean}{8.5e-08} & \textcolor{Cerulean}{3.2e-06} & \textcolor{Cerulean}{8.1e-04} & \textcolor{Cerulean}{37} & \textcolor{Cerulean}{0.0147} & \textcolor{Cerulean}{146} & \textcolor{Cerulean}{0.1246} \\ 
  &  &  &  &  & \textcolor{OliveGreen}{CG-K} & \textcolor{OliveGreen}{1.4e-08} & \textcolor{OliveGreen}{1.1e-08} & \textcolor{OliveGreen}{2.8e-06} & \textcolor{OliveGreen}{89} & \textcolor{OliveGreen}{0.0442} & \textcolor{OliveGreen}{481} & \textcolor{OliveGreen}{0.3244} \\ 
  &  &  &  &  & \textcolor{Violet}{MINRES-K} & \textcolor{Violet}{4.5e-07} & \textcolor{Violet}{1.1e-06} & \textcolor{Violet}{2.9e-04} & \textcolor{Violet}{69} & \textcolor{Violet}{0.0474} & \textcolor{Violet}{321} & \textcolor{Violet}{0.2883} \\ 
  &  &  &  &  & \textcolor{Magenta}{LSQR} & \textcolor{Magenta}{2.3e-07} & \textcolor{Magenta}{1.8e-07} & \textcolor{Magenta}{5.0e-05} & \textcolor{Magenta}{35} & \textcolor{Magenta}{0.0309} & \textcolor{Magenta}{185} & \textcolor{Magenta}{0.2043} \\ 
\hline 
\multirow{5}{*}{1000} & \multirow{5}{*}{19600} & \multirow{5}{*}{37365} & \multirow{5}{*}{3.92e-02} & \multirow{5}{*}{2.25e-04} & \textcolor{Red}{CG} & \textcolor{Red}{9.0e-09} & \textcolor{Red}{4.1e-08} & \textcolor{Red}{5.8e-05} & \textcolor{Red}{33} & \textcolor{Red}{0.0088} & \textcolor{Red}{335} & \textcolor{Red}{2.521} \\ 
  &  &  &  &  & \textcolor{Cerulean}{MINRES} & \textcolor{Cerulean}{7.5e-08} & \textcolor{Cerulean}{1.0e-05} & \textcolor{Cerulean}{1.4e-02} & \textcolor{Cerulean}{30} & \textcolor{Cerulean}{0.0133} & \textcolor{Cerulean}{200} & \textcolor{Cerulean}{1.594} \\ 
  &  &  &  &  & \textcolor{OliveGreen}{CG-K} & \textcolor{OliveGreen}{1.3e-08} & \textcolor{OliveGreen}{8.3e-09} & \textcolor{OliveGreen}{9.7e-06} & \textcolor{OliveGreen}{69} & \textcolor{OliveGreen}{0.0641} & \textcolor{OliveGreen}{769} & \textcolor{OliveGreen}{2.833} \\ 
  &  &  &  &  & \textcolor{Violet}{MINRES-K} & \textcolor{Violet}{2.8e-06} & \textcolor{Violet}{2.9e-06} & \textcolor{Violet}{4.1e-03} & \textcolor{Violet}{49} & \textcolor{Violet}{0.0658} & \textcolor{Violet}{469} & \textcolor{Violet}{1.922} \\ 
  &  &  &  &  & \textcolor{Magenta}{LSQR} & \textcolor{Magenta}{4.0e-07} & \textcolor{Magenta}{1.3e-07} & \textcolor{Magenta}{2.0e-04} & \textcolor{Magenta}{27} & \textcolor{Magenta}{0.0363} & \textcolor{Magenta}{309} & \textcolor{Magenta}{1.237} \\ 
\hline 
\multirow{5}{*}{1000} & \multirow{5}{*}{38400} & \multirow{5}{*}{202731} & \multirow{5}{*}{7.69e-02} & \multirow{5}{*}{1.22e-03} & \textcolor{Red}{CG} & \textcolor{Red}{5.7e-09} & \textcolor{Red}{3.0e-08} & \textcolor{Red}{5.2e-04} & \textcolor{Red}{27} & \textcolor{Red}{0.0096} & \textcolor{Red}{211} & \textcolor{Red}{20.65} \\ 
  &  &  &  &  & \textcolor{Cerulean}{MINRES} & \textcolor{Cerulean}{5.5e-08} & \textcolor{Cerulean}{1.3e-05} & \textcolor{Cerulean}{2.2e-01} & \textcolor{Cerulean}{24} & \textcolor{Cerulean}{0.0127} & \textcolor{Cerulean}{120} & \textcolor{Cerulean}{13.12} \\ 
  &  &  &  &  & \textcolor{OliveGreen}{CG-K} & \textcolor{OliveGreen}{4.2e-08} & \textcolor{OliveGreen}{5.4e-09} & \textcolor{OliveGreen}{8.3e-05} & \textcolor{OliveGreen}{55} & \textcolor{OliveGreen}{0.1219} & \textcolor{OliveGreen}{479} & \textcolor{OliveGreen}{9.466} \\ 
  &  &  &  &  & \textcolor{Violet}{MINRES-K} & \textcolor{Violet}{9.9e-06} & \textcolor{Violet}{2.8e-06} & \textcolor{Violet}{5.0e-02} & \textcolor{Violet}{35} & \textcolor{Violet}{0.0943} & \textcolor{Violet}{289} & \textcolor{Violet}{9.651} \\ 
  &  &  &  &  & \textcolor{Magenta}{LSQR} & \textcolor{Magenta}{4.8e-07} & \textcolor{Magenta}{3.8e-08} & \textcolor{Magenta}{8.9e-04} & \textcolor{Magenta}{21} & \textcolor{Magenta}{0.0421} & \textcolor{Magenta}{208} & \textcolor{Magenta}{6.108} \\ 
\hline 
\end{tabular} 

  \caption{\BA~graphs.}
  \label{tab:bakrylv} 
\end{table}

\clearpage

\section{Timing and error using algebraic
  multigrid} \label{appndx:amg}
 
The next three tables show the results of numerical experiments using
algebraic multigrid on simplicial complexes generated from \ER, \WS,
and \BA~graphs. We used the PyAMG implementation of algebraic
multigrid~\cite{BeOlSc2011}. As in the Krylov tables earlier, the
cases where the homology of the simplicial 2-complex is trivial are
marked by an asterisk ($\ast$). The cases where PyAMG failed in the
setup phase are marked by a dagger symbol ($\dagger$), and the cases
where algebraic multigrid reached maximum number of specified
iterations are marked by the double dagger symbol ($\ddagger$). In
Table~\ref{tab:eramg} the entries marked with a dash symbol (--) are
those for which Algorithm~\ref{alg:schur} (Schur complement and
algebraic multigrid algorithm, see Section~\ref{subsec:schur}) cannot
be applied because $\laplacian_2$ has no simple sparse/dense
partitioning. See Section~\ref{subsec:amg} for more details.

\bigskip

\begin{table}[h] \scriptsize
 \centering
 \begin{tabular}{|ccccc|c|ccc|cc|cc|}
\hline
\multirow{2}{*}{$N_0$} & \multirow{2}{*}{$N_1$} & \multirow{2}{*}{$N_2$} & Edge & Triangle & Algorithm / & \multicolumn{3}{|c|}{Relative Error} & \multicolumn{2}{|c|}{$\alpha$} & \multicolumn{2}{|c|}{$\beta$} \\
  &  &  & Density & Density & Formulation & $\|\boundary_1^T \alpha\|$ & $\|\boundary_2\beta\|$ & $\|h\|$ & iter. & sec. & iter. & sec. \\ 
 \hline 
\multirow{2}{*}{100} & \multirow{2}{*}{380} & \multirow{2}{*}{52} & \multirow{2}{*}{7.68e-02} & \multirow{2}{*}{3.22e-04} & \textcolor{Red}{AMG (SA)} & \textcolor{Red}{5.0e-09} & \textcolor{Red}{3.8e-09} & \textcolor{Red}{9.5e-08} & \textcolor{Red}{2} & \textcolor{Red}{0.0577} & \textcolor{Red}{2} & \textcolor{Red}{0.0030} \\ 
  &  &  &  &  & \textcolor{Cerulean}{AMG (LA)} & \textcolor{Cerulean}{5.0e-09} & \textcolor{Cerulean}{3.8e-09} & \textcolor{Cerulean}{9.5e-08} & \textcolor{Cerulean}{2} & \textcolor{Cerulean}{0.0078} & \textcolor{Cerulean}{2} & \textcolor{Cerulean}{0.0023} \\ 
  &  &  &  &  & \textcolor{OliveGreen}{Schur} & \textcolor{OliveGreen}{5.0e-09} & \textcolor{OliveGreen}{--} & \textcolor{OliveGreen}{--} & \textcolor{OliveGreen}{3} & \textcolor{OliveGreen}{0.0430} & \textcolor{OliveGreen}{--} & \textcolor{OliveGreen}{--} \\ 
\hline 
\multirow{2}{*}{100} & \multirow{2}{*}{494} & \multirow{2}{*}{144} & \multirow{2}{*}{9.98e-02} & \multirow{2}{*}{8.91e-04} & \textcolor{Red}{AMG (SA)} & \textcolor{Red}{2.8e-09} & \textcolor{Red}{3.8e-09} & \textcolor{Red}{1.0e-07} & \textcolor{Red}{2} & \textcolor{Red}{0.0078} & \textcolor{Red}{2} & \textcolor{Red}{0.0137} \\ 
  &  &  &  &  & \textcolor{Cerulean}{AMG (LA)} & \textcolor{Cerulean}{2.8e-09} & \textcolor{Cerulean}{3.8e-09} & \textcolor{Cerulean}{1.0e-07} & \textcolor{Cerulean}{2} & \textcolor{Cerulean}{0.0082} & \textcolor{Cerulean}{2} & \textcolor{Cerulean}{0.0131} \\ 
  &  &  &  &  & \textcolor{OliveGreen}{Schur} & \textcolor{OliveGreen}{2.8e-09} & \textcolor{OliveGreen}{--} & \textcolor{OliveGreen}{--} & \textcolor{OliveGreen}{3} & \textcolor{OliveGreen}{0.0228} & \textcolor{OliveGreen}{--} & \textcolor{OliveGreen}{--} \\ 
\hline 
\multirow{2}{*}{100} & \multirow{2}{*}{1212} & \multirow{2}{*}{2359} & \multirow{2}{*}{2.45e-01} & \multirow{2}{*}{1.46e-02} & \textcolor{Red}{AMG (SA)} & \textcolor{Red}{1.2e-10} & \textcolor{Red}{4.8e-08} & \textcolor{Red}{1.1e-04} & \textcolor{Red}{2} & \textcolor{Red}{0.0079} & \textcolor{Red}{32} & \textcolor{Red}{0.6259} \\ 
  &  &  &  &  & \textcolor{Cerulean}{AMG (LA)} & \textcolor{Cerulean}{1.2e-10} & \textcolor{Cerulean}{3.4e-08} & \textcolor{Cerulean}{7.5e-05} & \textcolor{Cerulean}{2} & \textcolor{Cerulean}{0.0079} & \textcolor{Cerulean}{30} & \textcolor{Cerulean}{0.1426} \\ 
  &  &  &  &  & \textcolor{OliveGreen}{Schur} & \textcolor{OliveGreen}{1.2e-10} & \textcolor{OliveGreen}{--} & \textcolor{OliveGreen}{--} & \textcolor{OliveGreen}{3} & \textcolor{OliveGreen}{0.0194} & \textcolor{OliveGreen}{--} & \textcolor{OliveGreen}{--} \\ 
\hline 
\multirow{2}{*}{100} & \multirow{2}{*}{2530} & \multirow{2}{*}{21494} & \multirow{2}{*}{5.11e-01} & \multirow{2}{*}{1.33e-01} & \textcolor{Red}{AMG (SA)} & \textcolor{Red}{6.5e-10} & \textcolor{Red}{6.9e-09} & \textcolor{Red}{2.3e-06$^{\ast}$} & \textcolor{Red}{2} & \textcolor{Red}{0.0142} & \textcolor{Red}{27} & \textcolor{Red}{2.102} \\ 
  &  &  &  &  & \textcolor{Cerulean}{AMG (LA)} & \textcolor{Cerulean}{6.5e-10} & \textcolor{Cerulean}{4.8e-09} & \textcolor{Cerulean}{1.6e-06$^{\ast}$} & \textcolor{Cerulean}{2} & \textcolor{Cerulean}{0.0184} & \textcolor{Cerulean}{19} & \textcolor{Cerulean}{15.89} \\ 
  &  &  &  &  & \textcolor{OliveGreen}{Schur} & \textcolor{OliveGreen}{6.5e-10} & \textcolor{OliveGreen}{--} & \textcolor{OliveGreen}{--} & \textcolor{OliveGreen}{3} & \textcolor{OliveGreen}{0.0185} & \textcolor{OliveGreen}{--} & \textcolor{OliveGreen}{--} \\ 
\hline 
\multirow{2}{*}{100} & \multirow{2}{*}{3706} & \multirow{2}{*}{67865} & \multirow{2}{*}{7.49e-01} & \multirow{2}{*}{4.20e-01} & \textcolor{Red}{AMG (SA)} & \textcolor{Red}{2.8e-10} & \textcolor{Red}{7.7e-09} & \textcolor{Red}{6.6e-06$^{\ast}$} & \textcolor{Red}{2} & \textcolor{Red}{0.0689} & \textcolor{Red}{39} & \textcolor{Red}{16.74} \\ 
  &  &  &  &  & \textcolor{Cerulean}{AMG (LA)} & \textcolor{Cerulean}{2.8e-10} & \textcolor{Cerulean}{7.3e-09} & \textcolor{Cerulean}{6.2e-06$^{\ast}$} & \textcolor{Cerulean}{2} & \textcolor{Cerulean}{0.0356} & \textcolor{Cerulean}{14} & \textcolor{Cerulean}{408.5} \\ 
  &  &  &  &  & \textcolor{OliveGreen}{Schur} & \textcolor{OliveGreen}{2.8e-10} & \textcolor{OliveGreen}{--} & \textcolor{OliveGreen}{--} & \textcolor{OliveGreen}{3} & \textcolor{OliveGreen}{0.3740} & \textcolor{OliveGreen}{--} & \textcolor{OliveGreen}{--} \\ 
\hline 
\multirow{2}{*}{500} & \multirow{2}{*}{1290} & \multirow{2}{*}{21} & \multirow{2}{*}{1.03e-02} & \multirow{2}{*}{1.01e-06} & \textcolor{Red}{AMG (SA)} & \textcolor{Red}{3.8e-09} & \textcolor{Red}{1.9e-09} & \textcolor{Red}{9.8e-08} & \textcolor{Red}{2} & \textcolor{Red}{0.7113} & \textcolor{Red}{2} & \textcolor{Red}{0.0360} \\ 
  &  &  &  &  & \textcolor{Cerulean}{AMG (LA)} & \textcolor{Cerulean}{3.8e-09} & \textcolor{Cerulean}{1.9e-09} & \textcolor{Cerulean}{9.8e-08} & \textcolor{Cerulean}{2} & \textcolor{Cerulean}{0.3080} & \textcolor{Cerulean}{2} & \textcolor{Cerulean}{0.0022} \\ 
  &  &  &  &  & \textcolor{OliveGreen}{Schur} & \textcolor{OliveGreen}{3.8e-09} & \textcolor{OliveGreen}{--} & \textcolor{OliveGreen}{--} & \textcolor{OliveGreen}{3} & \textcolor{OliveGreen}{0.2952} & \textcolor{OliveGreen}{--} & \textcolor{OliveGreen}{--} \\ 
\hline 
\multirow{2}{*}{500} & \multirow{2}{*}{12394} & \multirow{2}{*}{20315} & \multirow{2}{*}{9.94e-02} & \multirow{2}{*}{9.81e-04} & \textcolor{Red}{AMG (SA)} & \textcolor{Red}{7.2e-11} & \textcolor{Red}{7.7e-08} & \textcolor{Red}{2.5e-04} & \textcolor{Red}{2} & \textcolor{Red}{0.3305} & \textcolor{Red}{74} & \textcolor{Red}{2.711} \\ 
  &  &  &  &  & \textcolor{Cerulean}{AMG (LA)} & \textcolor{Cerulean}{7.2e-11} & \textcolor{Cerulean}{7.4e-08} & \textcolor{Cerulean}{2.5e-04} & \textcolor{Cerulean}{2} & \textcolor{Cerulean}{0.3293} & \textcolor{Cerulean}{71} & \textcolor{Cerulean}{8.358} \\ 
  &  &  &  &  & \textcolor{OliveGreen}{Schur} & \textcolor{OliveGreen}{7.2e-11} & \textcolor{OliveGreen}{--} & \textcolor{OliveGreen}{--} & \textcolor{OliveGreen}{3} & \textcolor{OliveGreen}{0.2836} & \textcolor{OliveGreen}{--} & \textcolor{OliveGreen}{--} \\ 
\hline 
\multirow{2}{*}{500} & \multirow{2}{*}{24788} & \multirow{2}{*}{162986} & \multirow{2}{*}{1.99e-01} & \multirow{2}{*}{7.87e-03} & \textcolor{Red}{AMG (SA)} & \textcolor{Red}{$\dagger$} & \textcolor{Red}{$\dagger$} & \textcolor{Red}{$\dagger$} & \textcolor{Red}{$\dagger$} & \textcolor{Red}{$\dagger$} & \textcolor{Red}{$\dagger$} & \textcolor{Red}{$\dagger$} \\ 
  &  &  &  &  & \textcolor{Cerulean}{AMG (LA)} & \textcolor{Cerulean}{$\dagger$} & \textcolor{Cerulean}{$\dagger$} & \textcolor{Cerulean}{$\dagger$} & \textcolor{Cerulean}{$\dagger$} & \textcolor{Cerulean}{$\dagger$} & \textcolor{Cerulean}{$\dagger$} & \textcolor{Cerulean}{$\dagger$} \\ 
  &  &  &  &  & \textcolor{OliveGreen}{Schur} & \textcolor{OliveGreen}{$\dagger$} & \textcolor{OliveGreen}{--} & \textcolor{OliveGreen}{--} & \textcolor{OliveGreen}{$\dagger$} & \textcolor{OliveGreen}{$\dagger$} & \textcolor{OliveGreen}{--} & \textcolor{OliveGreen}{--} \\ 
\hline
\multirow{2}{*}{1000} & \multirow{2}{*}{49690} & \multirow{2}{*}{163767} & \multirow{2}{*}{9.95e-02} & \multirow{2}{*}{9.86e-04} & \textcolor{Red}{AMG (SA)} & \textcolor{Red}{$\dagger$} & \textcolor{Red}{$\dagger$} & \textcolor{Red}{$\dagger$} & \textcolor{Red}{$\dagger$} & \textcolor{Red}{$\dagger$} & \textcolor{Red}{$\dagger$} & \textcolor{Red}{$\dagger$} \\ 
  &  &  &  &  & \textcolor{Cerulean}{AMG (LA)} & \textcolor{Cerulean}{$\dagger$} & \textcolor{Cerulean}{$\dagger$} & \textcolor{Cerulean}{$\dagger$} & \textcolor{Cerulean}{$\dagger$} & \textcolor{Cerulean}{$\dagger$} & \textcolor{Cerulean}{$\dagger$} & \textcolor{Cerulean}{$\dagger$} \\ 
  &  &  &  &  & \textcolor{OliveGreen}{Schur} & \textcolor{OliveGreen}{$\dagger$} & \textcolor{OliveGreen}{--} & \textcolor{OliveGreen}{--} & \textcolor{OliveGreen}{$\dagger$} & \textcolor{OliveGreen}{$\dagger$} & \textcolor{OliveGreen}{--} & \textcolor{OliveGreen}{--} \\ 
\hline
\end{tabular} 

 \caption{\ER~graphs.} 
\label{tab:eramg}
\end{table}

\begin{table}[H] \scriptsize
 \centering
 \begin{tabular}{|ccccc|c|ccc|cc|cc|}
\hline
\multirow{2}{*}{$N_0$} & \multirow{2}{*}{$N_1$} & \multirow{2}{*}{$N_2$} & Edge & Triangle & Algorithm / & \multicolumn{3}{|c|}{Relative Error} & \multicolumn{2}{|c|}{$\alpha$} & \multicolumn{2}{|c|}{$\beta$} \\
  &  &  & Density & Density & Formulation & $\|\boundary_1^T \alpha\|$ & $\|\boundary_2\beta\|$ & $\|h\|$ & iter. & sec. & iter. & sec. \\ 
 \hline 
\multirow{2}{*}{100} & \multirow{2}{*}{500} & \multirow{2}{*}{729} & \multirow{2}{*}{1.01e-01} & \multirow{2}{*}{4.51e-03} & \textcolor{Red}{AMG (SA)} & \textcolor{Red}{2.9e-09} & \textcolor{Red}{2.4e-08} & \textcolor{Red}{1.2e-06} & \textcolor{Red}{2} & \textcolor{Red}{0.0179} & \textcolor{Red}{20} & \textcolor{Red}{0.0577} \\ 
  &  &  &  &  & \textcolor{Cerulean}{AMG (LA)} & \textcolor{Cerulean}{2.9e-09} & \textcolor{Cerulean}{1.4e-08} & \textcolor{Cerulean}{6.9e-07} & \textcolor{Cerulean}{2} & \textcolor{Cerulean}{0.0009} & \textcolor{Cerulean}{18} & \textcolor{Cerulean}{0.0324} \\ 
\hline 
\multirow{2}{*}{100} & \multirow{2}{*}{1000} & \multirow{2}{*}{3655} & \multirow{2}{*}{2.02e-01} & \multirow{2}{*}{2.26e-02} & \textcolor{Red}{AMG (SA)} & \textcolor{Red}{1.7e-09} & \textcolor{Red}{3.1e-08} & \textcolor{Red}{4.6e-06} & \textcolor{Red}{2} & \textcolor{Red}{0.0009} & \textcolor{Red}{23} & \textcolor{Red}{0.1253} \\ 
  &  &  &  &  & \textcolor{Cerulean}{AMG (LA)} & \textcolor{Cerulean}{1.7e-09} & \textcolor{Cerulean}{2.7e-08} & \textcolor{Cerulean}{4.0e-06} & \textcolor{Cerulean}{2} & \textcolor{Cerulean}{0.0009} & \textcolor{Cerulean}{17} & \textcolor{Cerulean}{0.1601} \\ 
\hline 
\multirow{2}{*}{100} & \multirow{2}{*}{1500} & \multirow{2}{*}{8354} & \multirow{2}{*}{3.03e-01} & \multirow{2}{*}{5.17e-02} & \textcolor{Red}{AMG (SA)} & \textcolor{Red}{1.3e-10} & \textcolor{Red}{7.6e-08} & \textcolor{Red}{3.0e-04} & \textcolor{Red}{2} & \textcolor{Red}{0.0011} & \textcolor{Red}{99} & \textcolor{Red}{1.377} \\ 
  &  &  &  &  & \textcolor{Cerulean}{AMG (LA)} & \textcolor{Cerulean}{1.3e-10} & \textcolor{Cerulean}{7.1e-08} & \textcolor{Cerulean}{2.8e-04} & \textcolor{Cerulean}{2} & \textcolor{Cerulean}{0.0011} & \textcolor{Cerulean}{60} & \textcolor{Cerulean}{2.345} \\ 
\hline 
\multirow{2}{*}{100} & \multirow{2}{*}{2000} & \multirow{2}{*}{15530} & \multirow{2}{*}{4.04e-01} & \multirow{2}{*}{9.60e-02} & \textcolor{Red}{AMG (SA)} & \textcolor{Red}{3.8e-09} & \textcolor{Red}{4.6e-08} & \textcolor{Red}{1.1e-05$^{\ast}$} & \textcolor{Red}{2} & \textcolor{Red}{0.0015} & \textcolor{Red}{28} & \textcolor{Red}{1.282} \\ 
  &  &  &  &  & \textcolor{Cerulean}{AMG (LA)} & \textcolor{Cerulean}{3.8e-09} & \textcolor{Cerulean}{1.8e-08} & \textcolor{Cerulean}{4.3e-06$^{\ast}$} & \textcolor{Cerulean}{2} & \textcolor{Cerulean}{0.0013} & \textcolor{Cerulean}{15} & \textcolor{Cerulean}{5.039} \\ 
\hline 
\multirow{2}{*}{500} & \multirow{2}{*}{2500} & \multirow{2}{*}{3720} & \multirow{2}{*}{2.00e-02} & \multirow{2}{*}{1.80e-04} & \textcolor{Red}{AMG (SA)} & \textcolor{Red}{2.5e-09} & \textcolor{Red}{3.0e-08} & \textcolor{Red}{3.0e-06} & \textcolor{Red}{2} & \textcolor{Red}{0.0149} & \textcolor{Red}{15} & \textcolor{Red}{0.0650} \\ 
  &  &  &  &  & \textcolor{Cerulean}{AMG (LA)} & \textcolor{Cerulean}{2.5e-09} & \textcolor{Cerulean}{4.4e-08} & \textcolor{Cerulean}{4.3e-06} & \textcolor{Cerulean}{2} & \textcolor{Cerulean}{0.0136} & \textcolor{Cerulean}{21} & \textcolor{Cerulean}{0.0929} \\ 
\hline 
\multirow{2}{*}{500} & \multirow{2}{*}{5000} & \multirow{2}{*}{16948} & \multirow{2}{*}{4.01e-02} & \multirow{2}{*}{8.18e-04} & \textcolor{Red}{AMG (SA)} & \textcolor{Red}{6.0e-10} & \textcolor{Red}{6.9e-09} & \textcolor{Red}{1.9e-06} & \textcolor{Red}{2} & \textcolor{Red}{0.0139} & \textcolor{Red}{14} & \textcolor{Red}{0.4410} \\ 
  &  &  &  &  & \textcolor{Cerulean}{AMG (LA)} & \textcolor{Cerulean}{6.0e-10} & \textcolor{Cerulean}{3.9e-09} & \textcolor{Cerulean}{1.1e-06} & \textcolor{Cerulean}{2} & \textcolor{Cerulean}{0.0140} & \textcolor{Cerulean}{12} & \textcolor{Cerulean}{0.7969} \\ 
\hline 
\multirow{2}{*}{500} & \multirow{2}{*}{12500} & \multirow{2}{*}{110507} & \multirow{2}{*}{1.00e-01} & \multirow{2}{*}{5.34e-03} & \textcolor{Red}{AMG (SA)} & \textcolor{Red}{8.8e-10} & \textcolor{Red}{3.9e-07} & \textcolor{Red}{4.4e-04} & \textcolor{Red}{2} & \textcolor{Red}{0.0143} & \textcolor{Red}{101$^{\ddagger}$} & \textcolor{Red}{32.75} \\ 
  &  &  &  &  & \textcolor{Cerulean}{AMG (LA)} & \textcolor{Cerulean}{8.8e-10} & \textcolor{Cerulean}{7.3e-08} & \textcolor{Cerulean}{8.2e-05} & \textcolor{Cerulean}{2} & \textcolor{Cerulean}{0.0145} & \textcolor{Cerulean}{86} & \textcolor{Cerulean}{173.1} \\ 
\hline 
\multirow{2}{*}{1000} & \multirow{2}{*}{5000} & \multirow{2}{*}{7386} & \multirow{2}{*}{1.00e-02} & \multirow{2}{*}{4.44e-05} & \textcolor{Red}{AMG (SA)} & \textcolor{Red}{2.7e-08} & \textcolor{Red}{2.7e-08} & \textcolor{Red}{3.8e-06} & \textcolor{Red}{14} & \textcolor{Red}{0.0651} & \textcolor{Red}{20} & \textcolor{Red}{0.1325} \\ 
  &  &  &  &  & \textcolor{Cerulean}{AMG (LA)} & \textcolor{Cerulean}{2.0e-08} & \textcolor{Cerulean}{2.5e-08} & \textcolor{Cerulean}{3.5e-06} & \textcolor{Cerulean}{17} & \textcolor{Cerulean}{0.0646} & \textcolor{Cerulean}{19} & \textcolor{Cerulean}{0.1715} \\ 
\hline 
\multirow{2}{*}{1000} & \multirow{2}{*}{10000} & \multirow{2}{*}{33022} & \multirow{2}{*}{2.00e-02} & \multirow{2}{*}{1.99e-04} & \textcolor{Red}{AMG (SA)} & \textcolor{Red}{1.8e-08} & \textcolor{Red}{1.6e-08} & \textcolor{Red}{6.0e-06} & \textcolor{Red}{17} & \textcolor{Red}{0.0697} & \textcolor{Red}{17} & \textcolor{Red}{1.058} \\ 
  &  &  &  &  & \textcolor{Cerulean}{AMG (LA)} & \textcolor{Cerulean}{2.4e-08} & \textcolor{Cerulean}{1.3e-08} & \textcolor{Cerulean}{5.1e-06} & \textcolor{Cerulean}{10} & \textcolor{Cerulean}{0.0715} & \textcolor{Cerulean}{13} & \textcolor{Cerulean}{1.798} \\ 
\hline 
\multirow{2}{*}{1000} & \multirow{2}{*}{25000} & \multirow{2}{*}{220002} & \multirow{2}{*}{5.01e-02} & \multirow{2}{*}{1.32e-03} & \textcolor{Red}{AMG (SA)} & \textcolor{Red}{1.6e-08} & \textcolor{Red}{2.9e-08} & \textcolor{Red}{4.1e-05} & \textcolor{Red}{15} & \textcolor{Red}{0.0823} & \textcolor{Red}{29} & \textcolor{Red}{26.71} \\ 
  &  &  &  &  & \textcolor{Cerulean}{AMG (LA)} & \textcolor{Cerulean}{5.2e-09} & \textcolor{Cerulean}{5.3e-08} & \textcolor{Cerulean}{7.6e-05} & \textcolor{Cerulean}{6} & \textcolor{Cerulean}{0.1034} & \textcolor{Cerulean}{23} & \textcolor{Cerulean}{241.3} \\ 
\hline 
\end{tabular} 

 \caption{\WS~graphs.}
 \label{tab:wsamg} 
\end{table}

\begin{table}[h] \scriptsize
  \centering
  \begin{tabular}{|ccccc|c|ccc|cc|cc|}
\hline
\multirow{2}{*}{$N_0$} & \multirow{2}{*}{$N_1$} & \multirow{2}{*}{$N_2$} & Edge & Triangle & Algorithm / & \multicolumn{3}{|c|}{Relative Error} & \multicolumn{2}{|c|}{$\alpha$} & \multicolumn{2}{|c|}{$\beta$} \\
  &  &  & Density & Density & Formulation & $\|\boundary_1^T \alpha\|$ & $\|\boundary_2\beta\|$ & $\|h\|$ & iter. & sec. & iter. & sec. \\ 
 \hline 
\multirow{2}{*}{100} & \multirow{2}{*}{475} & \multirow{2}{*}{301} & \multirow{2}{*}{9.60e-02} & \multirow{2}{*}{1.86e-03} & \textcolor{Red}{AMG (SA)} & \textcolor{Red}{3.2e-09} & \textcolor{Red}{4.3e-09} & \textcolor{Red}{2.4e-07} & \textcolor{Red}{2} & \textcolor{Red}{0.0009} & \textcolor{Red}{2} & \textcolor{Red}{0.0038} \\ 
  &  &  &  &  & \textcolor{Cerulean}{AMG (LA)} & \textcolor{Cerulean}{3.2e-09} & \textcolor{Cerulean}{4.3e-09} & \textcolor{Cerulean}{2.4e-07} & \textcolor{Cerulean}{2} & \textcolor{Cerulean}{0.0011} & \textcolor{Cerulean}{2} & \textcolor{Cerulean}{0.0037} \\ 
\hline 
\multirow{2}{*}{100} & \multirow{2}{*}{900} & \multirow{2}{*}{1701} & \multirow{2}{*}{1.82e-01} & \multirow{2}{*}{1.05e-02} & \textcolor{Red}{AMG (SA)} & \textcolor{Red}{6.7e-10} & \textcolor{Red}{4.7e-08} & \textcolor{Red}{1.8e-05} & \textcolor{Red}{2} & \textcolor{Red}{0.0016} & \textcolor{Red}{50} & \textcolor{Red}{0.0995} \\ 
  &  &  &  &  & \textcolor{Cerulean}{AMG (LA)} & \textcolor{Cerulean}{6.7e-10} & \textcolor{Cerulean}{5.6e-08} & \textcolor{Cerulean}{2.2e-05} & \textcolor{Cerulean}{2} & \textcolor{Cerulean}{0.0009} & \textcolor{Cerulean}{44} & \textcolor{Cerulean}{0.1050} \\ 
\hline 
\multirow{2}{*}{100} & \multirow{2}{*}{1600} & \multirow{2}{*}{8105} & \multirow{2}{*}{3.23e-01} & \multirow{2}{*}{5.01e-02} & \textcolor{Red}{AMG (SA)} & \textcolor{Red}{2.0e-09} & \textcolor{Red}{1.2e-08} & \textcolor{Red}{1.9e-06$^{\ast}$} & \textcolor{Red}{2} & \textcolor{Red}{0.0009} & \textcolor{Red}{14} & \textcolor{Red}{0.3561} \\ 
  &  &  &  &  & \textcolor{Cerulean}{AMG (LA)} & \textcolor{Cerulean}{2.0e-09} & \textcolor{Cerulean}{5.9e-09} & \textcolor{Cerulean}{9.8e-07$^{\ast}$} & \textcolor{Cerulean}{2} & \textcolor{Cerulean}{0.0010} & \textcolor{Cerulean}{12} & \textcolor{Cerulean}{1.545} \\ 
\hline 
\multirow{2}{*}{100} & \multirow{2}{*}{2400} & \multirow{2}{*}{24497} & \multirow{2}{*}{4.85e-01} & \multirow{2}{*}{1.51e-01} & \textcolor{Red}{AMG (SA)} & \textcolor{Red}{1.1e-09} & \textcolor{Red}{7.3e-09} & \textcolor{Red}{2.6e-06$^{\ast}$} & \textcolor{Red}{2} & \textcolor{Red}{0.0017} & \textcolor{Red}{24} & \textcolor{Red}{2.632} \\ 
  &  &  &  &  & \textcolor{Cerulean}{AMG (LA)} & \textcolor{Cerulean}{1.1e-09} & \textcolor{Cerulean}{7.6e-09} & \textcolor{Cerulean}{2.7e-06$^{\ast}$} & \textcolor{Cerulean}{2} & \textcolor{Cerulean}{0.0013} & \textcolor{Cerulean}{13} & \textcolor{Cerulean}{22.41} \\ 
\hline 
\multirow{2}{*}{500} & \multirow{2}{*}{4900} & \multirow{2}{*}{4740} & \multirow{2}{*}{3.93e-02} & \multirow{2}{*}{2.29e-04} & \textcolor{Red}{AMG (SA)} & \textcolor{Red}{2.3e-09} & \textcolor{Red}{7.5e-08} & \textcolor{Red}{2.0e-05} & \textcolor{Red}{2} & \textcolor{Red}{0.0169} & \textcolor{Red}{89} & \textcolor{Red}{0.4003} \\ 
  &  &  &  &  & \textcolor{Cerulean}{AMG (LA)} & \textcolor{Cerulean}{2.3e-09} & \textcolor{Cerulean}{8.9e-08} & \textcolor{Cerulean}{2.4e-05} & \textcolor{Cerulean}{2} & \textcolor{Cerulean}{0.0137} & \textcolor{Cerulean}{95} & \textcolor{Cerulean}{0.5046} \\ 
\hline 
\multirow{2}{*}{500} & \multirow{2}{*}{9600} & \multirow{2}{*}{25016} & \multirow{2}{*}{7.70e-02} & \multirow{2}{*}{1.21e-03} & \textcolor{Red}{AMG (SA)} & \textcolor{Red}{4.9e-10} & \textcolor{Red}{7.3e-07} & \textcolor{Red}{1.4e-03} & \textcolor{Red}{2} & \textcolor{Red}{0.0142} & \textcolor{Red}{101$^{\ddagger}$} & \textcolor{Red}{4.627} \\ 
  &  &  &  &  & \textcolor{Cerulean}{AMG (LA)} & \textcolor{Cerulean}{4.9e-10} & \textcolor{Cerulean}{4.0e-07} & \textcolor{Cerulean}{7.8e-04} & \textcolor{Cerulean}{2} & \textcolor{Cerulean}{0.0142} & \textcolor{Cerulean}{101$^{\ddagger}$} & \textcolor{Cerulean}{29.74} \\ 
\hline 
\multirow{2}{*}{500} & \multirow{2}{*}{18400} & \multirow{2}{*}{133933} & \multirow{2}{*}{1.47e-01} & \multirow{2}{*}{6.47e-03} & \textcolor{Red}{AMG (SA)} & \textcolor{Red}{$\dagger$} & \textcolor{Red}{$\dagger$} & \textcolor{Red}{$\dagger$} & \textcolor{Red}{$\dagger$} & \textcolor{Red}{$\dagger$} & \textcolor{Red}{$\dagger$} & \textcolor{Red}{$\dagger$} \\ 
  &  &  &  &  & \textcolor{Cerulean}{AMG (LA)} & \textcolor{Cerulean}{$\dagger$} & \textcolor{Cerulean}{$\dagger$} & \textcolor{Cerulean}{$\dagger$} & \textcolor{Cerulean}{$\dagger$} & \textcolor{Cerulean}{$\dagger$} & \textcolor{Cerulean}{$\dagger$} & \textcolor{Cerulean}{$\dagger$} \\ 
\hline 
\multirow{2}{*}{1000} & \multirow{2}{*}{9900} & \multirow{2}{*}{6264} & \multirow{2}{*}{1.98e-02} & \multirow{2}{*}{3.77e-05} & \textcolor{Red}{AMG (SA)} & \textcolor{Red}{5.3e-09} & \textcolor{Red}{5.6e-06} & \textcolor{Red}{1.4e-03} & \textcolor{Red}{5} & \textcolor{Red}{0.0477} & \textcolor{Red}{101$^{\ddagger}$} & \textcolor{Red}{0.5723} \\ 
  &  &  &  &  & \textcolor{Cerulean}{AMG (LA)} & \textcolor{Cerulean}{3.5e-09} & \textcolor{Cerulean}{9.1e-06} & \textcolor{Cerulean}{2.3e-03} & \textcolor{Cerulean}{5} & \textcolor{Cerulean}{0.0682} & \textcolor{Cerulean}{101$^{\ddagger}$} & \textcolor{Cerulean}{1.22} \\ 
\hline 
\multirow{2}{*}{1000} & \multirow{2}{*}{19600} & \multirow{2}{*}{37365} & \multirow{2}{*}{3.92e-02} & \multirow{2}{*}{2.25e-04} & \textcolor{Red}{AMG (SA)} & \textcolor{Red}{7.5e-09} & \textcolor{Red}{1.1e-05} & \textcolor{Red}{1.5e-02} & \textcolor{Red}{4} & \textcolor{Red}{0.0393} & \textcolor{Red}{101$^{\ddagger}$} & \textcolor{Red}{10.61} \\ 
  &  &  &  &  & \textcolor{Cerulean}{AMG (LA)} & \textcolor{Cerulean}{5.6e-09} & \textcolor{Cerulean}{8.6e-06} & \textcolor{Cerulean}{1.2e-02} & \textcolor{Cerulean}{4} & \textcolor{Cerulean}{0.0823} & \textcolor{Cerulean}{101$^{\ddagger}$} & \textcolor{Cerulean}{63.01} \\ 
\hline 
\multirow{2}{*}{1000} & \multirow{2}{*}{38400} & \multirow{2}{*}{202731} & \multirow{2}{*}{7.69e-02} & \multirow{2}{*}{1.22e-03} & \textcolor{Red}{AMG (SA)} & \textcolor{Red}{$\dagger$} & \textcolor{Red}{$\dagger$} & \textcolor{Red}{$\dagger$} & \textcolor{Red}{$\dagger$} & \textcolor{Red}{$\dagger$} & \textcolor{Red}{$\dagger$} & \textcolor{Red}{$\dagger$} \\ 
  &  &  &  &  & \textcolor{Cerulean}{AMG (LA)} & \textcolor{Cerulean}{$\dagger$} & \textcolor{Cerulean}{$\dagger$} & \textcolor{Cerulean}{$\dagger$} & \textcolor{Cerulean}{$\dagger$} & \textcolor{Cerulean}{$\dagger$} & \textcolor{Cerulean}{$\dagger$} & \textcolor{Cerulean}{$\dagger$} \\ 
\hline 
\end{tabular} 

  \caption{\BA~graphs.}
  \label{tab:baamg} 
\end{table}
\clearpage

\section{Setup and solve timings for algebraic multigrid}
\label{appndx:breakup}

The setup phase of algebraic multigrid involves some
precomputation. In the next three tables we separate the timing
results shown in Tables~\ref{tab:eramg}--\ref{tab:baamg} into the time
required for the setup and solve phases.

\bigskip

\begin{table}[h] \scriptsize
 \centering
 \begin{tabular}{|ccccc|c|cccc|cccc|}
\hline
\multirow{2}{*}{$N_0$} & \multirow{2}{*}{$N_1$} & \multirow{2}{*}{$N_2$} & Edge & Triangle & Algorithm / & \multicolumn{4}{|c|}{$\alpha$} & \multicolumn{4}{|c|}{$\beta$} \\
  &  &  & Density & Density & Formulation &  iter. & setup & solve & total &  iter. & setup & solve & total \\ 
 \hline 
\multirow{2}{*}{100} & \multirow{2}{*}{380} & \multirow{2}{*}{52} & \multirow{2}{*}{7.68e-02} & \multirow{2}{*}{3.22e-04} & \textcolor{Red}{SA} & \textcolor{Red}{2} & \textcolor{Red}{0.0001} & \textcolor{Red}{0.0576} & \textcolor{Red}{0.0577} & \textcolor{Red}{2} & \textcolor{Red}{0.0001} & \textcolor{Red}{0.0029} & \textcolor{Red}{0.0030} \\ 
  &  &  &  &  & \textcolor{Cerulean}{LA} & \textcolor{Cerulean}{2} & \textcolor{Cerulean}{0.0001} & \textcolor{Cerulean}{0.0077} & \textcolor{Cerulean}{0.0078} & \textcolor{Cerulean}{2} & \textcolor{Cerulean}{0.0001} & \textcolor{Cerulean}{0.0022} & \textcolor{Cerulean}{0.0023} \\ 
\hline 
\multirow{2}{*}{100} & \multirow{2}{*}{494} & \multirow{2}{*}{144} & \multirow{2}{*}{9.98e-02} & \multirow{2}{*}{8.91e-04} & \textcolor{Red}{SA} & \textcolor{Red}{2} & \textcolor{Red}{0.0001} & \textcolor{Red}{0.0076} & \textcolor{Red}{0.0078} & \textcolor{Red}{2} & \textcolor{Red}{0.0001} & \textcolor{Red}{0.0136} & \textcolor{Red}{0.0137} \\ 
  &  &  &  &  & \textcolor{Cerulean}{LA} & \textcolor{Cerulean}{2} & \textcolor{Cerulean}{0.0001} & \textcolor{Cerulean}{0.0080} & \textcolor{Cerulean}{0.0082} & \textcolor{Cerulean}{2} & \textcolor{Cerulean}{0.0001} & \textcolor{Cerulean}{0.0129} & \textcolor{Cerulean}{0.0131} \\ 
\hline 
\multirow{2}{*}{100} & \multirow{2}{*}{1212} & \multirow{2}{*}{2359} & \multirow{2}{*}{2.45e-01} & \multirow{2}{*}{1.46e-02} & \textcolor{Red}{SA} & \textcolor{Red}{2} & \textcolor{Red}{0.0001} & \textcolor{Red}{0.0078} & \textcolor{Red}{0.0079} & \textcolor{Red}{32} & \textcolor{Red}{0.5708} & \textcolor{Red}{0.0551} & \textcolor{Red}{0.6259} \\ 
  &  &  &  &  & \textcolor{Cerulean}{LA} & \textcolor{Cerulean}{2} & \textcolor{Cerulean}{0.0002} & \textcolor{Cerulean}{0.0077} & \textcolor{Cerulean}{0.0079} & \textcolor{Cerulean}{30} & \textcolor{Cerulean}{0.0392} & \textcolor{Cerulean}{0.1034} & \textcolor{Cerulean}{0.1426} \\ 
\hline 
\multirow{2}{*}{100} & \multirow{2}{*}{2530} & \multirow{2}{*}{21494} & \multirow{2}{*}{5.11e-01} & \multirow{2}{*}{1.33e-01} & \textcolor{Red}{SA} & \textcolor{Red}{2} & \textcolor{Red}{0.0030} & \textcolor{Red}{0.0111} & \textcolor{Red}{0.0142} & \textcolor{Red}{27} & \textcolor{Red}{0.8655} & \textcolor{Red}{1.236} & \textcolor{Red}{2.102} \\ 
  &  &  &  &  & \textcolor{Cerulean}{LA} & \textcolor{Cerulean}{2} & \textcolor{Cerulean}{0.0048} & \textcolor{Cerulean}{0.0136} & \textcolor{Cerulean}{0.0184} & \textcolor{Cerulean}{19} & \textcolor{Cerulean}{12.36} & \textcolor{Cerulean}{3.528} & \textcolor{Cerulean}{15.89} \\ 
\hline 
\multirow{2}{*}{100} & \multirow{2}{*}{3706} & \multirow{2}{*}{67865} & \multirow{2}{*}{7.49e-01} & \multirow{2}{*}{4.20e-01} & \textcolor{Red}{SA} & \textcolor{Red}{2} & \textcolor{Red}{0.0303} & \textcolor{Red}{0.0386} & \textcolor{Red}{0.0689} & \textcolor{Red}{39} & \textcolor{Red}{5.652} & \textcolor{Red}{11.08} & \textcolor{Red}{16.74} \\ 
  &  &  &  &  & \textcolor{Cerulean}{LA} & \textcolor{Cerulean}{2} & \textcolor{Cerulean}{0.0138} & \textcolor{Cerulean}{0.0219} & \textcolor{Cerulean}{0.0356} & \textcolor{Cerulean}{14} & \textcolor{Cerulean}{375.7} & \textcolor{Cerulean}{32.85} & \textcolor{Cerulean}{408.5} \\ 
\hline 
\multirow{2}{*}{500} & \multirow{2}{*}{1290} & \multirow{2}{*}{21} & \multirow{2}{*}{1.03e-02} & \multirow{2}{*}{1.01e-06} & \textcolor{Red}{SA} & \textcolor{Red}{2} & \textcolor{Red}{0.1353} & \textcolor{Red}{0.5760} & \textcolor{Red}{0.7113} & \textcolor{Red}{2} & \textcolor{Red}{0.0071} & \textcolor{Red}{0.0289} & \textcolor{Red}{0.0360} \\ 
  &  &  &  &  & \textcolor{Cerulean}{LA} & \textcolor{Cerulean}{2} & \textcolor{Cerulean}{0.0001} & \textcolor{Cerulean}{0.3079} & \textcolor{Cerulean}{0.3080} & \textcolor{Cerulean}{2} & \textcolor{Cerulean}{0.0007} & \textcolor{Cerulean}{0.0015} & \textcolor{Cerulean}{0.0022} \\ 
\hline 
\multirow{2}{*}{500} & \multirow{2}{*}{12394} & \multirow{2}{*}{20315} & \multirow{2}{*}{9.94e-02} & \multirow{2}{*}{9.81e-04} & \textcolor{Red}{SA} & \textcolor{Red}{2} & \textcolor{Red}{0.0001} & \textcolor{Red}{0.3303} & \textcolor{Red}{0.3305} & \textcolor{Red}{74} & \textcolor{Red}{0.4907} & \textcolor{Red}{2.22} & \textcolor{Red}{2.711} \\ 
  &  &  &  &  & \textcolor{Cerulean}{LA} & \textcolor{Cerulean}{2} & \textcolor{Cerulean}{0.0041} & \textcolor{Cerulean}{0.3252} & \textcolor{Cerulean}{0.3293} & \textcolor{Cerulean}{71} & \textcolor{Cerulean}{3.082} & \textcolor{Cerulean}{5.277} & \textcolor{Cerulean}{8.358} \\ 
\hline 
\end{tabular} 

 \caption{\ER~graphs.}
 \label{tab:eramgsttstcs}
\end{table}

\begin{table}[h] \scriptsize
 \centering
 \begin{tabular}{|ccccc|c|cccc|cccc|}
\hline
\multirow{2}{*}{$N_0$} & \multirow{2}{*}{$N_1$} & \multirow{2}{*}{$N_2$} & Edge & Triangle & Algorithm / & \multicolumn{4}{|c|}{$\alpha$} & \multicolumn{4}{|c|}{$\beta$} \\
  &  &  & Density & Density & Formulation &  iter. & setup & solve & total &  iter. & setup & solve & total \\ 
 \hline 
\multirow{2}{*}{100} & \multirow{2}{*}{500} & \multirow{2}{*}{729} & \multirow{2}{*}{1.01e-01} & \multirow{2}{*}{4.51e-03} & \textcolor{Red}{SA} & \textcolor{Red}{2} & \textcolor{Red}{0.0001} & \textcolor{Red}{0.0178} & \textcolor{Red}{0.0179} & \textcolor{Red}{20} & \textcolor{Red}{0.0410} & \textcolor{Red}{0.0167} & \textcolor{Red}{0.0577} \\ 
  &  &  &  &  & \textcolor{Cerulean}{LA} & \textcolor{Cerulean}{2} & \textcolor{Cerulean}{0.0001} & \textcolor{Cerulean}{0.0008} & \textcolor{Cerulean}{0.0009} & \textcolor{Cerulean}{18} & \textcolor{Cerulean}{0.0172} & \textcolor{Cerulean}{0.0152} & \textcolor{Cerulean}{0.0324} \\ 
\hline 
\multirow{2}{*}{100} & \multirow{2}{*}{1000} & \multirow{2}{*}{3655} & \multirow{2}{*}{2.02e-01} & \multirow{2}{*}{2.26e-02} & \textcolor{Red}{SA} & \textcolor{Red}{2} & \textcolor{Red}{0.0001} & \textcolor{Red}{0.0008} & \textcolor{Red}{0.0009} & \textcolor{Red}{23} & \textcolor{Red}{0.0553} & \textcolor{Red}{0.0700} & \textcolor{Red}{0.1253} \\ 
  &  &  &  &  & \textcolor{Cerulean}{LA} & \textcolor{Cerulean}{2} & \textcolor{Cerulean}{0.0001} & \textcolor{Cerulean}{0.0008} & \textcolor{Cerulean}{0.0009} & \textcolor{Cerulean}{17} & \textcolor{Cerulean}{0.0939} & \textcolor{Cerulean}{0.0663} & \textcolor{Cerulean}{0.1601} \\ 
\hline 
\multirow{2}{*}{100} & \multirow{2}{*}{1500} & \multirow{2}{*}{8354} & \multirow{2}{*}{3.03e-01} & \multirow{2}{*}{5.17e-02} & \textcolor{Red}{SA} & \textcolor{Red}{2} & \textcolor{Red}{0.0002} & \textcolor{Red}{0.0009} & \textcolor{Red}{0.0011} & \textcolor{Red}{99} & \textcolor{Red}{0.1754} & \textcolor{Red}{1.202} & \textcolor{Red}{1.377} \\ 
  &  &  &  &  & \textcolor{Cerulean}{LA} & \textcolor{Cerulean}{2} & \textcolor{Cerulean}{0.0002} & \textcolor{Cerulean}{0.0009} & \textcolor{Cerulean}{0.0011} & \textcolor{Cerulean}{60} & \textcolor{Cerulean}{0.7271} & \textcolor{Cerulean}{1.618} & \textcolor{Cerulean}{2.345} \\ 
\hline 
\multirow{2}{*}{100} & \multirow{2}{*}{2000} & \multirow{2}{*}{15530} & \multirow{2}{*}{4.04e-01} & \multirow{2}{*}{9.60e-02} & \textcolor{Red}{SA} & \textcolor{Red}{2} & \textcolor{Red}{0.0004} & \textcolor{Red}{0.0011} & \textcolor{Red}{0.0015} & \textcolor{Red}{28} & \textcolor{Red}{0.4623} & \textcolor{Red}{0.8194} & \textcolor{Red}{1.282} \\ 
  &  &  &  &  & \textcolor{Cerulean}{LA} & \textcolor{Cerulean}{2} & \textcolor{Cerulean}{0.0003} & \textcolor{Cerulean}{0.0011} & \textcolor{Cerulean}{0.0013} & \textcolor{Cerulean}{15} & \textcolor{Cerulean}{3.793} & \textcolor{Cerulean}{1.246} & \textcolor{Cerulean}{5.039} \\ 
\hline 
\multirow{2}{*}{500} & \multirow{2}{*}{2500} & \multirow{2}{*}{3720} & \multirow{2}{*}{2.00e-02} & \multirow{2}{*}{1.80e-04} & \textcolor{Red}{SA} & \textcolor{Red}{2} & \textcolor{Red}{0.0007} & \textcolor{Red}{0.0142} & \textcolor{Red}{0.0149} & \textcolor{Red}{15} & \textcolor{Red}{0.0360} & \textcolor{Red}{0.0290} & \textcolor{Red}{0.0650} \\ 
  &  &  &  &  & \textcolor{Cerulean}{LA} & \textcolor{Cerulean}{2} & \textcolor{Cerulean}{0.0001} & \textcolor{Cerulean}{0.0135} & \textcolor{Cerulean}{0.0136} & \textcolor{Cerulean}{21} & \textcolor{Cerulean}{0.0432} & \textcolor{Cerulean}{0.0497} & \textcolor{Cerulean}{0.0929} \\ 
\hline 
\multirow{2}{*}{500} & \multirow{2}{*}{5000} & \multirow{2}{*}{16948} & \multirow{2}{*}{4.01e-02} & \multirow{2}{*}{8.18e-04} & \textcolor{Red}{SA} & \textcolor{Red}{2} & \textcolor{Red}{0.0001} & \textcolor{Red}{0.0138} & \textcolor{Red}{0.0139} & \textcolor{Red}{14} & \textcolor{Red}{0.2127} & \textcolor{Red}{0.2283} & \textcolor{Red}{0.4410} \\ 
  &  &  &  &  & \textcolor{Cerulean}{LA} & \textcolor{Cerulean}{2} & \textcolor{Cerulean}{0.0002} & \textcolor{Cerulean}{0.0138} & \textcolor{Cerulean}{0.0140} & \textcolor{Cerulean}{12} & \textcolor{Cerulean}{0.5014} & \textcolor{Cerulean}{0.2955} & \textcolor{Cerulean}{0.7969} \\ 
\hline 
\multirow{2}{*}{500} & \multirow{2}{*}{12500} & \multirow{2}{*}{110507} & \multirow{2}{*}{1.00e-01} & \multirow{2}{*}{5.34e-03} & \textcolor{Red}{SA} & \textcolor{Red}{2} & \textcolor{Red}{0.0003} & \textcolor{Red}{0.0140} & \textcolor{Red}{0.0143} & \textcolor{Red}{101} & \textcolor{Red}{4.96} & \textcolor{Red}{27.79} & \textcolor{Red}{32.75} \\ 
  &  &  &  &  & \textcolor{Cerulean}{LA} & \textcolor{Cerulean}{2} & \textcolor{Cerulean}{0.0004} & \textcolor{Cerulean}{0.0141} & \textcolor{Cerulean}{0.0145} & \textcolor{Cerulean}{86} & \textcolor{Cerulean}{86.51} & \textcolor{Cerulean}{86.57} & \textcolor{Cerulean}{173.1} \\ 
\hline 
\multirow{2}{*}{1000} & \multirow{2}{*}{5000} & \multirow{2}{*}{7386} & \multirow{2}{*}{1.00e-02} & \multirow{2}{*}{4.44e-05} & \textcolor{Red}{SA} & \textcolor{Red}{14} & \textcolor{Red}{0.0266} & \textcolor{Red}{0.0385} & \textcolor{Red}{0.0651} & \textcolor{Red}{20} & \textcolor{Red}{0.0577} & \textcolor{Red}{0.0748} & \textcolor{Red}{0.1325} \\ 
  &  &  &  &  & \textcolor{Cerulean}{LA} & \textcolor{Cerulean}{17} & \textcolor{Cerulean}{0.0250} & \textcolor{Cerulean}{0.0396} & \textcolor{Cerulean}{0.0646} & \textcolor{Cerulean}{19} & \textcolor{Cerulean}{0.0884} & \textcolor{Cerulean}{0.0831} & \textcolor{Cerulean}{0.1715} \\ 
\hline 
\multirow{2}{*}{1000} & \multirow{2}{*}{10000} & \multirow{2}{*}{33022} & \multirow{2}{*}{2.00e-02} & \multirow{2}{*}{1.99e-04} & \textcolor{Red}{SA} & \textcolor{Red}{17} & \textcolor{Red}{0.0264} & \textcolor{Red}{0.0433} & \textcolor{Red}{0.0697} & \textcolor{Red}{17} & \textcolor{Red}{0.4924} & \textcolor{Red}{0.5655} & \textcolor{Red}{1.058} \\ 
  &  &  &  &  & \textcolor{Cerulean}{LA} & \textcolor{Cerulean}{10} & \textcolor{Cerulean}{0.0307} & \textcolor{Cerulean}{0.0407} & \textcolor{Cerulean}{0.0715} & \textcolor{Cerulean}{13} & \textcolor{Cerulean}{1.152} & \textcolor{Cerulean}{0.6456} & \textcolor{Cerulean}{1.798} \\ 
\hline 
\multirow{2}{*}{1000} & \multirow{2}{*}{25000} & \multirow{2}{*}{220002} & \multirow{2}{*}{5.01e-02} & \multirow{2}{*}{1.32e-03} & \textcolor{Red}{SA} & \textcolor{Red}{15} & \textcolor{Red}{0.0308} & \textcolor{Red}{0.0515} & \textcolor{Red}{0.0823} & \textcolor{Red}{29} & \textcolor{Red}{10.21} & \textcolor{Red}{16.5} & \textcolor{Red}{26.71} \\ 
  &  &  &  &  & \textcolor{Cerulean}{LA} & \textcolor{Cerulean}{6} & \textcolor{Cerulean}{0.0474} & \textcolor{Cerulean}{0.0559} & \textcolor{Cerulean}{0.1034} & \textcolor{Cerulean}{23} & \textcolor{Cerulean}{192.5} & \textcolor{Cerulean}{48.85} & \textcolor{Cerulean}{241.3} \\ 
\hline 
\end{tabular} 

 \caption{\WS~graphs.}
 \label{tab:wsamgsttstcs} 
\end{table}

\begin{table}[h] \scriptsize
  \centering
  \begin{tabular}{|ccccc|c|cccc|cccc|}
\hline
\multirow{2}{*}{$N_0$} & \multirow{2}{*}{$N_1$} & \multirow{2}{*}{$N_2$} & Edge & Triangle & Algorithm / & \multicolumn{4}{|c|}{$\alpha$} & \multicolumn{4}{|c|}{$\beta$} \\
  &  &  & Density & Density & Formulation &  iter. & setup & solve & total &  iter. & setup & solve & total \\ 
 \hline 
\multirow{2}{*}{100} & \multirow{2}{*}{475} & \multirow{2}{*}{301} & \multirow{2}{*}{9.60e-02} & \multirow{2}{*}{1.86e-03} & \textcolor{Red}{SA} & \textcolor{Red}{2} & \textcolor{Red}{0.0001} & \textcolor{Red}{0.0008} & \textcolor{Red}{0.0009} & \textcolor{Red}{2} & \textcolor{Red}{0.0001} & \textcolor{Red}{0.0037} & \textcolor{Red}{0.0038} \\ 
  &  &  &  &  & \textcolor{Cerulean}{LA} & \textcolor{Cerulean}{2} & \textcolor{Cerulean}{0.0002} & \textcolor{Cerulean}{0.0009} & \textcolor{Cerulean}{0.0011} & \textcolor{Cerulean}{2} & \textcolor{Cerulean}{0.0001} & \textcolor{Cerulean}{0.0037} & \textcolor{Cerulean}{0.0037} \\ 
\hline 
\multirow{2}{*}{100} & \multirow{2}{*}{900} & \multirow{2}{*}{1701} & \multirow{2}{*}{1.82e-01} & \multirow{2}{*}{1.05e-02} & \textcolor{Red}{SA} & \textcolor{Red}{2} & \textcolor{Red}{0.0005} & \textcolor{Red}{0.0012} & \textcolor{Red}{0.0016} & \textcolor{Red}{50} & \textcolor{Red}{0.0242} & \textcolor{Red}{0.0753} & \textcolor{Red}{0.0995} \\ 
  &  &  &  &  & \textcolor{Cerulean}{LA} & \textcolor{Cerulean}{2} & \textcolor{Cerulean}{0.0001} & \textcolor{Cerulean}{0.0008} & \textcolor{Cerulean}{0.0009} & \textcolor{Cerulean}{44} & \textcolor{Cerulean}{0.0358} & \textcolor{Cerulean}{0.0692} & \textcolor{Cerulean}{0.1050} \\ 
\hline 
\multirow{2}{*}{100} & \multirow{2}{*}{1600} & \multirow{2}{*}{8105} & \multirow{2}{*}{3.23e-01} & \multirow{2}{*}{5.01e-02} & \textcolor{Red}{SA} & \textcolor{Red}{2} & \textcolor{Red}{0.0001} & \textcolor{Red}{0.0008} & \textcolor{Red}{0.0009} & \textcolor{Red}{14} & \textcolor{Red}{0.1842} & \textcolor{Red}{0.1720} & \textcolor{Red}{0.3561} \\ 
  &  &  &  &  & \textcolor{Cerulean}{LA} & \textcolor{Cerulean}{2} & \textcolor{Cerulean}{0.0001} & \textcolor{Cerulean}{0.0009} & \textcolor{Cerulean}{0.0010} & \textcolor{Cerulean}{12} & \textcolor{Cerulean}{1.149} & \textcolor{Cerulean}{0.3968} & \textcolor{Cerulean}{1.545} \\ 
\hline 
\multirow{2}{*}{100} & \multirow{2}{*}{2400} & \multirow{2}{*}{24497} & \multirow{2}{*}{4.85e-01} & \multirow{2}{*}{1.51e-01} & \textcolor{Red}{SA} & \textcolor{Red}{2} & \textcolor{Red}{0.0005} & \textcolor{Red}{0.0012} & \textcolor{Red}{0.0017} & \textcolor{Red}{24} & \textcolor{Red}{1.111} & \textcolor{Red}{1.521} & \textcolor{Red}{2.632} \\ 
  &  &  &  &  & \textcolor{Cerulean}{LA} & \textcolor{Cerulean}{2} & \textcolor{Cerulean}{0.0002} & \textcolor{Cerulean}{0.0011} & \textcolor{Cerulean}{0.0013} & \textcolor{Cerulean}{13} & \textcolor{Cerulean}{19.39} & \textcolor{Cerulean}{3.023} & \textcolor{Cerulean}{22.41} \\ 
\hline 
\multirow{2}{*}{500} & \multirow{2}{*}{4900} & \multirow{2}{*}{4740} & \multirow{2}{*}{3.93e-02} & \multirow{2}{*}{2.29e-04} & \textcolor{Red}{SA} & \textcolor{Red}{2} & \textcolor{Red}{0.0016} & \textcolor{Red}{0.0153} & \textcolor{Red}{0.0169} & \textcolor{Red}{89} & \textcolor{Red}{0.0742} & \textcolor{Red}{0.3261} & \textcolor{Red}{0.4003} \\ 
  &  &  &  &  & \textcolor{Cerulean}{LA} & \textcolor{Cerulean}{2} & \textcolor{Cerulean}{0.0001} & \textcolor{Cerulean}{0.0136} & \textcolor{Cerulean}{0.0137} & \textcolor{Cerulean}{95} & \textcolor{Cerulean}{0.1027} & \textcolor{Cerulean}{0.4019} & \textcolor{Cerulean}{0.5046} \\ 
\hline 
\multirow{2}{*}{500} & \multirow{2}{*}{9600} & \multirow{2}{*}{25016} & \multirow{2}{*}{7.70e-02} & \multirow{2}{*}{1.21e-03} & \textcolor{Red}{SA} & \textcolor{Red}{2} & \textcolor{Red}{0.0002} & \textcolor{Red}{0.0140} & \textcolor{Red}{0.0142} & \textcolor{Red}{101} & \textcolor{Red}{1.018} & \textcolor{Red}{3.608} & \textcolor{Red}{4.627} \\ 
  &  &  &  &  & \textcolor{Cerulean}{LA} & \textcolor{Cerulean}{2} & \textcolor{Cerulean}{0.0004} & \textcolor{Cerulean}{0.0138} & \textcolor{Cerulean}{0.0142} & \textcolor{Cerulean}{101} & \textcolor{Cerulean}{10.56} & \textcolor{Cerulean}{19.19} & \textcolor{Cerulean}{29.74} \\ 
\hline 
\multirow{2}{*}{1000} & \multirow{2}{*}{9900} & \multirow{2}{*}{6264} & \multirow{2}{*}{1.98e-02} & \multirow{2}{*}{3.77e-05} & \textcolor{Red}{SA} & \textcolor{Red}{5} & \textcolor{Red}{0.0215} & \textcolor{Red}{0.0261} & \textcolor{Red}{0.0477} & \textcolor{Red}{101} & \textcolor{Red}{0.0946} & \textcolor{Red}{0.4777} & \textcolor{Red}{0.5723} \\ 
  &  &  &  &  & \textcolor{Cerulean}{LA} & \textcolor{Cerulean}{5} & \textcolor{Cerulean}{0.0316} & \textcolor{Cerulean}{0.0366} & \textcolor{Cerulean}{0.0682} & \textcolor{Cerulean}{101} & \textcolor{Cerulean}{0.2080} & \textcolor{Cerulean}{1.012} & \textcolor{Cerulean}{1.22} \\ 
\hline 
\multirow{2}{*}{1000} & \multirow{2}{*}{19600} & \multirow{2}{*}{37365} & \multirow{2}{*}{3.92e-02} & \multirow{2}{*}{2.25e-04} & \textcolor{Red}{SA} & \textcolor{Red}{4} & \textcolor{Red}{0.0174} & \textcolor{Red}{0.0219} & \textcolor{Red}{0.0393} & \textcolor{Red}{101} & \textcolor{Red}{2.403} & \textcolor{Red}{8.212} & \textcolor{Red}{10.61} \\ 
  &  &  &  &  & \textcolor{Cerulean}{LA} & \textcolor{Cerulean}{4} & \textcolor{Cerulean}{0.0387} & \textcolor{Cerulean}{0.0436} & \textcolor{Cerulean}{0.0823} & \textcolor{Cerulean}{101} & \textcolor{Cerulean}{27.14} & \textcolor{Cerulean}{35.86} & \textcolor{Cerulean}{63.01} \\ 
\hline 
\end{tabular} 

  \caption{\BA~graphs.}
  \label{tab:baamgsttstcs} 
\end{table}

\clearpage

\section{Nonzeros in Laplacian and boundary matrices}
\label{appndx:nonzeros}

The next three tables show the number of nonzeros in the boundary and
Laplacian matrices. In addition we also show the number of nonzeros in
the coarse grid matrices that arise in algebraic multigrid. The latter
can be useful in interpreting and debugging algebraic multigrid
performance. The cases where PyAMG failed in the setup phase are
marked by a dagger symbol~($\dagger$).

\bigskip
\begin{table}[h] \scriptsize
  \centering
  \begin{tabular}{|ccccc|cccc|cccc|} 
\hline 
\multirow{3}{*}{$N_0$} & \multirow{3}{*}{$N_1$} & \multirow{3}{*}{$N_2$} & \multirow{2}{*}{Edge} & \multirow{2}{*}{Triangle} & \multicolumn{4}{|c|}{Krylov} & \multicolumn{4}{|c|}{AMG} \\
  &  &  & \multirow{2}{*}{Density} & \multirow{2}{*}{Density} & \multirow{2}{*}{$\boundary_1$} & \multirow{2}{*}{$\boundary_2$} & \multirow{2}{*}{$\laplacian_0$} & \multirow{2}{*}{$\laplacian_2$} & \multicolumn{2}{|c}{SA} & \multicolumn{2}{c|}{LA} \\
  &  &  &  &  &  &  &  &  & $\laplacian_0$ & $\laplacian_2$ & $\laplacian_0$ & $\laplacian_2$ \\ 
 \hline 
100 & 380 & 52 & 7.68e-02 & 3.22e-04 & 760 & 156 & 860 & 110 & 860 & 110 & 860 & 110 \\
\hline 
100 & 494 & 144 & 9.98e-02 & 8.91e-04 & 988 & 432 & 1088 & 528 & 1088 & 528 & 1088 & 528 \\
\hline 
\multirow{2}{*}{100} & \multirow{2}{*}{1212} & \multirow{2}{*}{2359} & \multirow{2}{*}{2.45e-01} & \multirow{2}{*}{1.46e-02} & \multirow{2}{*}{2424} & \multirow{2}{*}{7077} & \multirow{2}{*}{2524} & \multirow{2}{*}{42475} & 2524 & 42475 & 2524 & 42475 \\
  &  &  &  &  &  &  &  &  &  & 3721 &  & 49902 \\
\hline 
\multirow{3}{*}{100} & \multirow{3}{*}{2530} & \multirow{3}{*}{21494} & \multirow{3}{*}{5.11e-01} & \multirow{3}{*}{1.33e-01} & \multirow{3}{*}{5060} & \multirow{3}{*}{64482} & \multirow{3}{*}{5160} & \multirow{3}{*}{1645012} & 5160 & 1645012 & 5160 & 1645012 \\
  &  &  &  &  &  &  &  &  &  & 2500 &  & 4174983 \\
  &  &  &  &  &  &  &  &  &  &  &  & 45796 \\
\hline 
\multirow{4}{*}{100} & \multirow{4}{*}{3706} & \multirow{4}{*}{67865} & \multirow{4}{*}{7.49e-01} & \multirow{4}{*}{4.20e-01} & \multirow{4}{*}{7412} & \multirow{4}{*}{203595} & \multirow{4}{*}{7512} & \multirow{4}{*}{11134203} & 7512 & 11134203 & 7512 & 11134203 \\
  &  &  &  &  &  &  &  &  &  & 1681 &  & 34985035 \\
  &  &  &  &  &  &  &  &  &  &  &  & 459684 \\
  &  &  &  &  &  &  &  &  &  &  &  & 4489 \\
\hline 
500 & 1290 & 21 & 1.03e-02 & 1.01e-06 & 2580 & 63 & 3075 & 25 & 3075 & 25 & 3075 & 25 \\
\hline 
\multirow{3}{*}{500} & \multirow{3}{*}{12394} & \multirow{3}{*}{20315} & \multirow{3}{*}{9.94e-02} & \multirow{3}{*}{9.81e-04} & \multirow{3}{*}{24788} & \multirow{3}{*}{60945} & \multirow{3}{*}{25288} & \multirow{3}{*}{319503} & 25288 & 319503 & 25288 & 319503 \\
  &  &  &  &  &  &  &  &  &  & 454789 &  & 1739567 \\
  &  &  &  &  &  &  &  &  &  & 1 &  & 41209 \\
\hline 
500 & 24788 & 162986 & 1.99e-01 & 7.87e-03 & 49576 & 488958 & 50076 &
9807176 & $\dagger$ & $\dagger$ & $\dagger$ & $\dagger$ \\
\hline 
1000 & 49690 & 163767 & 9.95e-02 & 9.86e-04 & 99380 & 491301 & 100380 & 5033205 & $\dagger$ & $\dagger$ & $\dagger$ & $\dagger$ \\
\hline 
\end{tabular} 

  \caption{\ER~graphs.}
  \label{tab:ernnz} 
\end{table}

\begin{table}[h] \scriptsize
  \centering
  \begin{tabular}{|ccccc|cccc|cccc|} 
\hline 
\multirow{3}{*}{$N_0$} & \multirow{3}{*}{$N_1$} & \multirow{3}{*}{$N_2$} & \multirow{2}{*}{Edge} & \multirow{2}{*}{Triangle} & \multicolumn{4}{|c|}{Krylov} & \multicolumn{4}{|c|}{AMG} \\
  &  &  & \multirow{2}{*}{Density} & \multirow{2}{*}{Density} & \multirow{2}{*}{$\boundary_1$} & \multirow{2}{*}{$\boundary_2$} & \multirow{2}{*}{$\laplacian_0$} & \multirow{2}{*}{$\laplacian_2$} & \multicolumn{2}{|c}{SA} & \multicolumn{2}{c|}{LA} \\
  &  &  &  &  &  &  &  &  & $\laplacian_0$ & $\laplacian_2$ & $\laplacian_0$ & $\laplacian_2$ \\ 
 \hline 
\multirow{2}{*}{100} & \multirow{2}{*}{500} & \multirow{2}{*}{729} & \multirow{2}{*}{1.01e-01} & \multirow{2}{*}{4.51e-03} & \multirow{2}{*}{1000} & \multirow{2}{*}{2187} & \multirow{2}{*}{1100} & \multirow{2}{*}{10053} & 1100 & 10053 & 1100 & 10053 \\
  &  &  &  &  &  &  &  &  &  & 302 &  & 1145 \\
\hline 
\multirow{2}{*}{100} & \multirow{2}{*}{1000} & \multirow{2}{*}{3655} & \multirow{2}{*}{2.02e-01} & \multirow{2}{*}{2.26e-02} & \multirow{2}{*}{2000} & \multirow{2}{*}{10965} & \multirow{2}{*}{2100} & \multirow{2}{*}{127001} & 2100 & 127001 & 2100 & 127001 \\
  &  &  &  &  &  &  &  &  &  & 929 &  & 50045 \\
\hline 
\multirow{3}{*}{100} & \multirow{3}{*}{1500} & \multirow{3}{*}{8354} & \multirow{3}{*}{3.03e-01} & \multirow{3}{*}{5.17e-02} & \multirow{3}{*}{3000} & \multirow{3}{*}{25062} & \multirow{3}{*}{3100} & \multirow{3}{*}{443774} & 3100 & 443774 & 3100 & 443774 \\
  &  &  &  &  &  &  &  &  &  & 1921 &  & 373960 \\
  &  &  &  &  &  &  &  &  &  &  &  & 6889 \\
\hline 
\multirow{3}{*}{100} & \multirow{3}{*}{2000} & \multirow{3}{*}{15530} & \multirow{3}{*}{4.04e-01} & \multirow{3}{*}{9.60e-02} & \multirow{3}{*}{4000} & \multirow{3}{*}{46590} & \multirow{3}{*}{4100} & \multirow{3}{*}{1138022} & 4100 & 1138022 & 4100 & 1138022 \\
  &  &  &  &  &  &  &  &  &  & 2116 &  & 1622123 \\
  &  &  &  &  &  &  &  &  &  &  &  & 24025 \\
\hline 
\multirow{2}{*}{500} & \multirow{2}{*}{2500} & \multirow{2}{*}{3720} & \multirow{2}{*}{2.00e-02} & \multirow{2}{*}{1.80e-04} & \multirow{2}{*}{5000} & \multirow{2}{*}{11160} & \multirow{2}{*}{5500} & \multirow{2}{*}{52456} & 5500 & 52456 & 5500 & 52456 \\
  &  &  &  &  &  &  &  &  &  & 5500 &  & 6210 \\
\hline 
\multirow{3}{*}{500} & \multirow{3}{*}{5000} & \multirow{3}{*}{16948} & \multirow{3}{*}{4.01e-02} & \multirow{3}{*}{8.18e-04} & \multirow{3}{*}{10000} & \multirow{3}{*}{50844} & \multirow{3}{*}{10500} & \multirow{3}{*}{568252}  & 5500 & 568252 & 10500 & 568252 \\
  &  &  &  &  &  &  &  &  &  & 3808 &  & 205516 \\
  &  &  &  &  &  &  &  &  &  &  &  & 10015 \\
\hline 
\multirow{4}{*}{500} & \multirow{4}{*}{12500} & \multirow{4}{*}{110507} & \multirow{4}{*}{1.00e-01} & \multirow{4}{*}{5.34e-03} & \multirow{4}{*}{25000} & \multirow{4}{*}{331521} & \multirow{4}{*}{25500} & \multirow{4}{*}{9870149} & 25500 & 9870149 & 25500 & 9870149 \\
  &  &  &  &  &  &  &  &  &  & 38306 &  & 16959388 \\
  &  &  &  &  &  &  &  &  &  &  &  & 1221025 \\
  &  &  &  &  &  &  &  &  &  &  &  & 12100 \\
\hline 
\multirow{3}{*}{1000} & \multirow{3}{*}{5000} & \multirow{3}{*}{7386} & \multirow{3}{*}{1.00e-02} & \multirow{3}{*}{4.44e-05} & \multirow{3}{*}{10000} & \multirow{3}{*}{22158} & \multirow{3}{*}{11000} & \multirow{3}{*}{103946} & 11000 & 103946 & 11000 & 103946 \\
  &  &  &  &  &  &  &  &  & 5231 & 2787 & 8110 & 12339 \\
  &  &  &  &  &  &  &  &  &  &  &  & 631 \\
\hline 
\multirow{3}{*}{1000} & \multirow{3}{*}{10000} & \multirow{3}{*}{33022} & \multirow{3}{*}{2.00e-02} & \multirow{3}{*}{1.99e-04} & \multirow{3}{*}{20000} & \multirow{3}{*}{99066} & \multirow{3}{*}{21000} & \multirow{3}{*}{1086816} & 21000 & 1086816 & 21000 & 1086816 \\
  &  &  &  &  &  &  &  &  & 900 & 7442 & 10000 & 393322 \\
  &  &  &  &  &  &  &  &  &  &  &  & 22318 \\
\hline 
\multirow{4}{*}{1000} & \multirow{4}{*}{25000} & \multirow{4}{*}{220002} & \multirow{4}{*}{5.01e-02} & \multirow{4}{*}{1.32e-03} & \multirow{4}{*}{50000} & \multirow{4}{*}{660006} & \multirow{4}{*}{51000} & \multirow{4}{*}{19688290} & 51000 & 19688290 & 51000 & 19688290 \\
  &  &  &  &  &  &  &  &  & 121 & 42253 & 10000 & 32579049 \\
  &  &  &  &  &  &  &  &  &  &  &  & 4819514 \\
  &  &  &  &  &  &  &  &  &  &  &  & 48400 \\
\hline 
\end{tabular} 

  \caption{\WS~graphs.}
  \label{tab:wsnnz} 
\end{table}

\begin{table}[h] \scriptsize
  \centering
  \begin{tabular}{|ccccc|cccc|cccc|} 
\hline 
\multirow{3}{*}{$N_0$} & \multirow{3}{*}{$N_1$} & \multirow{3}{*}{$N_2$} & \multirow{2}{*}{Edge} & \multirow{2}{*}{Triangle} & \multicolumn{4}{|c|}{Krylov} & \multicolumn{4}{|c|}{AMG} \\
  &  &  & \multirow{2}{*}{Density} & \multirow{2}{*}{Density} & \multirow{2}{*}{$\boundary_1$} & \multirow{2}{*}{$\boundary_2$} & \multirow{2}{*}{$\laplacian_0$} & \multirow{2}{*}{$\laplacian_2$} & \multicolumn{2}{|c}{SA} & \multicolumn{2}{c|}{LA} \\
  &  &  &  &  &  &  &  &  & $\laplacian_0$ & $\laplacian_2$ & $\laplacian_0$ & $\laplacian_2$ \\ 
\hline 
100 & 475 & 301 & 9.60e-02 & 1.86e-03 & 950 & 903 & 1050 & 3189 & 1050 & 3189 & 1050 & 3189 \\
\hline 
\multirow{2}{*}{100} & \multirow{2}{*}{900} & \multirow{2}{*}{1701} & \multirow{2}{*}{1.82e-01} & \multirow{2}{*}{1.05e-02} & \multirow{2}{*}{1800} & \multirow{2}{*}{5103} & \multirow{2}{*}{1900} & \multirow{2}{*}{42465} & 1900 & 42465 & 1900 & 42465 \\
  &  &  &  &  &  &  &  &  &  & 1837 &  & 26324 \\
\hline 
\multirow{3}{*}{100} & \multirow{3}{*}{1600} & \multirow{3}{*}{8105} & \multirow{3}{*}{3.23e-01} & \multirow{3}{*}{5.01e-02} & \multirow{3}{*}{3200} & \multirow{3}{*}{24315} & \multirow{3}{*}{3300} & \multirow{3}{*}{472503} & 3300 & 472503 & 3300 & 472503 \\
  &  &  &  &  &  &  &  &  &  & 2401 &  & 607290 \\
  &  &  &  &  &  &  &  &  &  &  &  & 6561 \\
\hline 
\multirow{3}{*}{100} & \multirow{3}{*}{2400} & \multirow{3}{*}{24497} & \multirow{3}{*}{4.85e-01} & \multirow{3}{*}{1.51e-01} & \multirow{3}{*}{4800} & \multirow{3}{*}{73491} & \multirow{3}{*}{4900} & \multirow{3}{*}{2610185} & 4900 & 2610185 & 4900 & 2610185 \\
  &  &  &  &  &  &  &  &  &  & 1936 &  & 5244756 \\
  &  &  &  &  &  &  &  &  &  &  &  & 59536 \\
\hline 
\multirow{2}{*}{500} & \multirow{2}{*}{4900} & \multirow{2}{*}{4740} & \multirow{2}{*}{3.93e-02} & \multirow{2}{*}{2.29e-04} & \multirow{2}{*}{9800} & \multirow{2}{*}{14220} & \multirow{2}{*}{10300} & \multirow{2}{*}{121426} & 10300 & 121426 & 10300 & 121426 \\
  &  &  &  &  &  &  &  &  &  & 44271 &  & 138744 \\
\hline 
\multirow{3}{*}{500} & \multirow{3}{*}{9600} & \multirow{3}{*}{25016} & \multirow{3}{*}{7.70e-02} & \multirow{3}{*}{1.21e-03} & \multirow{3}{*}{19200} & \multirow{3}{*}{75048} & \multirow{3}{*}{19700} & \multirow{3}{*}{1180338} & 19700 & 1180338 & 19700 & 1180338 \\
  &  &  &  &  &  &  &  &  &  & 161347 &  & 4353167 \\
  &  &  &  &  &  &  &  &  &  &  &  & 62500 \\
\hline 
500 & 18400 & 133933 & 1.47e-01 & 6.47e-03 & 36800 & 401799 & 37300 & 14303959 & $\dagger$ & $\dagger$ & $\dagger$ & $\dagger$ \\
\hline 
\multirow{3}{*}{1000} & \multirow{3}{*}{9900} & \multirow{3}{*}{6264} & \multirow{3}{*}{1.98e-02} & \multirow{3}{*}{3.77e-05} & \multirow{3}{*}{19800} & \multirow{3}{*}{18792} & \multirow{3}{*}{20800} & \multirow{3}{*}{132712} & 20800 & 132712 & 20800 & 132712 \\
  &  &  &  &  &  &  &  &  & 196 & 119206 & 10000 & 161391 \\
  &  &  &  &  &  &  &  &  &  &  &  & 3722 \\
\hline 
\multirow{3}{*}{1000} & \multirow{3}{*}{19600} & \multirow{3}{*}{37365} & \multirow{3}{*}{3.92e-02} & \multirow{3}{*}{2.25e-04} & \multirow{3}{*}{39200} & \multirow{3}{*}{112095} & \multirow{3}{*}{40200} & \multirow{3}{*}{1808827} & 40200 & 1808827 & 40200 & 1808827 \\
  &  &  &  &  &  &  &  &  & 16 & 774769 & 10000 & 8082839 \\
  &  &  &  &  &  &  &  &  &  & 1 &  & 139129 \\
\hline 
1000 & 38400 & 202731 & 7.69e-02 & 1.22e-03 & 76800 & 608193 & 77800 & 19768919 & $\dagger$ & $\dagger$ & $\dagger$ & $\dagger$ \\
\hline 
\end{tabular} 

  \caption{\BA~graphs.}
  \label{tab:bannz} 
\end{table}

\clearpage

\section{Nonzero patterns in Laplacians} \label{appndx:spy} 

\begin{figure}[p]
  \centering
  \includegraphics[scale=1, trim=1.2in 6.4in 1.1in 0.9in, clip]
  {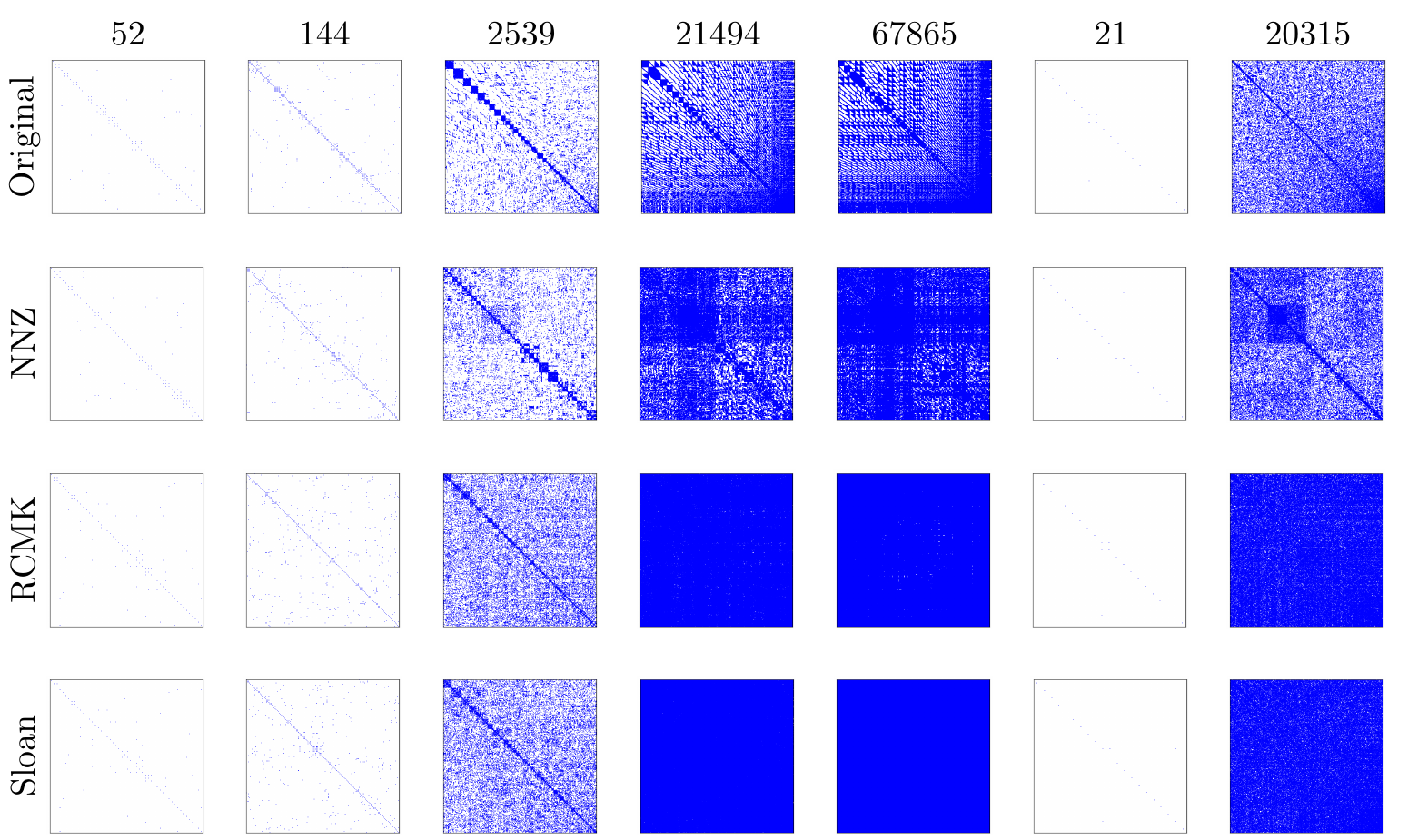}
 \caption{$\laplacian_2$ for \ER~graphs.}
  \label{fig:erlplcln2rrdrng}
\end{figure}

The next nine figures show the nonzero patterns for the Laplacian
matrices $\laplacian_0$ and $\laplacian_2$ and their various
reorderings for different types of graphs.  The four rows from top to
bottom in some of the figures show the pattern of nonzeros in the
$\laplacian_0$ or $\laplacian_2$ matrix and their reordering. The
reorderings are by vertex degree (or number of nonzeros), reverse
Cuthill-McKee algorithm, and Sloan's algorithm. The rows are labelled
as ``Original'', ``Degree'' (or ``NNZ''), ``RCMK'', and ``Sloan''. In
each column the size of matrix is shown in the top row. In each
collection all matrices are drawn with the same size.  This results in
some visual artifacts in the $\laplacian_2$ figures where the matrices
are of vastly different sizes. For example, in the middle column of
Figure~\ref{fig:erlplcln2rrdrng} the top matrix and the bottom matrix
have the same number of nonzeros even though it does not appear that
way in the figure.

\vspace{-0.2in}
\begin{figure}[h]
  \centering
  \includegraphics[scale=1, trim=1.2in 6.4in 1.1in 0.9in, clip]
  {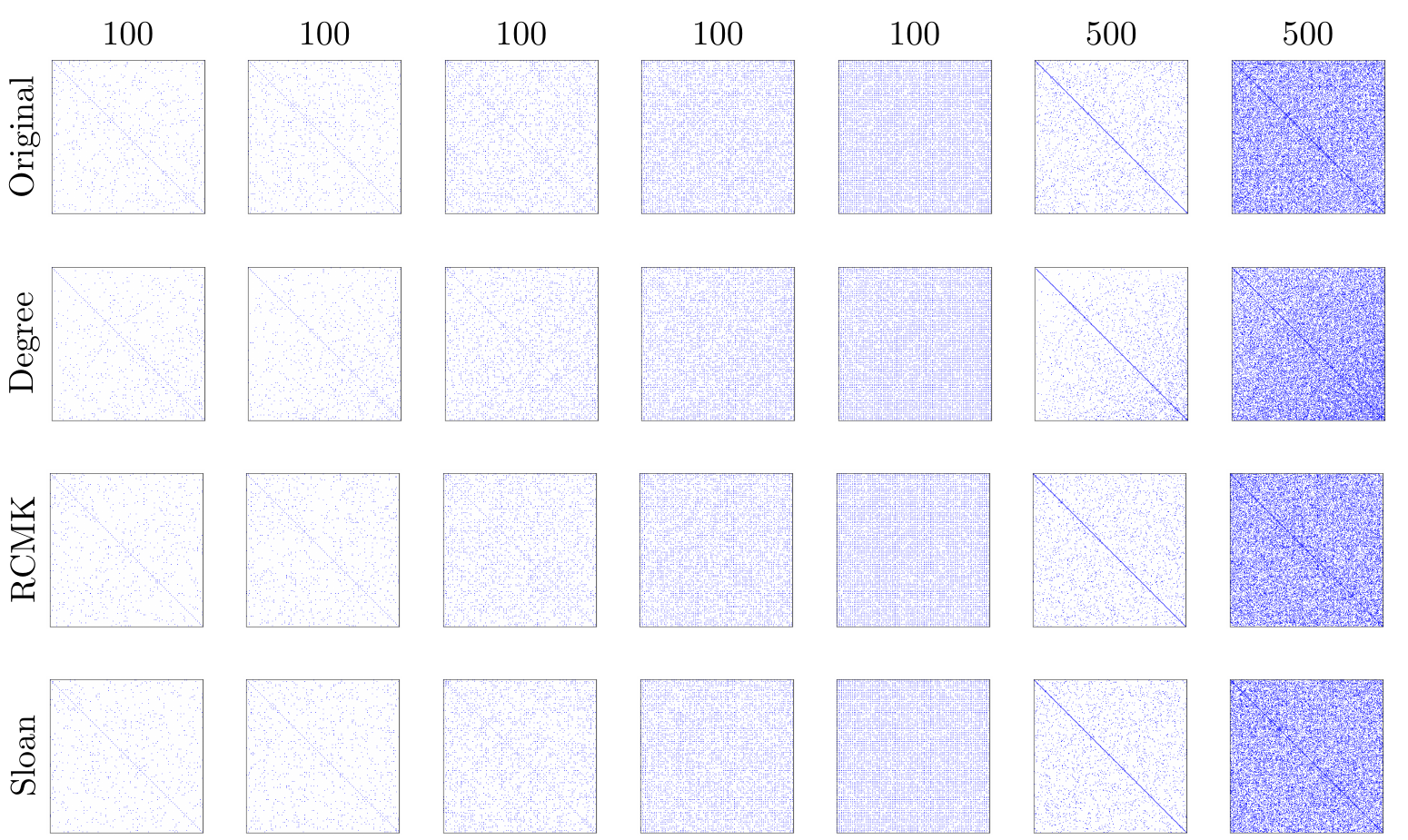}
  \caption{$\laplacian_0$ for \ER~graphs.}
\label{fig:erlplcln0rrdrng}
\end{figure}

\vspace{-0.3in}
\begin{figure}[H]
  \centering
  \begin{tabular}{ccc}
    \imagetop{\begin{xyoverpic}
        {(1, 1)}{scale=0.5}{spylplcn/er/lplcn0/n500p0.1zm.png},
        (0.5, 1.03)*{\text{Original}}
      \end{xyoverpic}} &
    \imagetop{\begin{xyoverpic}
        {(1, 1)}{scale=0.5}{spylplcn/er/lplcn0/n500p0.1dgrzm.png},
        (0.5, 1.03)*{\text{Vertex degree}}
      \end{xyoverpic}} &
    \imagetop{\begin{xyoverpic}
        {(1, 1)}{scale=0.5}{spylplcn/er/lplcn0/n500p0.1clrbr.png},
        (0.5, 1.05)*{\text{~}},
        (1.9, 0.98)*{10^2},
        (1.9, 0.685)*{10^1},
        (1.9, 0.37)*{10^0},
        (2.1, 0.055)*{10^{-1}}
      \end{xyoverpic}}
  \end{tabular}
  \caption{Enlarged view of the matrix corresponding to the
    last column in Figure~\ref{fig:erlplcln0rrdrng}.}
  \label{fig:erlplcln0zm}
\end{figure}

\clearpage

\begin{figure}[p]
  \centering
  \begin{tabular}{ccc}
    \imagetop{\begin{xyoverpic}
        {(1, 1)}{scale=0.5}{spylplcn/er/lplcn2/n500p0.1zm.png},
        (0.5, 1.03)*{\text{Original}}
      \end{xyoverpic}} &
    \imagetop{\begin{xyoverpic}
        {(1, 1)}{scale=0.5}{spylplcn/er/lplcn2/n500p0.1nnzzm.png},
        (0.5, 1.03)*{\text{Number of nonzeros}}
      \end{xyoverpic}} &
    \imagetop{\begin{xyoverpic}
        {(1, 1)}{scale=0.5}{spylplcn/er/lplcn2/n500p0.1clrbr.png},
        (0.5, 1.05)*{\text{~}},
        (1.9, 0.98)*{10^2},
        (1.9, 0.685)*{10^1},
        (1.9, 0.37)*{10^0},
        (2.1, 0.055)*{10^{-1}}
      \end{xyoverpic}}
  \end{tabular}
\caption{Enlarged view of the matrix corresponding to the
    last column in Figure~\ref{fig:erlplcln2rrdrng}.}
\label{fig:erlplcln2zm}
\end{figure}

\clearpage

\begin{figure}[p]
  \centering
  \begin{tabular}{c}
  \includegraphics[scale=0.94, trim=1.2in 6.5in 1.1in 0.98in, clip]
  {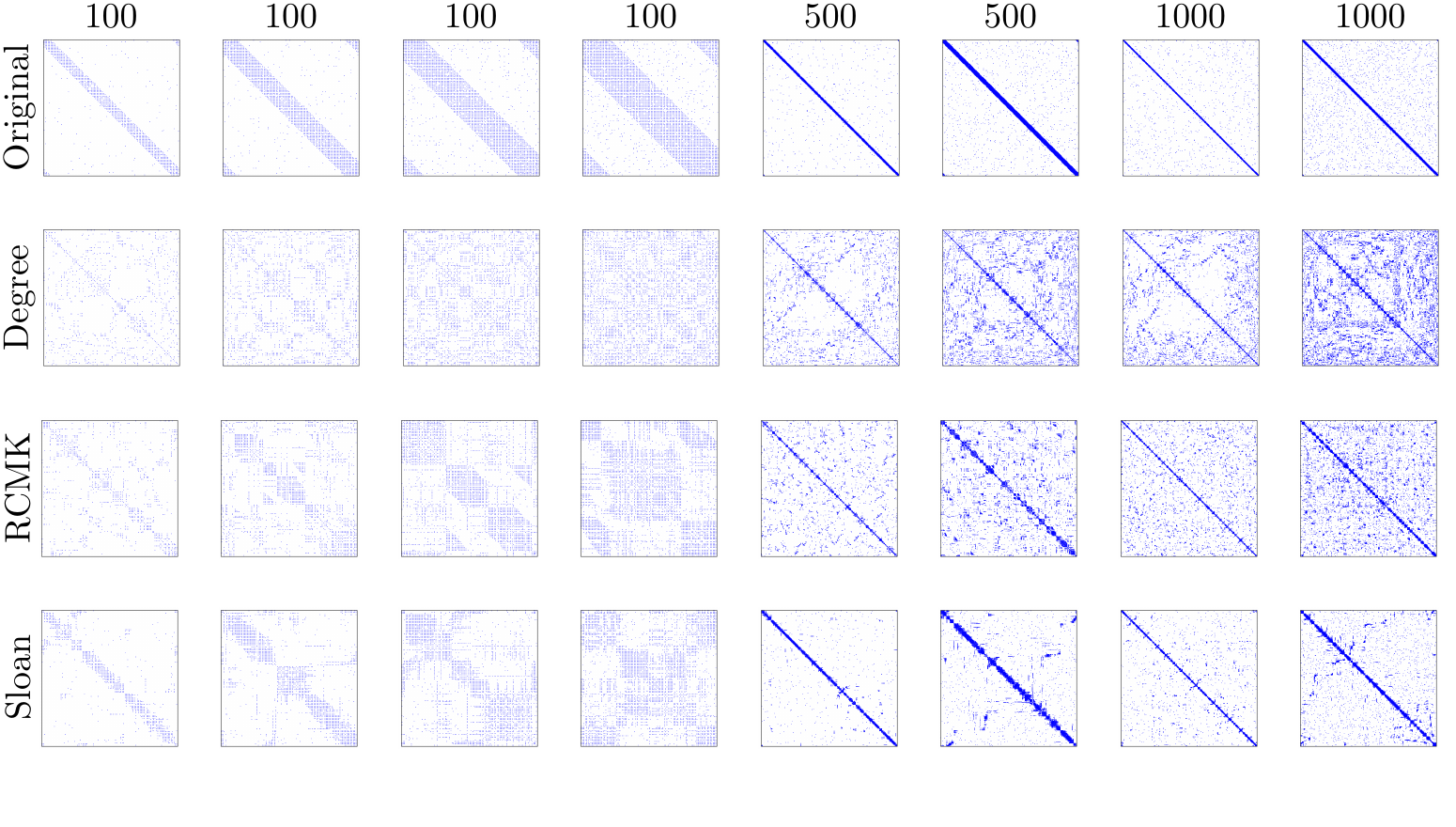} \\
  \hline \\
  \includegraphics[scale=0.94, trim=1.2in 6.5in 1.1in 0.98in, clip]
  {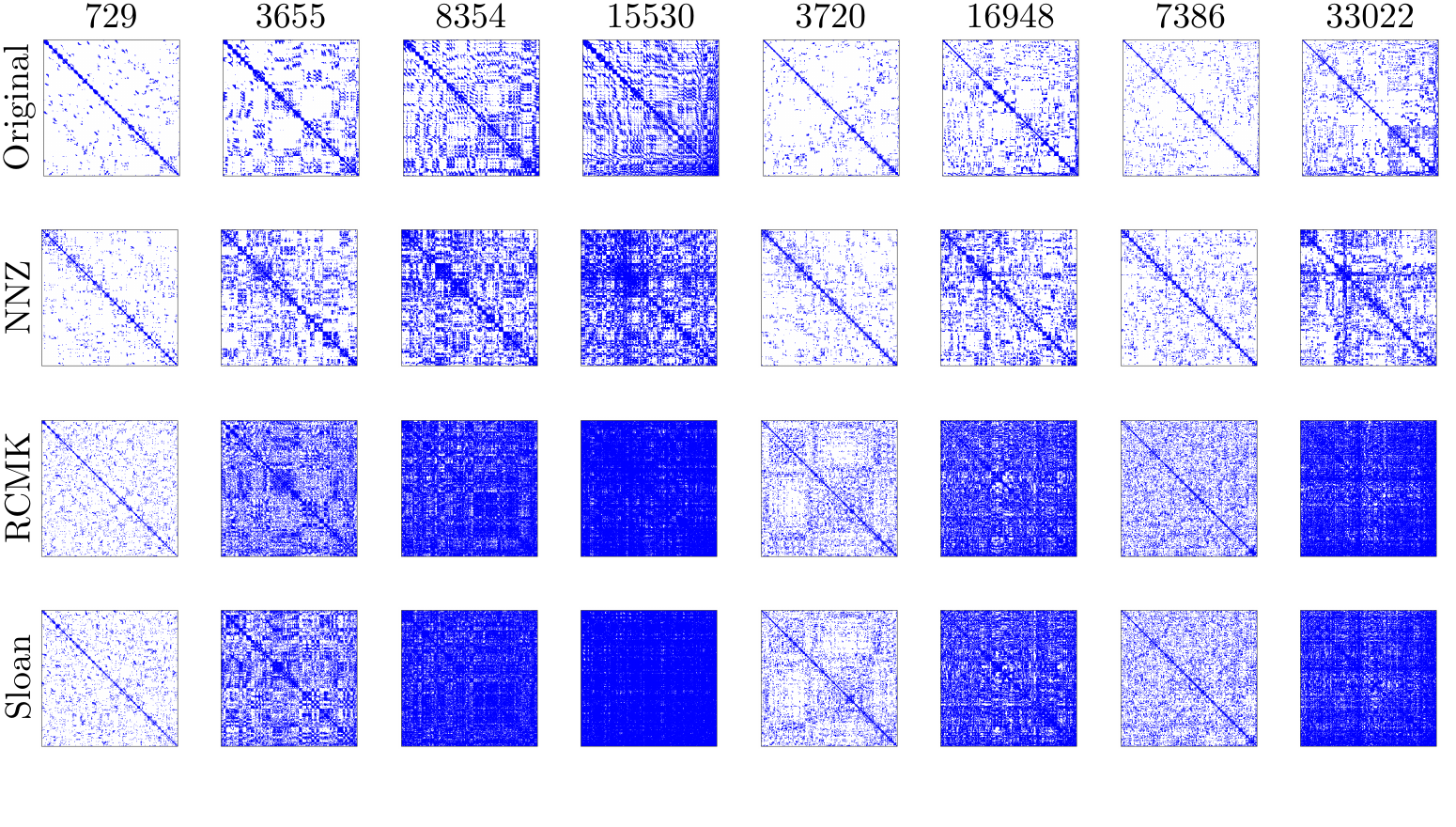}
  \end{tabular}
  \caption{\textit{Top four rows :} $\laplacian_0$ for
    \WS~graphs. \textit{Bottom four rows :} $\laplacian_2$ for
    \WS~graphs.}
  \label{fig:wslplcln2rrdrng}
\end{figure}

\clearpage

\begin{figure}[p]
  \centering
  \includegraphics[scale=1, trim=1.2in 6.5in 1.1in 0.9in, clip]
  {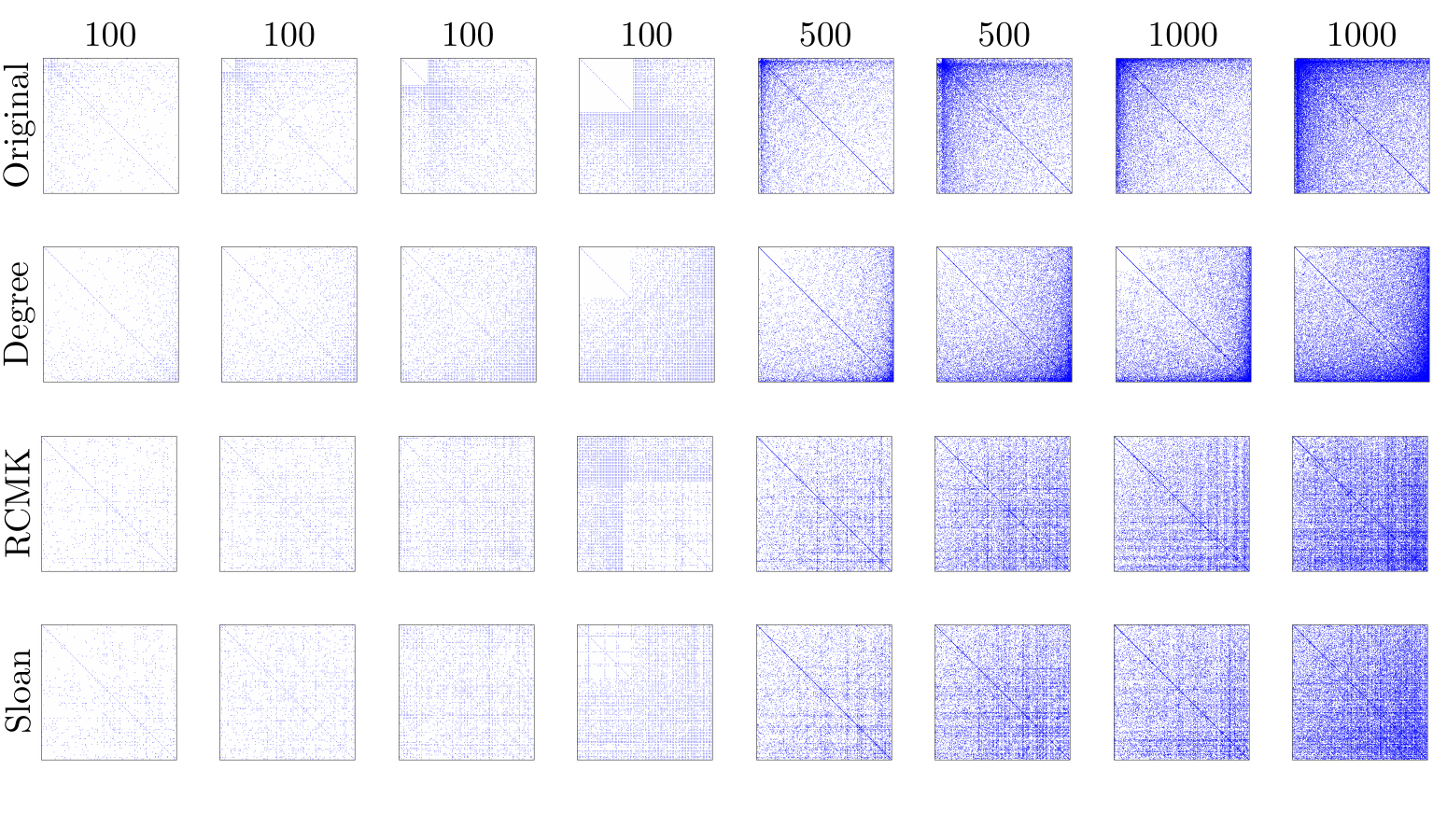}
  \caption{$\laplacian_0$ for \BA~graphs.}
  \label{fig:balplcln0rrdrng}
\end{figure}

\begin{figure}[p]
  \centering
  \begin{tabular}{ccc}
    \imagetop{\begin{xyoverpic}
        {(1, 1)}{scale=0.5}{spylplcn/ba/lplcn0/n1000m20zm.png},
        (0.5, 1.03)*{\text{Original}}
      \end{xyoverpic}} &
    \imagetop{\begin{xyoverpic}
        {(1, 1)}{scale=0.5}{spylplcn/ba/lplcn0/n1000m20dgrzm.png},
        (0.5, 1.03)*{\text{Vertex degree}}
      \end{xyoverpic}} &
    \imagetop{\begin{xyoverpic}
        {(1, 1)}{scale=0.5}{spylplcn/ba/lplcn0/n1000m20clrbr.png},
        (0.5, 1.05)*{\text{~}},
        (1.9, 0.98)*{10^2},
        (1.9, 0.685)*{10^1},
        (1.9, 0.37)*{10^0},
        (2.1, 0.055)*{10^{-1}}
      \end{xyoverpic}}
  \end{tabular}
  \caption{Enlarged view of the matrix corresponding to the
    last column in Figure~\ref{fig:balplcln0rrdrng}.}
  \label{fig:balplcln0zm}
\end{figure}

\clearpage

\begin{figure}[p]
  \centering
  \includegraphics[scale=1, trim=1.2in 6.5in 1.1in 0.9in, clip]
  {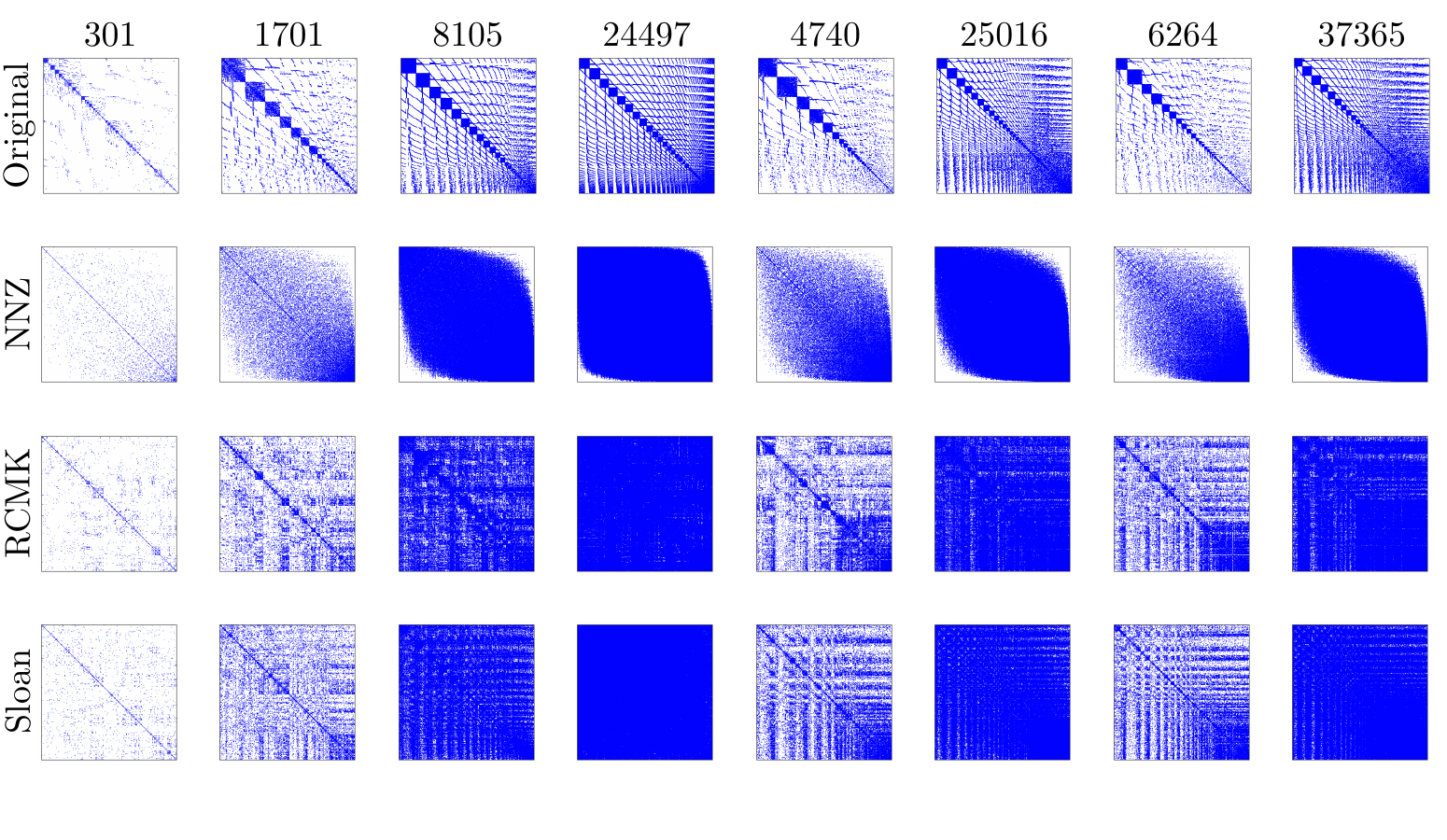}
  \caption{$\laplacian_2$ for \BA~graphs.}
  \label{fig:balplcln2rrdrng}
\end{figure}

\begin{figure}[p]
  \centering
  \begin{tabular}{ccc}
    \imagetop{\begin{xyoverpic}
        {(1, 1)}{scale=0.5}{spylplcn/ba/lplcn2/n1000m20zm.png},
        (0.5, 1.03)*{\text{Original}}
      \end{xyoverpic}} &
    \imagetop{\begin{xyoverpic}
        {(1, 1)}{scale=0.5}{spylplcn/ba/lplcn2/n1000m20nnzzm.png},
        (0.5, 1.03)*{\text{Number of nonzeros}}
      \end{xyoverpic}} &
    \imagetop{\begin{xyoverpic}
        {(1, 1)}{scale=0.5}{spylplcn/ba/lplcn2/n1000m20clrbr.png},
        (0.5, 1.05)*{\text{~}},
        (1.9, 0.98)*{10^2},
        (1.9, 0.685)*{10^1},
        (1.9, 0.37)*{10^0},
        (2.1, 0.055)*{10^{-1}}
      \end{xyoverpic}}
  \end{tabular}
  \caption{Enlarged view of the matrix corresponding to the
    last column in Figure~\ref{fig:balplcln2rrdrng}.}
  \label{fig:balplcln2zm}
\end{figure}

\clearpage

\end{document}